\documentclass[twoside,11pt]{amsbook}
\usepackage{amsmath,amsthm,amscd,amssymb,amsfonts}
\usepackage[all]{xy}
\makeindex
\UseComputerModernTips
\theoremstyle{definition}
\theoremstyle{plain}
	\newtheorem{theorem}{Theorem}[chapter]
	\newtheorem{maintheorem}{Theorem}

	\newtheorem{lemma}[theorem]{Lemma}
	\newtheorem{proposition}[theorem]{Proposition}
	\newtheorem{corollary}[theorem]{Corollary}
	\newtheorem{question}[theorem]{Question}

\theoremstyle{definition}
	\newtheorem{example}[theorem]{Example}
	\newtheorem{definition}[theorem]{Definition}
	\newtheorem{remark}[theorem]{Remark}
	\newtheorem{caveat}[theorem]{Caveat}
	\newtheorem{notation}[theorem]{Notation}

\def\mapright#1{\buildrel #1 \over \longrightarrow}
\def\functor{\hspace{0.5mm}\underline{\hspace{2mm}}\hspace{0.5mm}} 
\def\ab{{\rm ab}}
\def\cy#1{\mathcal{#1}} 
\def\mf#1{{\mathfrak #1}} 
\def\I{\cy{I}}
\def\C{\mathbb{C}}
\def\Z{\mathbb{Z}}
\def\N{\mathbb{N}}
\def\R{\mathbb{R}}
\def\Q{\mathbb{Q}}

\def\dash{\mbox{--}}

\DeclareMathOperator{\hcf}{hcf}

\DeclareMathOperator{\Krulldim}{Krull-dim}
\DeclareMathOperator{\Coker}{Coker}
\DeclareMathOperator{\Ker}{Ker}
\def\Im{{\rm Im}}

\DeclareMathOperator{\Aut}{Aut}
\DeclareMathOperator{\Hom}{Hom}
\DeclareMathOperator{\End}{End}

\DeclareMathOperator{\res}{Res}
\DeclareMathOperator{\Res}{Res}

\DeclareMathOperator{\GL}{GL}
\DeclareMathOperator{\id}{id}

\DeclareMathOperator{\Link}{Lk}

\DeclareMathOperator{\Gal}{Gal}

\DeclareMathOperator{\Trace}{Trace}
\DeclareMathOperator{\trace}{Trace}

\DeclareMathOperator\Tr{{\rm Tr}}
\DeclareMathOperator{\Det}{Det}
\DeclareMathOperator{\Support}{Support}
\def\rproj{\text{{\rm Proj-}}}
\def\proj{\text{{\rm -Proj}}}

\def\mod{\text{{\rm-Mod}}}
\newtheorem{Kexample}[theorem]{Knot Theory Example}

\def\sdM{\overline{\cy{M}}}
\def\ssdM{\sdM^{\raisebox{-.35ex}{$\scriptstyle{s}$}}}
\begin{document}
\title{Invariants of Boundary Link Cobordism}
\author{Desmond Sheiham}
\begin{abstract}
An $n$-dimensional $\mu$-component boundary link is a codimension $2$
embedding of spheres
\begin{equation*} 
L=\bigsqcup_{\mu}S^n \subset S^{n+2}
\end{equation*}
such that there exist $\mu$ disjoint oriented embedded
$(n+1)$-manifolds which span the components of~$L$.
An $F_\mu$-link is a boundary link together with a cobordism class of
such spanning manifolds.

The $F_\mu$-link cobordism group $C_n(F_\mu)$
is known to be trivial when $n$ is even but not finitely generated when $n$
is odd. Our main result is an algorithm to
decide whether two odd-dimensional $F_\mu$-links
represent the same cobordism class in~$C_{2q-1}(F_\mu)$ assuming
$q>1$. We proceed to compute the isomorphism class
of~$C_{2q-1}(F_\mu)$, generalizing Levine's computation of the knot
cobordism group $C_{2q-1}(F_1)$.

Our starting point is the algebraic formulation of Levine, Ko and Mio
who identify $C_{2q-1}(F_\mu)$ with a surgery obstruction group, the 
Witt group $G^{(-1)^q,\mu}(\Z)$ of $\mu$-component Seifert
matrices. We obtain a complete set of torsion-free invariants by passing from 
integer coefficients to complex coefficients and by applying the
algebraic machinery of Quebbemann, Scharlau and Schulte. Signatures
correspond to `algebraically integral' simple self-dual
representations of a certain quiver (directed graph with loops). These 
representations, in turn, correspond to algebraic integers on an infinite
disjoint union of real affine varieties.

To distinguish torsion classes, we consider rational coefficients in
 place of complex coefficients, expressing $G^{(-1)^q,\mu}(\Q)$ as
an infinite direct sum of Witt groups of finite-dimensional division
$\Q$-algebras with involution. The Witt group of every such algebra
appears as a summand infinitely often.

The theory of symmetric and hermitian forms
over these division algebras is well-developed. There are five classes
of algebras to be considered; complete Witt invariants are available
for four classes, those for which the local-global principle applies.
An algebra in the fifth class, namely a quaternion algebra with non-standard
involution, requires an additional Witt invariant which is defined if all
the local invariants vanish.
\end{abstract}
\maketitle
\newpage
\thanks{
I am very grateful to my PhD adviser, Andrew
Ranicki, for the energy, ideas, encouragement and
copious enthusiasm he has shared with me throughout my time as a
graduate student in Edinburgh. I would also like to thank Michael Farber for
valuable suggestions at an early stage of this thesis project.

Among many others with whom I have enjoyed interesting and useful
mathematical conversations, I would particularly like to mention
Richard Hill, David
Lewis, Kent Orr and Peter Teichner who made it possible for me to
spend time at their respective universities and were
generous to me with their time and expertise.
}
\tableofcontents
\chapter{Introduction}
\label{chapter:introduction}
The classification problem which we address in this volume concerns
$n$-dimensional spheres $S^n$ knotted and linked inside an $(n+2)$-dimensional
space. This
$(n+2)$-dimensional space could be~$\R^{n+2}$, but following standard
convention let us assume that it is the one point compactification,~$S^{n+2}$.

We focus on boundary links, links whose
components are the boundaries of {\it disjoint} $(n+1)$-dimensional
manifolds inside $S^{n+2}$. In fact, our basic object of study is the
$F_\mu$-link, a refinement of the ($\mu$-component) boundary link 
whose definition involves the free group $F_\mu$ on $\mu$
non-commuting generators; see definition~\ref{define_Flink} below.

Whilst all three notions - `link', `boundary link' and
`$F_\mu$-link' - are generalizations of `knot',
the theory of $F_\mu$-links is most directly analogous to knot theory.

The aim of the present work is to provide the means to calculate
whether two arbitrary $F_\mu$-links are `the same' or `different', up
to an equivalence relation known as cobordism\footnote{Some authors prefer the
synonym `concordance' to `cobordism'.}~\cite{FoxMil66}. 
In the definition of cobordism one regards the ambient
sphere~$S^{n+2}$ as the boundary of an $(n+3)$-dimensional disk, $D^{n+3}$.
To say that a knot, for
example, is cobordant to the trivial knot is {\it not} to say that one can
untie the knot in $S^{n+2}$, but that one can untie it in
$D^{n+3}$ by contracting it concentrically through some
$(n+1)$-dimensional disk.

If $n\geq2$ then the set $C_n(F_\mu)$ of cobordism classes of
$F_\mu$-links is an abelian group; one adds two links by
`ambient connected sum', joining corresponding components of the two
links with narrow tubes. The secondary aim of this work is to compute the
isomorphism class of $C_n(F_\mu)$. 

A detailed computation of the knot cobordism group $C_n(F_1)$ was achieved
in the 1960's and 1970's for all $n\geq2$; only the case $n=1$
remains an open problem. The even-dimensional groups
turned out to be trivial~\cite{Ker65}, whereas the
odd-dimensional groups are not even finitely generated. 

The odd-dimensional computation emerged in two stages. Firstly,
surgery methods were used to identify $C_n(F_1)$ with a 
group defined in purely algebraic
terms~\cite{Lev69,CapSha74,Kea75',Kea75,Par76,Par77,Ran81,Smi81}.  
Secondly, numerical invariants powerful enough to
distinguish all the odd-dimensional knot cobordism classes were
defined and the cobordism group was computed~\cite{Mil69,Lev69B,Ker71,Sto77}.

The first stage of the $F_\mu$-link cobordism computation,
reduction to algebra, has also been achieved.
The even-dimensional groups
are trivial while the odd-dimensional groups are `even larger' than in
knot theory. 
This first chapter is a short exposition of the reduction to algebra
of knot cobordism and $F_\mu$-link cobordism.
We recall the different notions of link cobordism and state
three (equivalent) identifications of $C_n(F_\mu)$ with surgery
obstruction groups.

Our main results, which we state in chapter~\ref{chapter:main_results},
concern the second stage of the $F_\mu$-link cobordism problem. We
define an algorithmic procedure to decide whether or not two $F_\mu$-links
are cobordant (assuming $q>1$)
and proceed to compute the
isomorphism class of $C_{2q-1}(F_\mu)$. 
The definitions of our
invariants also apply when $q=1$ but they are far from a
complete set in this case.

Previous odd-dimensional link cobordism signatures include those of
S.Cappell and J.Shaneson~\cite[p46]{CapSha80} and of
J.Levine~\cite{Lev92,Lev94} who obtained some $F_\mu$-link
signatures via jumps in the $\eta$-invariants associated to unitary
representations $F_\mu\to U(m)$ of the free group.

We employ quite different methods taking as a starting point the
algebraic formulation of $C_{2q-1}(F_\mu)$ in terms of `Seifert
matrices' (sections~\ref{section:introduce_seifert_form}
and~\ref{section:Flink_seifert_form} below).

\section{Isotopy and Cobordism; the Link}
Let us first explain what is meant by
cobordism of links.\index{Link}
A link is an embedding of disjoint $n$-dimensional spheres\footnote{
\label{smooth_PL_TOP_footnote}
More precisely, a link will be an embedding of disjoint smooth manifolds
$\Sigma_1\sqcup\cdots\sqcup\Sigma_\mu\hookrightarrow S^{n+2}$ such
that each component $\Sigma_i$ is homeomorphic to $S^n$. Alternatively, one can
interpret everything in the piecewise linear category,
if one assumes that links, isotopies and cobordisms are locally flat
PL embeddings. If one were to consider locally flat topological
embeddings the theory could differ only when $n=3$ or $4$.
(cf~Novikov~\cite[Theorem 6]{Nov68} and Cappell and
Shaneson~\cite{CapSha73}).} 
in an $(n+2)$-dimensional sphere:
\vspace{-2mm}
\begin{equation*}
L~\cong~\overbrace{S^n\sqcup\cdots \sqcup S^n}^{\mu} ~\subset~ S^{n+2}~.
\end{equation*}
Each component of a link may be knotted; indeed, a $1$-component link is
called a {\it knot}.\index{Knot} 

Two links are called isotopic and are usually
considered to be `the same' if one of the links can be transformed
into the other through embeddings:
\begin{definition}\index{Isotopy}
Links $L^0$ and $L^1$ are {\it isotopic} if they can be joined\footnote{
More precisely it is required that there exists a smooth
oriented submanifold $LI$ of $S^{n+2}\times I$ such that $LI$ is homeomorphic
to $(S^n\sqcup\cdots\sqcup S^n)\times [0,1]$ and meets
$S^{n+2}\times\{0\}$ and 
$S^{n+2}\times\{1\}$ transversely at $L^0$ and $L^1$ respectively.} 
in $S^{n+2}\times [0,1]$ by an embedding  
\begin{align*}
LI\cong(S^n\sqcup\cdots \sqcup S^n)  \times [0,1]
~&\subset~S^{n+2}\times [0,1] \\
\intertext{such that, for each $i\in[0,1]$,} 
(S^n\sqcup\cdots \sqcup S^n)\times \{i\}&\subset S^{n+2}\times \{i\}.
\end{align*}
\end{definition}
\noindent A link $L$ is in fact isotopic to the trivial link if and
only if $L$ is the boundary of $\mu$ disjoint disks embedded in
$S^{n+2}$~(see for example Hirsch~\cite[Thm8.1.5]{Hir76}).

Cobordism is a weaker equivalence relation first defined
in the context of classical knots $S^1\subset S^3$ by Fox and
Milnor~\cite{FoxMil66}. One merely omits the condition that the
embedding should be `level-preserving':
\begin{definition}\index{Cobordism!of links}\index{$C(n,\mu)$}\index{Link!cobordism}\index{Slice}\index{Null-cobordant}
\label{define_cobordism}
Links $L^0$ and $L^1$ are {\it cobordant} if they may be joined\footnotemark[\value{footnote}]
by an embedding 
\begin{equation*}
LI\cong(S^n\sqcup\cdots \sqcup S^n) \times [0,1] \subset S^{n+2}\times [0,1]~.
\end{equation*}
The set of $n$-dimensional $\mu$-component link
cobordism classes is denoted $C(n,\mu)$.
\end{definition}
A cobordism between $L$ and the trivial link can be `capped' at the
trivial end by $\mu$ disjoint disks $D^{n+1}\subset S^{n+2}$.   
Thus a link is cobordant to the trivial link if and only if it
bounds~$\mu$ disjoint disks embedded in $D^{n+3}$: 
\begin{equation}
\label{slice_link}
\begin{matrix}
{S^n\sqcup\cdots\sqcup S^n} & \subset & S^{n+2} \\
\cap & & \cap \\
{D^{n+1}\sqcup\cdots\sqcup D^{n+1}} & \subset & D^{n+3}
\end{matrix}
\end{equation}
A link which is cobordant to a trivial link is called {\it
null-cobordant} or {\it slice}.
  
The general problem of classifying knots and links up
to isotopy appears very difficult, although much is known about
certain restricted classes.
On the other hand a detailed classification of knot cobordism $C(n,1)$
was achieved in the 1960's and 1970's for all $n\geq2$.
For $\mu\geq 2$ the computation of $C(n,\mu)$ is at the time of
writing an important open problem. 
\section{Knot Theory}
\label{section:spot_of_knot}
Before introducing boundary links, 
let us review a little high-dimensional knot theory.
To the uninitiated reader we recommend Kervaire and 
Weber~\cite{KerWeb78} and the (low-dimensional) books of
Lickorish~\cite{Lic97} and Rolfsen~\cite{Rol76}. Further
expositions in high dimension include Farber~\cite{Far83}, Levine and
Orr~\cite{LevOrr00} and Ranicki~\cite[Introduction]{Ran98}.
\subsection{Seifert Surfaces}\index{Seifert!surface}
\label{section:seifert_existence_uniqueness}
It is well-known that every knot $S^n\subset S^{n+2}$ is the boundary
of a connected and oriented $(n+1)$-manifold 
\begin{equation}
\label{seifert_surface}
S^n=\partial V^{n+1}~\subset~V^{n+1}~\subset~S^{n+2}
\end{equation}
called a {\it Seifert surface}. Although there are (infinitely) many
possible choices for~$V^{n+1}$, a `Seifert matrix' of linking numbers associated to
the embedding of any Seifert surface in $S^{n+2}$ enable one to define
and compute
valuable isotopy invariants and
cobordism invariants of the knot.  
Indeed, M.Kervaire~\cite{Ker65} and 
J.Levine~\cite{Lev69} proved that when $n\geq2$
complete cobordism information is obtained. We elaborate in
section~\ref{section:algebraic_knot_cobordism} below. 

One can explain both the existence of Seifert surfaces and the
independence of associated invariants to one's choice of Seifert surface 
as follows:
Let~$X$ denote the exterior of a knot $K$, that is to say, the
complement in $S^{n+2}$ of an open tubular neighbourhood of $K$.\label{knot!exterior}\label{link!exterior} 
By Alexander duality $H_*(X)\cong H_*(S^1)$ for all knots $K$ 
so there is a canonical surjection
\begin{equation}
\label{canonical_knot_group_surjection}
\pi_1(X)\to \pi_1(X)^{\ab} \cong H_1(X) \cong \Z~.
\end{equation}
This surjection is induced by a map
\begin{equation*}
\theta:X\to S^1
\end{equation*}
which sends a meridian of the knot to a generator of $\pi_1(S^1)=\Z$
and is unique up to homotopy.
If $x\in S^1$ and $\theta$ is chosen judiciously then the inverse image
$\theta^{-1}(x)$ (together with a small collar) is a Seifert surface for $K$.

\subsection{The Infinite Cyclic Cover}
\label{section:infinite_cyclic_cover}
Although Seifert surfaces facilitate computation,
a knot invariant is perhaps more elegant if its definition does not
require arbitrary choices.
J.Milnor~\cite{Mil68} and R.Blanchfield~\cite{Bla57} considered
the algebraic topology of an infinite cyclic cover $\overline{X}\to
X$, the pull-back of $\theta:X\to S^1$ 
\begin{equation*}
\xymatrix{
{\overline{X}} \ar[r]\ar[d] & {\R} \ar[d] \\
X \ar[r] & S^1~.
}
\end{equation*}
The space $\overline{X}$ enjoys a free action by the
group $\Z$ of deck transformations. Poincar\'e duality
of $\overline{X}$, defined $\Z$-equivariantly, yields knot invariants
which are equivalent in sensitivity to `Seifert matrix' invariants. In
particular, when $n\geq2$ one can obtain complete knot cobordism
invariants. We shall outline
the various approaches to knot cobordism in
section~\ref{section:algebraic_knot_cobordism}.
\section{Boundary Links}
\label{section:boundary_links}
One would naturally like to extend the knot theory we have outlined to links;
the reader is referred to the recent book of J.Hillman~\cite{Hil02}
for a broad treatment of the algebraic theory.

Certainly one can obtain link invariants by means of a map
$\theta:X\to S^1$ where $X$ denotes the link exterior.
If $\theta$ is chosen suitably then the preimage of a point $x\in S^1$ 
is a connected oriented $(n+1)$-manifold which spans the link
(compare~Tristram~\cite{Tri69}). 

A link admits not only an infinite cyclic cover, but also a free
abelian cover. By Alexander duality $H_1(X)\cong \Z^\mu$ 
so there is a natural surjection
\begin{equation*}
\pi_1(X)\to \pi_1(X)^{\ab}\cong H_1(X) \cong \Z^\mu~.
\end{equation*}
which leads to further link invariants (e.g.~Hillman~\cite{Hil81} and
Sato~\cite{Sat81',Sat84}).

Closer parallels to high-dimensional knot theory become possible
if one assumes that there is a suitable homomorphism from $\pi_1(X)$
to the (non-abelian) free group $F_\mu$ on $\mu$-generators. 
Such a homomorphism exists if and only if the components of the link bound
{\it disjoint} oriented $(n+1)$-manifolds.
\begin{definition}\index{Link!boundary}\index{Boundary
link}\index{Seifert!surface}
A link $L$ is a {\it boundary link} if the components bound
disjoint connected oriented $(n+1)$-manifolds
\begin{equation*}
L=\partial(V^{n+1}_1 \sqcup \cdots \sqcup V^{n+1}_\mu)~\subset~V^{n+1}_1\sqcup
\cdots\sqcup V^{n+1}_\mu~\subset~S^{n+2}~.
\end{equation*} 
The union $V$ of these $(n+1)$-manifolds is called a {\it Seifert surface}.
\end{definition} 
\begin{proposition}[Gutierrez~\cite{Gut72}, Smythe~\cite{Smy66}] 
\label{characterize_boundary_links}
A link $L$ is a boundary link if and only if there is a group homomorphism
\begin{equation*}
\theta:\pi_1X \to F_\mu
\end{equation*}
which sends some choice of meridians $m_1,\cdots,m_\mu$ of the
components of $L$ to a distinguished basis 
$z_1,\cdots,z_\mu$ for $F_\mu$. 
\end{proposition}
\begin{definition}
\label{define_Flink}\index{F-link@$F_\mu$-link}\index{Link!F-@$F_\mu$-}
A pair $(L,\theta)$ as in theorem~\ref{characterize_boundary_links} is
called an {\it $F_\mu$-link}.
\end{definition}
In particular, the map~(\ref{canonical_knot_group_surjection}) above implies
that every knot is a boundary link, as we discussed in
section~\ref{section:seifert_existence_uniqueness}.
A knot can also be regarded as an $F_1$-link because
(\ref{canonical_knot_group_surjection}) is unique.

Meridians $m_i\in \pi_1(X)$ in
theorem~\ref{characterize_boundary_links} can be defined by choosing a
suitable path from a base-point of $X$ to a small meridinal circle
around each component of the link. 
Two surjections $\theta$ and $\theta'$ corresponding to
distinct choices of meridians are related 
by an equation $\theta=\alpha\theta'$ where $\alpha$ is an
automorphism of $F_\mu$ which sends each generator $z_i$ to some
conjugate~$g_iz_ig_i^{-1}\in F_\mu$. 

Returning briefly to knot theory,
the Thom-Pontrjagin construction implies that every
Seifert surface is $\theta^{-1}(x)$ for some
$\theta:X\to S^1$ and $x\in S^1$. 
Given two Seifert surfaces $V=\theta^{-1}(x)$ and
$V'={\theta'}^{-1}(x)$ one can find a homotopy $h:X\times [0,1] \to
S^1$ from $\theta$ to $\theta'$. If $h$ is chosen judiciously, then
$h^{-1}(x)$ is a cobordism between $V$ and $V'$ relative the
knot.\index{Cobordism!of Seifert surfaces}

Boundary link theory is similar: One constructs a
Seifert surface from a map not to a circle but to a one-point union of
$\mu$ circles
\begin{equation*}
X\to S^1\vee\cdots\vee S^1.
\end{equation*}
The inverse image of $\mu$ regular points chosen on distinct copies of
$S^1$ is a $\mu$-component Seifert surface. Conversely, given a 
Seifert surface the Thom-Pontrjagin construction gives a 
map $X\to S^1\vee\cdots \vee S^1$.
The choices of surjection~$\theta:\pi_1(X)\to F_\mu$ correspond
bijectively with the cobordism classes (rel L) of Seifert surfaces.

\subsection{Cobordisms}
Just as `link' is not the only generalization of `knot', `link
cobordism' is not the only generalization of `knot cobordism'.
\begin{definition}\index{Cobordism!of boundary links}\index{$B(n,\mu)$}\index{Slice!boundary-}
Two boundary links $L^0$ and $L^1$ are {\it boundary cobordant} if they can be
joined in $S^{n+2}\times [0,1]$ by a link cobordism $LI$ 
whose components bound disjoint oriented $(n+2)$-manifolds. If a
boundary link is boundary cobordant to the trivial link, it is 
called {\it boundary-slice}. The set of boundary cobordism classes of
boundary links is denoted $B(n,\mu)$. 
\end{definition}
\begin{definition}\index{Cobordism!of F-links@of $F_\mu$-links}\index{$C_n(F_\mu)$}
Two $F_\mu$-links $(L^0,\theta_0)$ and $(L^1,\theta_1)$ are
{\it cobordant} if they may be joined in $S^{n+2}\times [0,1]$ by a
pair $(LI, \Theta)$ where $LI$ is a link
cobordism and
\begin{equation*}
\Theta:\pi_1(S^{n+2}\times [0,1]~\backslash~LI) \to F_\mu
\end{equation*}
agrees with $\theta$ and $\theta'$ up to inner automorphism. The set
of cobordism classes of $F_\mu$-links is denoted $C_n(F_\mu)$.
\end{definition}
In the case of knot theory, $\mu=1$, every knot cobordism extends to a
cobordism between Seifert surfaces. In symbols
\begin{equation*}
C_n(F_1)=B(n,1)=C(n,1)~.
\end{equation*}

Returning to links, there are canonical forgetful maps
\begin{equation*}
C_n(F_\mu)\to B(n,\mu)\to C(n,\mu)
\end{equation*}
but they are not in general bijective. T.Cochran and K.Orr
proved~\cite{CocOrr90,CocOrr93} that if $n$ is odd and $\mu\geq2$ then
$B(n,\mu)\to C(n,\mu)$ is not surjective. It is an open problem to
determine whether this map is injective.

The map $C_n(F_\mu)\to B(n,\mu)$ is easier to understand. 
As we discussed in section~\ref{section:boundary_links}, a boundary
link has a splitting  
$\theta:\pi_1(X)\to F_\mu$ defined uniquely up to composition by
generator conjugating automorphisms of the free
group. If~$A_\mu$\index{$A_\mu$}
denotes the group of such automorphisms (modulo inner automorphisms) then
$B(n,\mu)$ is the set of $A_\mu$-orbits
\begin{equation*}
B(n,\mu)= C_n(F_\mu) / A_\mu~.
\end{equation*}
Note that $A_\mu$ is a non-trivial group for all $\mu\geq3$.

Ko proved further~\cite[Theorem 2.7]{Ko87} that there is a bijective
correspondence between cobordism classes of $F_\mu$-links and cobordism classes
of pairs\index{Cobordism!of Seifert surfaces}
\begin{equation*}
(V^{n+1}\subset S^{n+2},~L=\partial V)~=~(\mbox{Seifert
Surface},~\mbox{Boundary Link}).
\end{equation*} 
He also gave a geometric interpretation of the action of $A_\mu$ in terms
of Seifert surfaces.

The present work uses Seifert surface methods to define 
invariants which distinguish cobordism classes of $F_\mu$-links.
\subsection{Addition of
$F_\mu$-Links}\index{Addition of knots and $F_\mu$-links}
It is well-known that one can add two knots by performing
an ambient connected sum. Picturing the two knots embedded, but
separated, in a single ambient space $S^{n+2}$, one chooses an arc in
$S^{n+2}$ joining the two knots. 
One cuts out a small disk $D^n$ from each knot $S^n$ and attaches the
knots by a narrow tube $S^{n-1}\times[0,1]$ which surrounds the arc.
The isotopy class of this sum is independent of the
choice of connecting arc and the set of isotopy classes of knots
$S^n\subset S^{n+2}$ becomes a commutative semigroup. 

Knot cobordism respects addition of
knots; in fact $C(n,1)$ is an abelian group. The inverse of a
knot $S^n\subset S^{n+2}$ is obtained by reversing the orientations
both of~$S^n$ and of the ambient space~$S^{n+2}$.

Addition of links is more delicate. Given links $L^0$ and
$L^1$ one can join each component of $L^0$ to the corresponding
component of $L^1$ by an arc, but the isotopy class and the
cobordism class of the sum usually depend on the choice of
arcs. In general, neither $C(n,\mu)$ nor $B(n,\mu)$ have an obvious
group structure.

On the other hand, there is a well-defined notion of addition for 
two $F_\mu$-links $(L^0,\theta_0)$ and $(L^1,\theta_1)$ if one assumes
that $\theta_0$ and $\theta_1$ are isomorphisms.
Denoting $z_1,\cdots,z_\mu$ a distinguished basis of $F_\mu$,
one can use preferred meridians $\theta_0^{-1}(z_j)$ and
$\theta_1^{-1}(z_j)$ to define a connecting arc between the $j$th
component of $L^0$ and the $j$th component of $L^1$ 
(cf~Le~Dimet~\cite{LeD88}).

Equivalently, one may choose simply-connected
Seifert surfaces corresponding to $\theta_0$ and
$\theta_1$ and then join the components of $L^0$ to components of $L^1$ along
arcs which avoid both Seifert surfaces.   

In the classical dimension, $C_1(F_\mu)$ does not have an obvious
group structure if $\mu\geq2$. However, if $n\geq2$ then every
$F_\mu$-link $(L,\theta)$ is cobordant to some $(L',\theta')$ such
that $\theta'$ is an isomorphism. Moreover, cobordism respects addition of
$F_\mu$-links so $C_n(F_\mu)$ is an abelian group for all
$n\geq2$. For further details see Ko~\cite[Prop 2.11]{Ko87} or
Mio~\cite[p260]{Mio87}. 
\subsection{Split $F_\mu$-links}\index{Link!split|(}
Given $\mu$ knots there is a canonical way to make a link; one puts
all the knots in a single ambient space $S^{n+2}$, far apart from each
other so that the components lie in disjoint disks $D^{n+2}$. 
A link which can be constructed in this way is called a {\it
split} link. Choosing a Seifert surface for each knot, one obtains a
Seifert surface for the split link. It follows that a split link
is, canonically, an $F_\mu$-link $(L,\theta)$ where
$\theta:\pi_1(X)\to F_\mu$ is the free product of the canonical
homomorphisms~(\ref{canonical_knot_group_surjection}) associated to
the component knots. We thus have a canonical map
\begin{equation}
\label{include_split_links}
\prod_{i=1}^\mu C_n(F_1) \to C_n(F_\mu)
\end{equation}
which is a group homomorphism if $n\geq2$. 

Conversely, there are $\mu$ forgetful maps $C_n(F_\mu)\to C_n(F_1)$,
which, together, split~(\ref{include_split_links}). If $n\geq2$ there
is therefore a decomposition
\begin{equation*}
C_n(F_\mu)\cong \left(\bigoplus_{i=1}^\mu C_n(F_1)\right) \oplus
\widetilde{C}_n(F_\mu)~.
\end{equation*}\index{Link!split|)}

We shall return to the $F_\mu$-link cobordism group in
section~\ref{section:algebraic_Flink_cobordism}, for we must first describe
the reduction to algebra of high-dimensional knot cobordism~$C_n(F_1)$.
\section{Surgery Obstructions}\index{Surgery|(}
\label{section:algebraic_knot_cobordism}
Given a high-dimensional knot, one can attempt to construct a
cobordism to the trivial knot by performing surgeries. 
It turns out that, in even dimensions, such a feat is always possible;
every knot is null-cobordant.

In odd dimensions, obstructions exist but the attempt leads to an
identification of $C_{2q-1}(F_1)$ with a `surgery obstruction group'.

\subsection{Surgery on a Seifert Surface}\index{Surgery!on a Seifert surface|(}
\label{surgery_on_seifert_surface}
\begin{theorem}[Kervaire 1965]
\label{even_dimensional_slice}
Every even-dimensional knot $K\cong S^{2q}\subset S^{2q+2}$ is null-cobordant:
\begin{equation*}
C_{2q}(F_1)=0~\mbox{for all $q\geq1$}.
\end{equation*}
\end{theorem}
\begin{proof}[Sketch Proof] 
Starting with any Seifert surface 
\begin{equation*}
V^{2q+1}\subset S^{2q+2},\quad
\partial V=K
\end{equation*}
one performs surgery on $V^{2q+1}$ killing homology
classes of degree at most~$q$, until one has turned $V$ into a
disk~$D^{2q+1}$. One must check that the surgery operations can be performed
ambiently, not in $S^{2q+2}$, but in a disk $D^{2q+3}$ whose boundary
is $S^{2q+2}$. 

The effect of each surgery operation is to cut out a copy of
$S^i\times D^{2q-i}$ in 
the interior of $V^{2q}$ and graft a copy of $D^{i+1}\times
S^{2q-i-1}$ along the boundary $S^i\times S^{2q-i-1}$. At core, such
surgeries are possible because the embedded spheres $S^i$ bound
disks $D^{i+1}\subset D^{2q+3}$ which are disjoint and intersect $V$
only at $S^i$. 
\end{proof}
\vspace{-.2cm}
The same method shows that every even-dimensional
boundary link is boundary-slice\index{Slice!boundary-}, a result which was proved by 
Cappell and Shaneson~\cite{CapSha80} using homology surgery (see
sections~\ref{section:knot_homology_surgery}
and~\ref{section:Flink_homology_surgery}). In symbols,
\begin{equation*}
C_{2q}(F_\mu)=B(2q,\mu)=0~\mbox{for all $\mu\geq1$ and all $q\geq1$}.
\end{equation*}

By contrast not all odd-dimensional knots $K=S^{2q-1}\subset S^{2q+1}$
are slice. In fact the odd-dimensional knot cobordism groups are
not even finitely generated. If $V^{2q}$ is a Seifert surface for $K$
one can perform ambient surgery to kill all homology in degree strictly
less than~$q$, but $q$-spheres, if they are linked in $S^{2q+1}$, do
not bound disjoint disks $D^{q+1}$ in $D^{2q+1}$. The Seifert
form measures this obstruction to surgery in dimension $q$ and was used
by J.Levine~\cite{Lev69} to obtain an algebraic description of
odd-dimensional knot cobordism groups
(theorem~\ref{knot_cobordism:seifert} below).

\subsection{The Seifert Form}\index{Seifert!form}
\label{section:introduce_seifert_form}
By Alexander duality there is a
non-singular pairing which measures linking between 
$V^{2q}$ and $S^{2q+1}\backslash V^{2q}$:
\begin{align*}
H_q(V)\times H_q(S^{2q+1}\backslash V) &\to \Z \\
(x,x') &\mapsto \Link(x,x')~.
\end{align*}

We may assume, after preliminary surgeries below the middle dimension,
that $H_q(V^{2q})$ is a free abelian group. 
Let $x_1,\cdots,x_m$ be a basis for $H_q(V)$ and let
$i_+,i_-:H_q(V^{2q})\to H_q(S^{2q+1}\backslash V^{2q})$ be 
small translations in the positive and negative normal directions
to~$V^{2q}$. These directions are determined by the orientation of $V$. 

\begin{definition}\index{Seifert!form}
The {\it Seifert matrix} $S$ associated to the Seifert surface $V$
with respect to $x_1,\cdots, x_m$ is the matrix
\begin{equation*}
S_{ij}=\Link(i_+x_i,x_j)=\Link(x_i,i_-x_j)~.
\end{equation*}
\end{definition}
\noindent Addition of knots $K+K'$ by
connected sum induces a block sum 
$\left(\begin{smallmatrix}
S & 0 \\
0 & S'
\end{smallmatrix}\right)$
of Seifert matrices.
Although the Seifert matrix $S$ is not symmetric or skew-symmetric,
$S+(-1)^qS^t$ is certainly $(-1)^q$-symmetric (where $S^t$ denotes the
transpose of~$S$). It is not difficult to see that 
\begin{equation*}
S_{ij}+(-1)^qS_{ji}=\Link((i_+-i_-)x_i,x_j)
\end{equation*}
which is the intersection pairing
on~$H_q(V^{2q})$. By Poincar\'e duality this intersection matrix
$S+(-1)^qS^t$ is non-singular\footnote{The exact 
sequence for the pair $(V,\partial V)$ implies that $H_q(V,\partial
V)\cong H_q(V)$.}.

Kervaire proved~\cite[Th\'eor\`eme II.3]{Ker65} that every matrix $S$ such that
$S+(-1)^q S^t$ is invertible is the Seifert matrix of some Seifert
surface of some knot $S^{2q-1}\subset S^{2q+1}$ (when
$q\neq2$). Note that if one working in the category of smooth
manifolds the differentiable structure on $S^{2q-1}$ may be exotic
(see footnote~\ref{smooth_PL_TOP_footnote} on
page~\pageref{smooth_PL_TOP_footnote}).
To outline Kervaire's argument, one can first construct a stably
parallelizable manifold $V^{2q}$ with intersection pairing $S+(-1)^q
S^t$ by attaching $q$-handles $D^q\times D^q$ to the boundary of a
zero-handle $D^{2q}$. Since 
$S+(-1)^q S^t$ is invertible, the boundary of $V^{2q}$ is a homotopy
sphere, and hence is homeomorphic to a sphere. One can
complete the construction by adjusting an embedding
$V^{2q}\hookrightarrow S^{2q+1}$ until the handles are linked in the manner
dictated by $S$.

If a $(2q-1)$-dimensional knot $K$ is slice then one can find a basis
for $H_q(V)$ half of which 
is completely unlinked in $S^{2q+1}$. In other words, the Seifert
matrix has the appearance
\begin{equation}
\left(\begin{matrix}
0 & * \\
* & *
\end{matrix}\right)
\end{equation}
where each $*$ denotes some square matrix with half the rank
of $H_q(V)$. Such a Seifert matrix is called {\it
metabolic}.\index{Metabolic} It follows that there is a group
homomorphism from $C_{2q-1}(F_1)$ to the Witt
group\index{$G^{\epsilon,\mu}(A)$}\index{Witt group!of Seifert forms}
$G^{(-1)^q,1}(\Z)$ of Seifert matrices modulo metabolic matrices (see
section~\ref{section:define_seifert_form} of
chapter~\ref{chapter:preliminaries} for a more precise definition). 

If $q\geq2$ the converse is true: Any knot $S^{2q-1}\subset
S^{2q+1}$ which has a metabolic Seifert form is a slice knot.\index{Slice}
\begin{theorem}[Levine 1969]
\label{knot_cobordism:seifert}
If $q\geq3$ then $C_{2q-1}(F_1)\cong G^{(-1)^q,1}(\Z)$. 
\end{theorem} 
In the case $q=2$, the cobordism group $C_3(F_1)$ is isomorphic to an
index two subgroup of $G^{1,1}(\Z)$. If $q=1$ there is a surjection
$C_1(F_1)\twoheadrightarrow G^{-1,1}(\Z)$ but A.Casson and C.Gordon
defined knot invariants~\cite{CasGor78,CasGor86} which show that the kernel is
non-trivial (see also Gilmer~\cite{Gil83,Gil93}, Kirk and
Livingston~\cite{KirLiv99} and Letsche~\cite{Let00}).   
T.Cochran, K.Orr and P.Teichner~\cite{COT99} have defined an infinite 
tower of obstructions to slicing a knot and have recovered all the earlier
cobordism invariants in the early layers. At the time of writing, it
remains an open problem to complete the cobordism classification of classical
knots.\index{Surgery!on a Seifert surface|)}
 
\subsection{The Blanchfield Form}\index{Blanchfield form|(}
As we indicated in section~\ref{section:infinite_cyclic_cover} 
the infinite cyclic cover $\overline X^{2q+1}$ of an odd-dimensional 
knot exterior $X^{2q+1}$ can be used to define an algebraic
obstruction to slicing a knot without reference to a choice of Seifert
surface.
 
Since $\Z$ acts on $\overline{X}$, the homology  
$H_i(\overline X)$ is a module over the group ring
$\Z[\Z]=\Z[z,z^{-1}]$.
We may assume, if necessary by taking a quotient by the $\Z$-torsion
subgroup, that $H_q(\overline{X})$ is $\Z$-torsion-free and has a
one-dimensional presentation
\begin{equation}
\label{knot_module_presentation}
0\to \Z[z,z^{-1}]^m\xrightarrow{\sigma} \Z[z,z^{-1}]^m \to
H_q{\overline X} \to 0.
\end{equation} 
The determinant $p=\det(\sigma)\in \Z[z,z^{-1}]$ of such a presentation
is called the {\it Alexander polynomial}\index{Alexander polynomial} of the knot and is well-defined up
to multiplication by units $\pm z^i$. Since $H_*(X;\Z)\cong
H_*(S^1;\Z)$, one finds that $\sigma$ becomes an invertible matrix
over $\Z$ if one sets $z=1$. In other words $p(1)=\pm1$. 

By Poincar\'e duality\footnote{
Note that $H_q(\overline{X},\partial \overline{X})\cong H_q(\overline{X})$.
} 
there is a non-singular `Blanchfield linking form': 
\begin{equation}
\label{blanchfield_form} 
\beta:H_q(\overline X)\times H_q(\overline X) \to 
\frac{P^{-1}\Z[z,z^{-1}]}{\Z[z,z^{-1}]}
\end{equation}  
where $P=\{p\in \Z[z,z^{-1}]\mid p(1)=\pm1\}$ is the set of
Alexander polynomials. 
The ring $\Z[z,z^{-1}]$ has an involution $\overline{z}=z^{-1}$ with
respect to which the Blanchfield form is $(-1)^{q+1}$-hermitian:
\begin{equation*} 
\beta(x,x')=(-1)^{q+1}\overline{\beta(x',x)}~\mbox{for all $x,x'\in H_q(\overline{X})$.}
\end{equation*}

C.Kearton showed~\cite{Kea75} that a Seifert form for a knot is
metabolic if and only if the Blanchfield form of the knot is
metabolic. He also proved~\cite[Theorem 8.3]{Kea75'} that if $q\neq2$ then
every $\Z[z,z^{-1}]$-module with
presentation~(\ref{knot_module_presentation}) and a linking
pairing~(\ref{blanchfield_form}) is the Blanchfield form of some knot
$S^{2q-1}\subset S^{2q+1}$. Consequently,
\begin{theorem}[Kearton, Levine 1975]
\label{knot_cobordism_blanchfield}
If $q\geq3$ then $C_{2q-1}(F_1)$ is isomorphic to the Witt group of
$(-1)^{q+1}$-hermitian Blanchfield forms.\index{Witt group!of
Blanchfield forms}
\end{theorem}\index{Blanchfield form|)}
\subsection{Surgery and Homology Surgery}\index{Surgery!homology}
\label{section:knot_homology_surgery}
A third approach to knot cobordism was pioneered by S.Cappell and 
J.Shaneson~\cite{CapSha74}. We shall give a simplified outline of their
method, but we must first review a little of the original surgery theory of
Browder, Novikov, Sullivan, Wall and many others.\index{Surgery|textbf}

One of the central questions in surgery theory is the following
\begin{question}
\label{usual_surgery_question}
Given a homotopy equivalence $f:N^n\to M^n$ between two manifolds, is
$f$ is homotopic to a homeomorphism? 
\end{question}

One can attack the problem in two stages. Firstly one seeks a (normal) bordism
$F:W^{n+1}\to M\times[0,1]$ such that $\partial W=N\sqcup N'$ and such that 
$F$ restricts to $f:N\to M\times\{0\}$ and to a homeomorphism $N'\to
M\times\{1\}$. Secondly, if such a bordism
$F$ exists, one tries to change it, by surgeries, to make it a
(simple) homotopy equivalence. In other words, one asks whether $F$ is
bordant, relative the boundary, to a homotopy equivalence. 
If these two steps can both be achieved then the $s$-cobordism theorem
implies that $f$ is homotopic to a homeomorphism. 

Without dwelling on the potential obstruction to the first stage, let
us note that there is an obstruction $\sigma(F)$ to the second stage
of the program which takes values in a group $L_{n+1}(\Z[\pi_1(M)])$\index{$L_n(A)$}
introduced by C.T.C.Wall~\cite{Wall70}. The map $f$ is
homotopic to a homeomorphism if and only if one can find a bordism $F$ 
such that $\sigma(F)=0\in L_{n+1}(\Z[\pi_1(M)])$.

In A.Ranicki's algebraic surgery theory~\cite{Ran80,Ran80',Ran81},
$L_{n+1}(\Z[\pi])$ is the cobordism group of
$(n+1)$-dimensional free $\Z[\pi]$-module chain complexes with (quadratic)
Poincar\'e duality structure. 
Even-dimensional $L$-groups are isomorphic to Witt groups of
quadratic forms whereas the odd-dimensional groups can
be expressed in terms of `formations'.

Returning to knot theory, the question whether a knot is null-cobordant has
a formulation analogous to~\ref{usual_surgery_question}. 
The exterior $X^{n+2}$ of a knot $K$ is homotopic to the trivial knot
exterior~$X^{n+2}_0$ if and only if $K$ is trivial~(Levine~\cite{Lev65},
Wall~\cite[Theorem 16.4]{Wall70}). 
On the other hand, {\it any} knot exterior~$X^{n+2}$ admits
a homology equivalence\footnote{
To be more precise, one begins with a degree one
normal map  whose restriction to boundaries
is a homotopy equivalence~$\partial X\simeq \partial X_0 \simeq
S^n\times S^1$.
}  
$\theta:X\to X_0$ which restricts to a homotopy equivalence $\partial
X\to \partial X_0$.

Conversely, suppose we are given a manifold $X^{n+2}$ with the properties that
$\partial X=S^1\times S^n$ and $\pi_1(X)$ is
normally generated by $\pi_1(\partial X)=\Z$. If $H_*(X)=H_*(X_0)$
then one can regard $X^{n+2}$ as a knot exterior by gluing a solid
torus $D^2\times S^n$ along the boundary to form a sphere $S^{n+2}$,
in which the core of the torus $\{0\}\times S^n$ is a knot. The
exterior of a knot cobordism has a similar characterization so the
question whether a knot $K$ is null-cobordant translates as follows:
\begin{question}
Given a homology equivalence $\theta:X\to X_0$, is $\theta$ homology bordant
to a homotopy equivalence?
\end{question}
A homology bordism between $\theta:X\to X_0$ and $\theta':X'\to X_0$
is by definition a homology equivalence $\Theta:W\to X_0$ such that
$\partial W=X\sqcup X'$, 
and such that $\Theta|_X=\theta$ and $\Theta|_{X'}=\theta'$.

One can tackle the problem in two stages.
Firstly one seeks any (normal) bordism $\Theta:W\to X_0$ between
$\theta:X\to X_0$ and a homotopy equivalence 
$\theta':X'\to X_0$. Secondly one tries to change the
bordism $\Theta$, by surgeries,
to make it a homology equivalence. There are obstructions at each
stage; the first takes 
values in 
\begin{equation*}
L_{n+2}(\Z[\pi_1(X_0)])=L_{n+2}(\Z[z,z^{-1}])
\end{equation*}
while the second lies in a group 
\begin{equation*}
\Gamma_{n+3}(\epsilon:\Z[z,z^{-1}]\to\Z)
\end{equation*}
defined by Cappell and Shaneson.
Here, $\epsilon$ is the augmentation map which evaluates polynomials at $z=1$. 

In the framework of Ranicki~\cite[\S2.4]{Ran81}, $\Gamma_n(\Z[z,z^{-1}]\to
\Z)$\index{$Gamma_n$@$\Gamma_n(\Z[\pi]\to\Z)$} is the cobordism group of free
$n$-dimensional $\Z[z,z^{-1}]$-module chain complexes with a quadratic
structure such that the induced structure over $\Z$ is Poincar\'e. 

The knot $K$ is null-cobordant if and only if one can find a bordism
$\Theta$ between $\theta$ and a homotopy equivalence such that the
surgery obstruction
$\sigma(\Theta)\in\Gamma_{n+3}(\epsilon:\Z[z,z^{-1}]\to\Z)$ vanishes.
These obstructions, together with an appropriate realization theorem,
imply that the various cobordism groups fit together in an exact sequence:
\begin{equation}
\label{knot_homology_surgery_sequence}
\begin{split}
\cdots\to L_{n+3}(\Z[z,z^{-1}]) \to
\Gamma_{n+3}&\left(\Z[z,z^{-1}]\xrightarrow{\epsilon} \Z\right) \\
&\to C_n(F_1) \to L_{n+2}(\Z[z,z^{-1}]) \to\cdots 
\end{split}
\end{equation}
Identifying $L_*(\Z[z,z^{-1}])\cong
\Gamma_*(\Z[z,z^{-1}]\to\Z[z,z^{-1}])$ one obtains:
\begin{theorem}[Cappell and Shaneson 1974]
\label{knot_homology_surgery}
If $n\neq1$ or $3$ then $C_n(F_1)$ is isomorphic to a relative $\Gamma$-group:
\begin{equation*}
C_n(F_1)~\cong~
\Gamma_{n+3}\left(\vcenter{\xymatrix{ 
{\Z[z,z^{-1}]} \ar[r]^{\displaystyle{\id}}\ar[d]_{\displaystyle{\id}}
& {\Z[z,z^{-1}]}\ar[d]^{\displaystyle{\epsilon}} \\ 
{\Z[z,z^{-1}]} \ar[r]_{\displaystyle{\epsilon}} & {\Z}
}}\right)~.
\end{equation*}
\end{theorem}
Cappell and Shaneson also gave a geometric interpretation of the
$4$-periodicity $C_n(F_1)\cong C_{n+4}(F_1)$, $n\geq4$ of
high-dimensional knot cobordism groups~\cite[Theorem 13.3]{CapSha74}.

The relationship between theorems~\ref{knot_homology_surgery}
and~\ref{knot_cobordism_blanchfield} was explained
by Pardon~\cite{Par76,Par77}, Ranicki~\cite[\S7.9]{Ran81} and
Smith~\cite{Smi81}.
Recalling that $P=\epsilon^{-1}\{\pm1\}\subset \Z[z,z^{-1}]$ is the set of
Alexander polynomials, the Witt group of Blanchfield forms in
the\nobreak orem~\ref{knot_cobordism_blanchfield} is isomorphic to a relative 
$L$-group:
\begin{equation*}
C_n(F_1)\cong L_{n+3}(\Z[z,z^{-1}],P)~\mbox{for all $n\neq1$ or $3$}
\end{equation*} 
and fits into a localization exact sequence
\begin{equation}
\label{knot_exact_L_sequence}\index{Localization exact sequence}
\begin{split} 
\cdots \to L_{n+3}(\Z[z,z^{-1}]) \to L_{n+3}(&P^{-1}\Z[z,z^{-1}]) \\ 
&\to C_n(F_1)\to L_{n+2}(\Z[z,z^{-1}])\to \cdots
\end{split}
\end{equation}
This sequence is isomorphic
to~(\ref{knot_homology_surgery_sequence}) 
for there is a natural isomorphism
\begin{equation}
\label{localization_is_Gamma}
L_*(P^{-1}\Z[z,z^{-1}]) \cong \Gamma_*(\Z[z,z^{-1}]\to\Z)~.
\end{equation}
A systematic approach to the algebraic $K$- and $L$-theory which
pertains to high-dimensional knots and other manifold embeddings can
be found in Ranicki~\cite{Ran98}\footnote{To 
translate notation, the Witt group $G^{(-1)^q,1}(A)$ of Seifert forms
is denoted ${\rm LIso}^{2q}(A)$ in~\cite{Ran98}.}.
\section{Algebraic $F_\mu$-link Cobordism}
\label{section:algebraic_Flink_cobordism} 
\subsection{Homology Surgery}\index{Surgery!homology}
\label{section:Flink_homology_surgery}
The first identification of $F_\mu$-link cobordism groups with surgery
obstruction groups was obtained by S.Cappell and
J.Shaneson~\cite{CapSha80} using homology surgery methods.

The exterior $X$ of a boundary link $L$
is homotopy equivalent to the trivial link exterior $X_0$ if
and only if $L$ is isotopic to the trivial link
(Gutierrez~\cite{Gut72}). On the other hand, all link
exteriors are homology equivalent, by Alexander duality, and a boundary
link exterior admits a homology equivalence $\theta:X\to X_0$ which restricts
to a homotopy equivalence on each boundary component $S^1\times S^n$
(cf~theorem~\ref{characterize_boundary_links}).

Conversely, suppose we are given a manifold $X^{n+2}$ with boundary
$\partial X=\bigsqcup_\mu S^1\times S^n$ such that $\pi_1(X)$ is normally
generated by the fundamental groups of the boundary components.
If $\theta:X\to X_0$ is a homology equivalence which is a homotopy
equivalence on boundary components, one can recover an
$F_\mu$-link by gluing $\mu$ copies of $D^2\times S^n$ onto the
boundary of $X$.

The question whether an $F_\mu$-link is cobordant to the trivial link
translates into the question whether $\theta:X\to X_0$ is homology
bordant to a homotopy equivalence.
\begin{theorem}[Cappell and Shaneson 1980]
\label{Flink_cobordism:Hsurgery}
If $n\neq1,3$ and $\mu\geq1$ there is an isomorphism
\begin{equation*}
C_n(F_\mu)\cong
\Gamma_{n+3}\left(\vcenter{\xymatrix{ 
{\Z[F_\mu]} \ar[r]^{\displaystyle{\id}}\ar[d]_{\displaystyle{\id}}
& {\Z[F_\mu]}\ar[d]^{\displaystyle{\epsilon}} \\ 
{\Z[F_\mu]} \ar[r]_{\displaystyle{\epsilon}} & {\Z}
}}\right)
\end{equation*} 
where the augmentation $\epsilon:\Z[F_\mu]\to \Z$ sends every element
of the free group to~$1$. 
\end{theorem}  
In other words, $C_n(F_\mu)$ fits into a long exact sequence
\begin{equation}
\label{Flink_homology_surgery_sequence}\index{$Gamma_n$@$\Gamma_n(\Z[\pi]\to\Z)$} 
\begin{split}
\cdots\to L_{n+3}(\Z[F_\mu]) \to
\Gamma_{n+3}&\left(\Z[F_\mu]\xrightarrow{\epsilon} \Z\right) \\
&\to C_n(F_\mu) \to L_{n+2}(\Z[\mu]) \to\cdots 
\end{split}
\end{equation}
There is a two-stage obstruction to the existence of a
homology bordism between $\theta$ and a homotopy equivalence just as
in knot theory (section~\ref{section:knot_homology_surgery}). The first is  
the Wall obstruction
\begin{equation*}
\sigma(\theta)\in L_{n+2}(\Z[\pi_1(X_0)])~=~L_{n+2}(\Z[F_\mu])~=~\Gamma_{n+2}(\Z[F_\mu]\to\Z[F_\mu])
\end{equation*}
to the existence of any (normal) bordism $\Theta$ while the second is
the homology surgery obstruction
\begin{equation*}
\sigma(\Theta)\in \Gamma_{n+3}(\Z[F_\mu]\xrightarrow{\epsilon}\Z).
\end{equation*}
\subsection{The $\mu$-component Seifert Form}\index{Seifert!form|(}
\label{section:Flink_seifert_form}
K.H.Ko~\cite{Ko87} and W.Mio~\cite{Mio87} described the
odd-dimensional groups $C_{2q-1}(F_\mu)$ in terms
of $\mu$-component Seifert forms. Their results will be our starting
point. To illustrate the definitions let us
consider boundary links of two components for the general case is very
similar (see section~\ref{section:define_seifert_form} in chapter~\ref{chapter:preliminaries}).

The homology of a two component Seifert surface is a direct sum
$H_q(V)=H_q(V_1)\oplus H_q(V_2)$. If one chooses a basis for each
summand $H_q(V_i)$ then one obtains a Seifert matrix
\begin{equation*}
S=\left(\begin{array}{c|c}
S^{11} & S^{12} \\ \hline
S^{21} & S^{22} 
\end{array}\right)
\end{equation*}
where the entries of $S^{ij}$ are linking numbers between elements of
$H_q(V_i)$ and $H_q(V_j)$ in the sphere $S^{2q+1}$. 
The intersection pairing of the Seifert surface is given, as before, by the
$(-1)^q$-symmetric non-singular matrix
\begin{equation*}
T=\left(\begin{array}{c|c}
T_1 & 0 \\ \hline
0 & T_2  
\end{array}\right)= S+(-1)^qS^t~.
\end{equation*}

In knot theory, two Seifert matrices $S$ and $S'$ are congruent
if and only if $S=PS'P^t$ for some invertible matrix $P$. In boundary
link theory, any change of basis must preserve the component
structure. Two $2$-component Seifert matrices $S$ and $S'$ are congruent 
if and only if $S=PS'P^t$ where
$P=\left(\begin{array}{c|c}
P_1 & 0 \\ \hline
0 & P_2
\end{array}\right)$ 
and~$P_1$ and~$P_2$ are invertible.

The definition of metabolic\index{Metabolic} must be refined
accordingly. A Seifert matrix is metabolic if and only if it is
congruent to a matrix of the form
\begin{equation*}
PSP'=\left(\begin{array}{cc|cc}
0 & * & 0 & * \\
* & * & * & * \\ \hline
0 & * & 0 & * \\
* & * & * & * 
\end{array}\right)
\end{equation*}
where the four matrices in the $ij$th block each have half as many columns and
half as many rows as $S^{ij}$.
The Witt group of $\mu$-component Seifert matrices will be
denoted~$G^{(-1)^q,\mu}(\Z)$.\index{$G^{\epsilon,\mu}(A)$}

A handlebody construction shows that if $q\geq3$ then every
$\mu$-component Seifert matrix 
is realised by some
$(2q-1)$-dimensional boundary link. Levine's
theorem~\ref{knot_cobordism:seifert} generalizes as follows:
\begin{theorem}[Ko, Mio 1987]
\label{Flink_cobordism:seifert}
For all $\mu\geq1$ and for all $q\geq3$, 
\begin{equation*}
C_{2q-1}(F_\mu)\cong G^{(-1)^q,\mu}(\Z).
\end{equation*}
\end{theorem}
Ko also studied the cases $q=2$ and $q=1$
showing that $C_3(F_\mu)$ is isomorphic to an index $2^\mu$ subgroup
of $G^{1,\mu}(\Z)$ while $C_1(F_\mu)$ maps onto $G^{-1,\mu}(\Z)$.
He explored geometrically the relation between
theorems~\ref{Flink_cobordism:seifert}
and~\ref{Flink_cobordism:Hsurgery} in~\cite{Ko89}.\index{Seifert!form|)}

\subsection{The Blanchfield-Duval Form}\index{Blanchfield form}
J.Duval~\cite{Duv86} generalized to $F_\mu$-links the Blanchfield
forms of odd-dimensional knot theory.
In place of the infinite cyclic cover of a knot exterior one studies
the covering $\overline X$ of a boundary link exterior induced by
$\theta:\pi_1(X)\to F_\mu$. 
Unlike knot theory, the middle homology $H_q(\overline X)$ is such
that $\Q\otimes_\Z H_q(\overline X)$ is 
infinite-dimensional unless it is trivial. However, taking a quotient
by the torsion subgroup is necessary we may assume that $H_q(\overline
X)$ is $\Z$-torsion-free and has a presentation 
\begin{equation*}
0\to \Z[F_\mu]^m\xrightarrow{\sigma}\Z[F_\mu]^m\to H_q(\overline X)\to 0
\end{equation*}
where $\sigma$ is a square matrix which becomes invertible under the
augmentation $\epsilon:\Z[F_\mu]\to\Z$ (cf
(\ref{knot_module_presentation}) above and Sato~\cite{Sat81}). 
Let $\Sigma$ denote the set of
such $\epsilon$-invertible matrices.

Duval showed that the $F_\mu$-equivariant Poincar\'e duality of
$\overline X$ gives rise to a linking form (compare Farber~\cite{Far92B})
\begin{equation*}
H_q(\overline X) \times H_q(\overline X) \to
\frac{\Sigma^{-1}\Z[F_\mu]}{\Z[F_\mu]}
\end{equation*}
where $\Sigma^{-1}\Z[F_\mu]$ is the Cohn localization, defined by
formally adjoining inverses to all the matrices in $\Sigma$.
M.Farber and P.Vogel~\cite{FarVog92} subsequently identified
$\Sigma^{-1}\Z[F_\mu]$ with a ring of rational power series in $\mu$ non-commuting indeterminates.
\begin{theorem}[Duval 1986]
\label{Flink_cobordism:duval}
If $q\geq3$ then $C_{2q-1}(F_\mu)$ is isomorphic to the Witt group of
Blanchfield-Duval forms.\index{Blanchfield-Duval form}
\end{theorem}
It follows that $C_{2q-1}(F_\mu)$ is isomorphic
 to~$L_{2q+2}(\Z[F_\mu],\Sigma)$,
 the relative $L$-group in the localization exact
sequence of Vogel~\cite{Vog80,Vog82}:
\begin{equation*}
\label{Flink_exact_L_sequence}\index{Localization exact sequence}
\cdots \to L_{n+3}(\Z[F_\mu]) \to L_{n+3}(\Sigma^{-1}\Z[F_\mu]) 
\to C_n(F_\mu)\to L_{n+2}(\Z[F_\mu])\to \cdots
\end{equation*}
This $L$-theory exact sequence is isomorphic to the exact
sequence~(\ref{Flink_homology_surgery_sequence}) above
(see~Duval~\cite[pp633-634]{Duv86}).

By algebraic means M.Farber constructed~\cite{Far91,Far92} self-dual
finite rank lattices inside a link module $H_q(\overline X)$ to mimic the
geometric relationship between the Seifert surfaces of a boundary link
and the free cover of the link exterior.\index{Surgery|)}
\chapter{Main Results}
\label{chapter:main_results}
The title of the present volume echoes that of a paper by
J.Levine~\cite{Lev69B} which gave an algorithm to
decide whether or not two $(2q-1)$-dimensional knots are cobordant, assuming
$q>1$; we shall define an algorithm to decide whether or not two
$(2q-1)$-dimensional $F_\mu$-links are cobordant, again for $q>1$.

It is a theorem of Kervaire that all the even-dimensional groups are
trivial (see chapter~\ref{chapter:introduction},
section~\ref{surgery_on_seifert_surface}). In odd dimensions Levine
used his knot invariants to prove:
\begin{theorem}[Levine 1969]
\label{knot_answer_up_to_isomorphism}
In every odd dimension $2q-1>1$, the knot cobordism group is
isomorphic to a countable direct sum of cyclic groups of orders $2$,
$4$ and~$\infty$:
\begin{equation*}
C_{2q-1}(F_1)\cong \Z^{\oplus\infty}\oplus
\left(\frac{\Z}{2\Z}\right)^{\oplus\infty}
\oplus \left(\frac{\Z}{4\Z}\right)^{\oplus\infty}.
\end{equation*}
\end{theorem}
\noindent Note that any two decompositions of a group into direct sums of
cyclic groups of infinite and prime power order are isomorphic
(e.g.~Fuchs~\cite[Theorem 17.4]{Fuc70}).

The following theorem is a corollary of our main result,
theorem~\ref{maintheorem_invariants}:
\begin{maintheorem}
\label{the_answer_up_to_isomorphism}
If $\mu\geq2$ and $2q-1>1$ then the $F_\mu$-link cobordism group is
isomorphic to a countable direct sum of cyclic groups of orders $2$,
$4$, $8$ and~$\infty$:
\begin{equation*}
C_{2q-1}(F_\mu)\cong \Z^{\oplus\infty}\oplus
\left(\frac{\Z}{2\Z}\right)^{\oplus\infty}\oplus\left(\frac{\Z}{4\Z}\right)^{\oplus\infty}
\oplus\left(\frac{\Z}{8\Z}\right)^{\oplus\infty}.
\end{equation*}
\end{maintheorem}

Levine utilized work of J.Milnor~\cite{Mil69} `On isometries of inner
product spaces' to obtain Witt invariants of
($1$-component) Seifert forms. We shall employ more general 
non-commutative methods formulated by Quebbemann, Scharlau and
Schulte~\cite{QSS79} to define Witt-invariants of $\mu$-component
Seifert forms.  
Seifert form invariants correspond to cobordism invariants of knots and
$F_\mu$-links by theorems~\ref{knot_cobordism:seifert}
and~\ref{Flink_cobordism:seifert} above. 

One could define $F_\mu$-link cobordism invariants using
Blanchfield-Duval\index{Blanchfield-Duval form} forms
(theorem~\ref{Flink_cobordism:duval}) in place of Seifert forms. The
advantage of the Seifert form approach pursued here for
computational purposes is that the middle homology of a Seifert
surface is finitely generated as a $\Z$-module whereas the middle
homology of a free cover of a boundary link complement is not (unless
it is zero).

Although the algebraic techniques we use are a direct generalization
of Milnor's, our $F_\mu$-link invariants do not directly generalize
Levine's knot cobordism invariants; our approach is closer to that of
Kervaire~\cite{Ker71} (see also Kervaire and
Weber~\cite[p107-111]{KerWeb78} and Ranicki~\cite[Chapter 42]{Ran98}).

We shall pay particular attention to signature\index{Signature} invariants, which
detect the torsion-free part of the $F_\mu$-link cobordism group. In knot
theory, there is a signature for each algebraic integer on the
line $\left\{\frac{1}{2}+ai\mid a\in\R\right\}$, or, equivalently, one
for each point on the unit circle which is a root of an Alexander
polynomial\index{Alexander polynomial}.
We shall define an $F_\mu$-link signature
for each `algebraic integer' on a disjoint union of real algebraic varieties.

\section{Signatures}\index{Signature}
Before stating results about $F_\mu$-link signatures, let us first
describe the corresponding knot theory.

Recall from section~\ref{section:introduce_seifert_form} of
chapter~\ref{chapter:introduction} that
a Seifert matrix~$S$ of linking 
numbers can be associated to each choice of Seifert surface $V^{2q}$
for a knot $S^{2q-1}\subset S^{2q+1}$. If $M=H_q(V)$ and
$\epsilon=(-1)^q$ then $S$ represents a bilinear `Seifert form'
$\lambda:M\to M^*$ such that $\lambda + \epsilon\lambda^*$ is invertible. 
The Witt group of these Seifert forms with integer entries is denoted
$G^{\epsilon,1}(\Z)$ and Levine's theorem~\ref{knot_cobordism:seifert}
states that $C_{2q-1}(F_1)\cong G^{\epsilon,1}(\Z)$ when $q\geq3$. 

\subsection{Complex Coefficients}
A complete set of {\it torsion-free} invariants of $G^{\epsilon,1}(\Z)$ can be 
obtained by replacing the coefficient ring~$\Z$ by
$\C^-$.\index{$C^-$@$\C^-$, $\C^+$}
The superscript indicates that we work
in the hermitian setting, the involution on $\C$ being complex
conjugation. `Tensoring by $\C^-$' defines a canonical   
map $G^{\epsilon,1}(\Z)\to G^{\epsilon,1}(\C^-)$.
\begin{Kexample}
\label{knot_signatures}
The kernel of the canonical map
\begin{equation}
\label{knot_Z_to_C}
G^{\epsilon,1}(\Z)\to G^{\epsilon,1}(\C^-)
\end{equation}
is $4$-torsion. Moreover, $G^{\epsilon,1}(\C^-)$ is a free abelian
group
\begin{equation*}
G^{\epsilon,1}(\C^-)\cong \Z^{\oplus\infty};
\end{equation*}
the points on the line $\{\frac{1}{2}+bi\mid b\in\R\}$ correspond to a
basis for $G^{\epsilon,1}(\C^-)$. The image of $G^{\epsilon,1}(\Z)$ is
contained in the subgroup with one basis
element for each algebraic integer $\frac{1}{2}+bi$.
The composition of~(\ref{knot_Z_to_C}) with a projection of
$G^{\epsilon,1}(\C^-)$ onto the subgroup generated by a basis element
$\omega\in \{\frac{1}{2}+bi\mid b\in\R\}$ is called a signature
\begin{equation*}
\sigma_\omega:G^{\epsilon,1}(\Z)\to \Z~.
\end{equation*}
Note that $\sigma_\omega=\sigma_{\overline{\omega}}$.
\end{Kexample}

We turn to signatures of $F_\mu$-links.
Recall that a $\mu$-component Seifert form is a
homomorphism $\lambda:M\to M^*$ together with a decomposition
$M\cong\pi_1M\oplus\cdots\oplus \pi_\mu M$. The component $\pi_i M$
denotes the middle homology of the $i$th component of a Seifert surface.
Recall further that $G^{\epsilon,\mu}(\Z)$ denotes the Witt group of
$\mu$-component Seifert forms. Theorem~\ref{Flink_cobordism:seifert}
states that $C_{2q-1}(F_\mu)\cong G^{\epsilon,\mu}(\Z)$ for all $q\geq3$.

\begin{proposition}
\label{Flink_signatures}
Let $\mu\geq1$. The kernel of the canonical map
\begin{equation}
\label{Flink_Z_to_C}
G^{\epsilon,\mu}(\Z)\to G^{\epsilon,\mu}(\C^-)
\end{equation}
is $8$-torsion. Moreover, $G^{\epsilon,\mu}(\C^-)$ is a free abelian
group.
\end{proposition}
\subsection{Varieties of Signatures}
Although the cardinality of $G^{\epsilon,\mu}(\C^-)$ does not depend
on~$\mu$, $G^{\epsilon,1}(\C^-)$ is much smaller than, say,
$G^{\epsilon,2}(\C^-)$ in a sense which we explain next.

Our proof of proposition~\ref{Flink_signatures} will establish a
bijection between a basis for $G^{\epsilon,\mu}(\C^-)$ as a free abelian group
and the isomorphism classes of self-dual simple finite-dimensional 
complex representations of a certain ring~$P_\mu$.\index{$P_\mu$} In
other words, we associate to each\footnote{
Note that $\sigma_{M,b}$ is equal to the complex conjugate
$\sigma_{\overline M,\overline b}$ so if $M\ncong \overline{M}$ then
we are really defining one signature to correspond to the pair
$\{M,\overline M\}$.}
such isomorphism class $M$ one signature invariant
\begin{equation*}
\sigma_{M,b}:C_{2q-1}(F_\mu)\to\Z,
\end{equation*}
the composition of~(\ref{Flink_Z_to_C}) with a projection onto a summand. 
The second subscript $b:M\to M^*$ is one of two possible forms and
serves to specify a choice of sign for the signature. Note also that a simple
representation $M$ must be `algebraically integral' if the
corresponding signature of an $F_\mu$-link is to be non-trivial.

The ring $P_\mu$, which is non-commutative when $\mu\geq2$, 
was introduced by M.Farber~\cite{Far92}; we explain its role
and outline the proof of proposition~\ref{Flink_signatures} 
in section~\ref{section:defining_invariants} below. 

We implicitly assume that all representations are finite-dimensional.
In the case $\mu=1$, $P_\mu$ is just the polynomial ring $\Z[s]$ with
the involution $s\mapsto 1-s$. By the fundamental theorem of algebra,
the simple complex representations of $\Z[s]$ are
precisely the one-dimensional representations where the action
of~$s$ is given by multiplication by a complex number $\nu$. The 
representation is self-dual if and only if
$\overline\nu=1-\nu$, i.e. if and only if $\nu$ lies on
the line\footnote{
From the point of view of Blanchfield forms,\index{Blanchfield form}
one considers representations of $\Z[z,z^{-1}]$ in place of $\Z[s]$ where the
involution on $\Z[z,z^{-1}]$ is $z\mapsto z^{-1}$. 
A simple complex representation is one-dimensional and the action of
$z$ is multiplication by a complex number $\nu'$. The representation is
self-dual if and only if $\overline{\nu'}={\nu'}^{-1}$ so one
signature is defined for each point on the unit circle.
To relate signatures of Seifert and Blanchfield forms,
the simple representation $s\mapsto\nu:\C\to\C$ gives rise to the
presentation~(\ref{knot_module_presentation}) with $\sigma=z\nu+(1-\nu)$.
The action of $z$ on $\Coker(z\nu+(1-\nu))$ is multiplication by 
$\nu'=\nu^{-1}(\nu-1)=1-\nu^{-1}$. The map $\nu\mapsto 1-\nu^{-1}$
sends the line $\{\frac{1}{2}+ib\mid b\in\R\}$ to the unit circle.
For further discussion of the relationship between Seifert forms and
Blanchfield forms see for example Kearton~\cite{Kea75}, Levine~\cite{Lev77} or
Ranicki~\cite[Chapter 32]{Ran98}.
}
 $\{\frac{1}{2}+bi\mid b\in\R\}$.

To obtain a corresponding geometric description when $\mu\geq2$ we
begin by observing that $P_\mu$ is the path ring of a certain
quiver.\index{Quiver} A quiver, by definition, is a directed graph
possibly with
loops and multiple edges. The path ring\index{Path ring} of a quiver
is free as a $\Z$-module with one basis element for each path in the quiver;
the product of two paths is their concatenation if that makes sense
or zero if it does not. See chapter~\ref{chapter:preliminaries}, section~\ref{section:quiver} for further details.

A representation\index{Representation} of a quiver is a collection of
vector spaces, one for each vertex in the quiver, and a collection of
linear maps, one for each arrow (i.e.~directed edge). {\it The
representation theory of a quiver is identical to the representation
theory of its path ring.}

The particular quiver of which $P_\mu$ is the path ring is the `complete quiver
on $\mu$ vertices'.\index{Quiver!complete} It contains $\mu^2$ arrows, 
one arrow for each ordered pair of vertices. For example, in the case
$\mu=2$ the quiver has the following appearance:
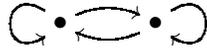
\begin{equation*}
\xymatrix{
\bullet \ar@/^/[r] \ar@(ul,dl)[] & \bullet \ar@/^/[l] \ar@(ur,dr)[]
}
\end{equation*}

To describe geometrically the simple
representations\index{Representation!simple (=irreducible)} of $P_\mu$ one
must first specify the dimensions of the vector spaces one wishes to
associate to the vertices. These dimensions are usually written as a `dimension
vector'\index{Dimension vector} $\alpha=(\alpha_1,\cdots,\alpha_\mu)$. Using
geometric invariant theory, the isomorphism classes of semisimple
dimension vector $\alpha$ representations of $P_\mu$ correspond to the points
on an (irreducible) affine algebraic variety
$\cy{M}(P_\mu,\alpha)$.\index{$M(R,\alpha)$@$\cy{M}(R,\alpha)$, $\sdM(R,\alpha)$} 
This variety is rarely smooth.\index{Variety of representations} 

If there are any simple representations of dimension vector~$\alpha$, then
almost all are simple. More precisely the simple isomorphism classes
are the points in a
(Zariski) open smooth subvariety of $\cy{M}(P_\mu,\alpha)$, the top
(Luna) stratum\index{Luna stratum} in a partition of $\cy{M}(P_\mu,\alpha)$ 
into `representation types' (Le~Bruyn and Procesi~\cite{LeBPro90}). In the case
$\mu=1$, only the dimension vector
$\alpha=(1)$ admits simple representations as we have discussed before. On
the other hand, when $\mu\geq2$, most dimension vectors admit
simple representations (see lemma~\ref{nonempty}).

There is a restriction to impose, for we wish to isolate {\it
self-dual} representations\index{Representation!self-dual}
of~$P_\mu$. The duality functor\index{Duality functor} $M\mapsto M^*$
induces an involution\index{Involution} on each variety
$\cy{M}(P_\mu,\alpha)$ and the
fixed point set $\sdM(P_\mu,\alpha)$ of this
involution turns out to be a real algebraic variety whose (real)
dimension coincides with the
(complex) dimension of $\cy{M}(P_\mu,\alpha)$. In summary,
\begin{proposition}
\label{variety_structure}
A basis of $G^{\epsilon,\mu}(\C^-)$ is in bijective correspondence with
a Zariski open subset of an infinite disjoint union of (absolutely
irreducible) real algebraic
varieties $\bigsqcup_{\alpha}{\sdM(P_\mu,\alpha)}$. The dimension of
$\sdM(P_\mu,\alpha)$ is 
$1+\sum_{1\leq i<j\leq\mu} 2\alpha_i\alpha_j~$.
\end{proposition}
In general an $F_\mu$-link may have non-zero signatures corresponding
to several different simple, self-dual representations, and they need not
lie on the same variety $\sdM(P_\mu,\alpha)$. If $(L,\theta)$ is a
split\index{Link!split} $F_\mu$-link then all its signatures lie on the
one-dimensional varieties $\sdM(P_\mu,\delta^i)$ where
$\displaystyle{\delta^i_j=\begin{cases}
1 & \mbox{if $i=j$}, \\
0 & \mbox{otherwise}.
\end{cases}
}$ Each of these varieties is a copy of the variety
$\{\frac{1}{2}+bi\mid b\in\R\}$ of knot signatures.
\subsection{Character}\index{Character}
Since we are associating signatures to representations of $P_\mu$ we
must learn to tell such representations apart. Moreover, if a signature
$\sigma_{M,b}:G^{\epsilon,\mu}(\Z)\to\Z$ is to be non-zero then $M$ must be a
summand of a representation which is induced up from an integral
representation:
\begin{equation*}
M\oplus M'\cong \C\otimes_\Z M_0.
\end{equation*}
We desire criteria to detect such `algebraically
integral'\index{Representation!algebraically integral}
representations $M$.

The character of a representation 
\begin{equation}
\label{character}
\chi_M:P_\mu\to \C^-;\quad r\mapsto \Trace(\rho(r)),
\end{equation}
which is an interpretation of M.Farber's
`trace invariant' (cf~\cite{Far92}, \cite{RRV99}), is helpful in both
respects. Firstly, a semisimple representation $M$ of
any ring over any field of characteristic zero is determined by its
character $\chi_M$ (see chapter~\ref{chapter:characters}). The
representation is self-dual if and only if $\chi_M$ respects involutions. 

Although $\chi$ may appear to be an infinite entity
one can decide whether or not two representations of~$P_\mu$ are isomorphic by
a finite number of comparisons. Two $m$-dimensional
complex representations $(M,\rho)$ and $(M',\rho')$ of~$P_\mu$ are
isomorphic if and only if $\chi_M(w)=\chi_{M'}(w)$ for all oriented
cycles $w$ of length at most $m^2$ in the quiver (see
Formanek~\cite{For86}).

The character also detects algebraically integral representations; a
complex representation (of any ring) is algebraically integral
if and only if the character takes values in the integers of some
algebraic number field.

\section{Number Theory}
Invariants which distinguish torsion elements of $C_{2q-1}(F_\mu)$
can be obtained by changing the ground ring not
to $\C^-$ but to $\Q$. 
\begin{Kexample}
The canonical map
\begin{equation*}
G^{\epsilon,1}(\Z)\to G^{\epsilon,1}(\Q)
\end{equation*}
is injective. Moreover, $G^{\epsilon,1}(\Q)$ is isomorphic to an
infinite direct sum of Witt groups of hermitian\footnote{
To be precise, there is one summand
$W^1\left(\frac{\Q[s]}{s-\frac{1}{2}}\right)=W^1(\Q)$ which obviously has 
trivial involution. However the projection of
$G^{\epsilon,1}(\Z)$ onto this summand is zero because $\frac{1}{2}$
is not an integer - see chapter~\ref{chapter:complete_invariants},
section~\ref{section:localization_exact_sequence}.}
forms over algebraic number fields (with non-trivial involution):
\begin{equation*}
G^{\epsilon,1}(\Q)\cong \bigoplus_p W\left(\frac{\Q[s]}{p}\right)
\end{equation*}
There is one summand for each maximal ideal $(p)\vartriangleleft\Q[s]$
which is self-dual in the sense that
$(p(s))=(p(1-s))\vartriangleleft\Q[s]$.

Hermitian forms over algebraic number fields were first classified by
Landherr~\cite{Lan36}. The Witt class of such a form is determined
by the dimension modulo~$2$, the signatures if any, and the
discriminant\index{Discriminant}. Up to sign, the discriminant of a
form is just the determinant of any matrix which represents the form.
Precise definitions can be found in section~\ref{define_discriminant}
of chapter~\ref{chapter:complete_invariants}.
\end{Kexample}

\begin{proposition}
\label{complete_invariants_outline}
The canonical map
\begin{equation*}
G^{\epsilon,\mu}(\Z)\to G^{\epsilon,\mu}(\Q)
\end{equation*}
is injective. Moreover, $G^{\epsilon,\mu}(\Q)$ is isomorphic to an
infinite direct sum of Witt groups of finite-dimensional division
$\Q$-algebras with involution:
\begin{equation}
\label{Flink_decomposition}
G^{\epsilon,\mu}(\Q)\cong \bigoplus_M W^1(\End(M)).
\end{equation}
\end{proposition}
There is one summand in~(\ref{Flink_decomposition}) for each
isomorphism class~$M$ of simple $\epsilon$-self-dual
finite-dimensional rational representations of~$P_\mu$ - see
section~\ref{section:defining_invariants} below.
Note that the involution on the endomorphism algebra
$\End(M)$\index{$End(M)$@$\End(M)$} and the
isomorphism~(\ref{Flink_decomposition}) both depend a choice of
$\epsilon$-symmetric form $b:M\to M^*$.

Fortunately, complete invariants for the Witt groups of the
division algebras $\End(M)$ are available in the literature - for
example~\cite{Alb39,HamMad93,Lam73,Lew82',Lew82,Scha85}. We summarize
these invariants in theorem~\ref{table_of_invariants}.
For all classes of algebras but one, some combination of the following
suffices to characterize the Witt class of a
form: dimension modulo 2;
signatures; discriminant; Hasse-Witt\index{Hasse-Witt invariant}
invariant. If $\End(M)$ is a quaternion algebra\index{Algebra!quaternion}
with a `non-standard'\index{Involution!non-standard}
involution one requires a secondary invariant such as the Lewis
$\theta$-invariant\index{Lewis $\theta$-invariant} which is defined if
the primary invariants vanish. We discuss all of these invariants in
chapter~\ref{chapter:complete_invariants}
below. Using~(\ref{Flink_decomposition}) they lead to an algorithm of
to decide whether or not two $F_\mu$-links are cobordant when $q>1$.

Proposition~\ref{complete_invariants_outline} also explains a
qualitative difference between knot and $F_\mu$-link cobordism, the difference 
between Levine's theorem~\ref{knot_answer_up_to_isomorphism} and
theorem~\ref{the_answer_up_to_isomorphism}.
Whereas knot theory leads us to study Witt groups of number
fields~$K_i$ which have non-trivial involution, it turns out that
every finite-dimensional division algebra with involution appears 
as $\End(M)$ in~(\ref{Flink_decomposition}), with infinite multiplicity. 
In particular the Witt groups of all number fields
with trivial involution are summands of
$G^{\epsilon,\mu}(\Q)$. One can find symmetric forms of order $8$ in the Witt 
group of a number field if and only if $-1$ can be expressed
as a sum of four squares, but not as a sum of fewer than four squares.

\section{Defining Invariants}
\label{section:defining_invariants}
Let us outline the proofs of propositions~\ref{Flink_signatures}
and~\ref{complete_invariants_outline}.
By studying the structure of $G^{\epsilon,\mu}(\C^-)$ we shall define
signatures of $C_{2q-1}(F_\mu)\cong G^{\epsilon,\mu}(\Z)$,
one for each self-dual simple complex representation of~$P_\mu$.
These signatures will be defined in four steps. The first step is a
`change of variables'; the second `devissage' step decomposes Seifert
forms over~$\C^-$ into simple constituents; the third step is a Morita
equivalence which allows us to replace these simple constituents by
their endomorphism rings; the fourth is Sylvester's
theorem\index{Sylvester's theorem} $W^1(\C^-)\cong\Z$. 

Steps two and three are part of a general theory of hermitian forms
formulated by Quebbemann, Scharlau and Schulte~(see~\cite{QSS79}
and~\cite[Chapter7]{Scha85}). We set out the details in
chapters~\ref{chapter:morita_equivalence} and~\ref{chapter:devissage} below.
The first three steps will be recycled, substituting $\Q$ for $\C^-$,
when we define invariants to distinguish the torsion classes of
$G^{\epsilon,\mu}(\Z)$ and prove proposition~\ref{complete_invariants_outline}.

\subsection{Sylvester's Theorem}\index{Sylvester's theorem}
Since the notion of signature\index{Signature} is so important in this
work, we begin at the `fourth step', the theory of hermitian forms over~$\C^-$,
however elementary the topic may be considered. Further details can be
found in algebra textbooks such as Lang~\cite[p577]{Lan93}.

Suppose $M$ is a finite-dimensional complex vector space and
$\phi:M\to M^*$ is a non-singular hermitian form. With
respect to any basis $x_1,\cdots x_m$, the form~$\phi$ is represented by an
invertible matrix $T_{ij}=\phi(x_i)(x_j)$ such that
${\overline{T}}^t=T$. An orthogonal basis can be chosen such that
\begin{equation*}
\phi(x_i)(x_j)=
\begin{cases} 
1~\mbox{or}~-1 &\mbox{if $i=j$} \\
0 &\mbox{if $i\neq j$}
\end{cases}.
\end{equation*}
In other words, $T$ is congruent to a diagonal matrix with diagonal
entries $\pm1$. The number $m_+(\phi)$ of positive entries and the number
$m_-(\phi)$ of negative entries do not
depend on the choice of basis; they are well-defined
invariants which determine $\phi$ uniquely up to isomorphism. 

For the purposes of cobordism computations one aims to classify forms not up to
isomorphism but up to Witt equivalence.
The Witt group $W(\C^-)$ is by definition the semigroup of
non-singular hermitian 
forms modulo the subsemigroup of metabolic forms.
As we saw in section~\ref{section:introduce_seifert_form} of
chapter~\ref{chapter:introduction} a form is
metabolic if it is represented by a matrix
$T=
\left(\begin{smallmatrix}
0 & * \\ 
* & *
\end{smallmatrix}\right)$ where each $*$ denotes some square matrix. 
In terms of invariants a
non-singular hermitian form~$\phi$ over~$\C^-$ is metabolic if and
only if $m_+(\phi)=m_-(\phi)$.
\begin{definition}
The {\it signature} of a non-singular hermitian form is
\begin{equation*}
\sigma(\phi)=m_+(\phi)-m_-(\phi).  
\end{equation*}
\end{definition}
\begin{theorem}[Sylvester]
\label{sylvesters_theorem}
The signature defines an isomorphism
\begin{equation*}
\sigma:W(\C^-)\xrightarrow{\cong} \Z.
\end{equation*}
\end{theorem}
\subsection{Change of Variables}
\label{section:change_of_variable}
The first step in our definition of signature for a Seifert form is a
`symmetrization'. A more precise treatment than that in the present
section can be found in chapter~\ref{chapter:preliminaries},
section~\ref{section:define_seifert_form}. Let us concentrate
first on knot theory.

A Seifert form $\lambda:M\to M^*$, which is asymmetric, can be
replaced by two entities: the $\epsilon$-symmetric form
$\phi=\lambda+\epsilon\lambda^*:M\to M^*$
and an endomorphism $s=\phi^{-1}\lambda:M\to M$. 
The intersection form~$\phi$ is intrinsic to the Seifert surface
$V^{2q}$ while $s$ encodes homological information 
about the embedding of $V^{2q}$ in $S^{2q+1}$. The two are related by
\begin{equation}
\label{relate_phi_and_s}
\phi(sx)(x')=\phi(x)((1-s)x')~\mbox{for all $x,x'\in M$}.
\end{equation}

To rephrase things a little, a Seifert form $(M,\lambda)$ is
replaced by an integral representation $\Z[s]\to\End_\Z M$ of the
polynomial ring $\Z[s]$ and
an $\epsilon$-symmetric form $\phi:M\to M^*$ which respects the representation.
The involution on $\Z[s]$ is given by $s\mapsto 1-s$.
This change of variables leads to an identity of Witt groups:
\begin{equation*}
G^{\epsilon,1}(\Z)\cong W^\epsilon(\Z[s]\dash\Z)~.
\end{equation*}

One can perform essentially the same trick for all $\mu\geq1$.
We start with slightly more data, a homomorphism
$\lambda:M\to M^*$ and a system of $\mu$ projections $\pi_i:M\to M$.
The intersection form $\phi=\lambda+\epsilon\lambda^*$ respects the
projections,
$\pi^*_i\phi=\phi\pi_i$,
for there is no intersection between the homology classes of
distinct components of a Seifert surface (the matrix $T$ represented
$\phi$ in chapter~\ref{chapter:introduction},
section~\ref{section:Flink_seifert_form}). We define
$s=\phi^{-1}\lambda:M\to M$ just as in knot theory, and observe that
equation~(\ref{relate_phi_and_s}) still holds.

Let us again rephrase matters using a representation. One can adjoin to
the polynomial ring $\Z[s]$ a system of orthogonal idempotents
$\pi_1,\cdots,\pi_\mu$ obtaining a non-commutative ring
\begin{align*}
P_\mu~ &\cong~\Z\left\langle s,\pi_1,\cdots,\pi_\mu
\biggm| \pi_i^2=\pi_i; \pi_i\pi_j=0 \text{\
for \mbox{$i\neq j$}};\sum_{i=1}^\mu
\pi_i=1 \right\rangle. \\ 
&\cong~\Z[s] *_\Z \left(\prod_\mu \Z\right)
\end{align*}
With the appropriate involution on $P_\mu$, 
\begin{equation*}
s\mapsto1-s;\quad \pi_i\mapsto\pi_i~\mbox{for $1\leq i\leq\mu$},
\end{equation*}
a $\mu$-component Seifert form $(M,\lambda,\{\pi_i\})$ becomes an
integral representation $P_\mu\to \End_\Z M$ together with a form
$\phi:M\to M^*$ which respects the representation
(cf~Farber~\cite{Far92}). The corresponding
isomorphism of Witt groups is the following:
\begin{equation*}
\kappa:G^{\epsilon,\mu}(\Z)\xrightarrow{\cong} W^\epsilon(P_\mu\dash\Z)
\end{equation*}
where $W^\epsilon(P_\mu\dash\Z)$ is the Witt group of triples
\begin{equation*}
(M~,~\rho:P_\mu\to \End_\Z M,~\phi:M\to M^*).
\end{equation*}

\subsection{Devissage}\index{Devissage}
\label{subsection:devissage}
The second step in the definition of signature is a hermitian version of the
Jordan-H\"older theorem.\index{Jordan-H\"older theorem} Further
details can be found in chapter~\ref{chapter:devissage}.
The Jordan-H\"older theorem says that in a context where
representations admit finite composition series
\begin{equation*}
0=N_0\subset N_1\subset\cdots\subset N_l=M
\end{equation*}
such that each subfactor $N_i/N_{i-1}$ is simple, the subfactors
of any two composition series for the same representation are
isomorphic, possibly after reordering. 

The hermitian version we need identifies an element of the Witt group
$W^\epsilon(P_\mu\dash\C^-)$ with a direct sum of Witt classes of forms
defined over simple representations, the subfactors. 
Of course, a simple representation $M$ must be self-dual if there is
to be a non-singular form $\phi:M\to M^*$. Indeed, it turns out
that any non-self-dual subfactors `cancel out' in pairs. The
appropriate uniqueness theorem is expressed by the canonical isomorphism
\begin{equation}
\label{devissage}
W^\epsilon(P_\mu\dash\C^-) \cong \bigoplus_{M}W^\epsilon_M(P_\mu\dash\C^-)
\end{equation}
with one summand for each isomorphism class of simple self-dual
complex representations $M$ of $P_\mu$.
The summand $W^\epsilon_M(P_\mu\dash\C^-)$ is by definition
the Witt group of $\epsilon$-hermitian forms 
\begin{equation*}
\phi:M\oplus\cdots\oplus M \to \left(M\oplus\cdots\oplus M\right)^*~.
\end{equation*}
over $M$-isotypic representations $M\oplus\cdots\oplus M$.
\begin{notation}
Let us denote by $p_M$ the canonical projection
\begin{equation*}
W^\epsilon(P_\mu\dash\C^-)\twoheadrightarrow W^\epsilon_M(P_\mu\dash\C^-).
\end{equation*}
\end{notation}

A complex representation of $P_\mu$ can of course be expressed as a
module over $\C^-\otimes_\Z P_\mu$. In knot theory, 
$P_1=\Z[s]$ and $\C\otimes_\Z\Z[s]\cong \C[s]$ is a principal ideal
domain. J.Milnor~\cite{Mil69} used the structure theorem for modules
over principal ideal domains to prove a version of~(\ref{devissage})
when $\mu=1$.
  
\subsection{Hermitian Morita Equivalence}\index{Morita equivalence}
\label{subsection:morita_equivalence}
The third step, which we treat more carefully in
chapter~\ref{chapter:morita_equivalence}, is the following isomorphism
\begin{equation*}
W^\epsilon_M(P_\mu\dash\C^-)\cong W^1(\End(M))
\end{equation*}
which replaces an isotypic representation $M\oplus\cdots\oplus M$ on
the left hand side by 
\begin{equation*}
\Hom(M,M\oplus\cdots\oplus M)
\end{equation*}
on the right hand side. The latter admits a natural `composition' action by the
endomorphism ring $\End(M)$.

Since $M$ is a simple representation, the endomorphism ring $\End(M)$
is a division ring (Schur's lemma). All the representations with which
we are concerned are finite-dimensional so $\End(M)$ is also
finite-dimensional and is therefore isomorphic to $\C^-$. 

\begin{notation}
Let us denote the Morita equivalence by
\begin{equation}
\label{morita_eq}
\Theta_{M,b}:W^\epsilon_M(P_\mu\dash\C^-)\xrightarrow{\cong} W^1(\C^-).
\end{equation}
\end{notation}
Unlike~(\ref{devissage}) above, the isomorphism $\Theta_{M,b}$ depends on a choice of
$\epsilon$-hermitian form $b:M\to M^*$. There are precisely two
$1$-dimensional forms over $\C^-$, namely $\langle1\rangle$ and
$\langle-1\rangle$. It follows that there are precisely two choices
for $b$; in other words one must make a choice of sign in the
definition of each signature.

To summarize all four steps we have the following definition:
\begin{definition}\index{Signature}
\label{signature_definitions}
If $M$ is a finite-dimensional simple complex representation
of~$P_\mu$ and $\epsilon b^*=b:M\to M^*$ 
is a non-singular $\epsilon$-hermitian form then the signature
$\sigma_{M,b}$ is the composite
\begin{equation*}
G^{\epsilon,\mu}(\Z)\to 
G^{\epsilon,\mu}(\C^-)\xrightarrow{\kappa}W^\epsilon(P_\mu\dash\C^-)
\xrightarrow{p_M}
W^\epsilon_M(P_\mu\dash\C^-)\xrightarrow{\Theta_{M,b}}W^\epsilon(\C^-)
\xrightarrow{\sigma}\Z
\end{equation*}
We define the signature $\sigma_{M,b}(L,\theta)$ of an $F_\mu$-link to
be the corresponding signature of any Seifert form for $(L,\theta)$.\index{$Sigma$@$\sigma_{M,b}(L,\theta)$}
\end{definition}

\subsection{Number Theory Invariants}
To distinguish torsion classes in $G^{\epsilon,\mu}(\Z)$, we apply the
first step, the `change of variables', exactly as before. We replace
$\C^-$ by $\Q$ and perform
devissage and hermitian Morita, obtaining
\begin{equation*}
G^{\epsilon,\mu}(\Z)\to G^{\epsilon,\mu}(\Q)\xrightarrow{\kappa}
W^\epsilon(P_\mu\dash\Q) \cong \bigoplus_M W^\epsilon_M(P_\mu\dash\Q) \cong
\bigoplus_M W^1(\End(M)).
\end{equation*}
The map $\kappa$ is an isomorphism (lemma~\ref{seifert_and_P}) so
$\displaystyle{G^{\epsilon,\mu}(\Q)\cong \bigoplus_M W^1(\End(M))}$.
We denote the projection onto the $M$th component 
\begin{equation*}
p_M:W^\epsilon(P_\mu\dash\Q) \twoheadrightarrow W^\epsilon_M(P_\mu\dash\Q) 
\end{equation*}
and the Morita equivalence
\begin{equation*}
\Theta_{M,b}:W^\epsilon_M(P_\mu\dash\Q) \xrightarrow{\cong} W^1(\End(M)).
\end{equation*}

A straightforward lemma - lemma~\ref{rationals_injection} -
shows that the canonical map $G^{\epsilon,\mu}(\Z)\to
G^{\epsilon,\mu}(\Q)$ is injective, completing the proof of
proposition~\ref{complete_invariants_outline}.
We can finally state our main result:
\begin{maintheorem}
\label{maintheorem_invariants}
Suppose $(L^0,\theta^0)$ and $(L^1,\theta^1)$ are $(2q-1)$-dimensional
$F_\mu$-links, where $q>1$. For $i=0,1$ let $V^i$ denote a
$\mu$-component Seifert surface for $(L^i,\theta^i)$ and let 
$N^i=H_q(V^i)/{\rm torsion}$. Let $\lambda^i:N^i\to (N^i)^*$ be the
Seifert form corresponding to $V^i$. 

The $F_\mu$-links $(L^0,\theta^0)$ and $(L^1,\theta^1)$ are cobordant
if and only if for each finite-dimensional $\epsilon$-self-dual simple
rational representation $M$ of~$P_\mu$, the dimension modulo~$2$, the
signatures, the discriminant, the Hasse-Witt invariant and the Lewis
$\theta$-invariant of 
\begin{equation*}
\Theta_{M,b}\,p_M\,\kappa\, [\Q\otimes_\Z(N^0\oplus N^1,\lambda^0\oplus-\lambda^1)] \in W^1(\End(M))
\end{equation*}
are trivial (if defined). 
\end{maintheorem}
Note that one need compute only a finite number of invariants to
establish that two given $F_\mu$-links are cobordant. Note
also that if all the invariants vanish for some choice of $b:M\to
M^*$ then they vanish for all possible choices.
\section{Chapter Summaries}
In chapter~\ref{chapter:preliminaries} we collect some definitions and
basic properties of quivers, Grothendieck groups, Witt groups and
Seifert forms. We treat the `change of variables' outlined in
section~\ref{section:change_of_variable}.

Chapters~\ref{chapter:morita_equivalence} and~\ref{chapter:devissage} 
develop the machinery required to decompose $G^{\epsilon,\mu}(\C^-)$
and $G^{\epsilon,\mu}(\Q)$ as direct sums of Witt groups of division rings.
Chapter~\ref{chapter:morita_equivalence} concerns Morita equivalence,
the `third step' described in
section~\ref{subsection:morita_equivalence} above. We work in an
arbitrary additive category with hermitian structure and split
idempotents. 
Chapter~\ref{chapter:devissage} treats the `second
step', the devissage technique of section~\ref{subsection:devissage} above. 
Here, we work in any abelian category which has a hermitian structure
and in which objects have finite composition series. 

Chapter~\ref{chapter:variety} describes the variety structures for the
signatures, proving proposition~\ref{variety_structure}. 

Chapter~\ref{chapter:pfister} gives an algebraic proof, adapted from
work of Scharlau, that the order of every element in the kernel of the
canonical map $G^{\epsilon,\mu}(\Q)\to G^{\epsilon,\mu}(\C^-)$
is a power of~$2$. This implies that the signatures in
definition~\ref{signature_definitions} are a complete set of
torsion-free invariants. In fact we prove a generalization which
applies to arbitrary fields.

Chapter~\ref{chapter:characters} contains a proof that a semisimple
representation is determined up to isomorphism by its character.
Chapter~\ref{chapter:rationality_of_representations} uses the
character to give criteria for a semisimple representation to be
a summand of a rational or integral representation.

In chapter~\ref{chapter:smalldim} we compute explicitly two
examples of the real varieties $\sdM(P_2,\alpha)$ and the Zariski open
subsets of simple self-dual representations.
One signature of $C_{2q-1}(F_2)$ has been defined for each algebraic
integer, or complex conjugate pair of algebraic integers,  in
these open subsets.

Chapter~\ref{chapter:complete_invariants} describes invariants of Witt
groups of finite-dimensional division $\Q$-algebras, completing the
proofs of propositions~\ref{Flink_signatures}
and~\ref{complete_invariants_outline} and of
theorem~\ref{maintheorem_invariants}. We go on to discuss the
localization exact sequence\index{Localization exact sequence}
\begin{equation*}
0\to W^\epsilon(P_\mu\dash\Z)\to W^\epsilon_\Z(P_\mu\dash\Q) \to
W^\epsilon(P_\mu\dash\Q/\Z)\to\cdots
\end{equation*}
Following the work of Stolzfus~\cite{Sto77} on knot cobordism, one 
could use this exact sequence in an attempt to
describe $C_{2q-1}(F_\mu)$ with greater precision, 
to characterize the subgroup $G^{\epsilon,\mu}(\Z)\subset
G^{\epsilon,\mu}(\Q)$ in terms of invariants.

Finally, in chapter~\ref{chapter:endrings} we construct examples of
integral representations~$M$ of~$P_\mu$
\index{Representation!integral}%
such that the endomorphism ring of $\Q\otimes_\Z M$ is {\it any}
prescribed finite-dimensional division algebra with involution. Since
every class of division algebras arises, all the invariants in the
previous chapter are warranted. Chapter~\ref{chapter:endrings}
concludes with a proof of theorem~\ref{the_answer_up_to_isomorphism}.

\chapter{Preliminaries}
\label{chapter:preliminaries}
In this chapter we recall the definitions of the Grothendieck group and
Witt group of an associative ring and, more generally, of an abstract
category with sufficient structure. Our main example is the
category of finitely generated projective representations of an
associative ring~$R$ with involution over another associative
ring~$A$. In particular, we are concerned with representations of the 
particular ring $P_\mu$ over $\Z$ or over a field.

In section~\ref{section:define_seifert_form} we define 
$\mu$-component Seifert forms and prove that the Witt group of such
forms is isomorphic to the Witt group of representations of $P_\mu$.  
\section{Representations}\index{Representation|textbf}
\label{section:representations}
Let $A$ and $R$ denote associative rings with identity.
A {\it representation} $(M,\rho)$\index{$(M,\rho)$} of $R$ over $A$
will be a ring homomorphism $\rho:R\to \End_A(M)$ where $M$ is a
finitely generated projective $A$-module. 
We often write $M$ in place of $(M,\rho)$.
 
Let $(R\dash A)\proj$\index{$(R\dash A)\proj$} 
denote the category of representations of $R$ over
$A$. By definition a morphism $\theta:(M,\rho)\to (M',\rho')$ in
$(R\dash A)\proj$ is an $A$-module map $\theta:M\to M'$ such that $\rho'(r)\theta =
\theta\rho(r)$ for all $r\in R$.

A ring homomorphism $A\to B$ induces a functor 
\begin{equation*}
(R\dash A)\proj\to (R\dash B)\proj
\end{equation*}
as follows: If $(M,\rho)$ is a representation of $R$ over
$A$ then $(B\otimes_AM, \rho_B)$ is a representation over $B$ where 
\begin{align*}
\rho_B : R &\to \End_B(B\otimes_A M) \\
r &\mapsto 1\otimes \rho(r).
\end{align*}

In our application to boundary links, $A$ will be a field $A=k$ or the
integers $A=\Z$ and $R$ will be a path ring for a quiver:
\subsection{Representations of Quivers}
\label{section:quiver}
For the convenience of the reader we summarize the basic definitions
and properties of quivers, paraphrasing Benson~\cite[p99]{Ben95}. 
Further background can be found in the book by Auslander, Reiten and
Smal\o~\cite[\S3.1]{ARS95}. 
\begin{definition}\index{Quiver|textbf}
A quiver $Q$ is a directed graph, possibly with loops and
multiple directed edges (arrows). It has a finite set $Q_0$ of
vertices and a set $Q_1$ of
arrows. Each arrow $e\in Q_1$ has a head $h(e)\in Q_0$ and a tail
$t(e)\in Q_0$.
\end{definition}
A representation\index{Representation!of a quiver} of a quiver will be a system $\{M_x\}_{x\in Q_0}$ 
of projective $A$-modules together with a system of homomorphisms
\begin{equation*}
\{f_e: M_{t(e)}\to M_{h(e)}\}_{e\in Q_1}.
\end{equation*} 
If $A$ is a field, say $A=k$, then each representation $(M,\rho)$ has a
{\it dimension vector}\index{Dimension vector|textbf} 
\begin{equation*}
\alpha=(\dim_k M_x)_{x\in Q_0}\in
\Z_{\geq0}^{Q_0}.  
\end{equation*}
Over any $A$, a representation of a quiver can
be interpreted as a representation $\rho:\Z Q\to \End_A(M)$ of the
path ring $\Z Q$ which we define next. 

\subsection{The Path Ring of a Quiver}\index{Path ring|textbf}
The path ring $\Z Q$ is a free $\Z$-module with one basis element
for each non-trivial path 
\begin{equation*}
p=\left(\bullet \mapright{e_1} \bullet \mapright{e_2} \cdots \mapright{e_n}
\bullet\right)
\end{equation*}
in~$Q$ and one basis element $\pi_x$ corresponding to the trivial path
at each vertex $x\in Q_0$. If $p$ is non-trivial
we write $h(p)=h(e_n)$ and $t(p)=t(e_1)$ whereas for trivial paths we set
$h(\pi_x)=t(\pi_x)=x$.

If $h(q)=t(p)$ then the product $p.q$ is by definition `path $q$ followed
by path $p$', the composite of the paths in reverse order. If $h(q)\neq
t(p)$ then one defines $p.q=0$.
Note that the trivial paths $\{\pi_x\}_{x\in Q_0}$ form a system of orthogonal
idempotents in $\Z Q$
\begin{equation*}
\sum_{x\in Q_0}\pi_x=1;\quad \pi_x\pi_y=
\begin{cases}
\pi_x & \mbox{if $x=y$} \\
0     & \mbox{if $x\neq y$}.
\end{cases}
\end{equation*}
\begin{Kexample}
Let $Q$ be the quiver with one vertex and one arrow.
The path ring $\Z Q$ is (isomorphic to) the ring of polynomials
$\Z[s]$ in an indeterminate~$s$.
\end{Kexample}

A representation $(\{M_x\}_{x\in Q_0}, \{f_e\}_{e\in Q_1})$ of any
quiver~$Q$ can be rewritten 
\begin{equation*}
 (\bigoplus_{x\in Q_0}M_x , \rho)
\end{equation*}
where $\rho:\Z Q\to \End(\bigoplus_{x\in Q_0}M_x)$ and   
\begin{equation*}
\rho\left(\bullet \mapright{e_1} \bullet \mapright{e_2} \cdots \mapright{e_n}
\bullet\right)
\end{equation*}
is the composite of the corresponding maps 
\begin{equation*}
\bigoplus M_x \twoheadrightarrow M_{t(e_1)} \mapright{f_{e_1}} \cdots \mapright{f_{e_n}}
M_{h(e_n)} \hookrightarrow \bigoplus M_x.
\end{equation*}
Inversely, given a representation $(M,\rho)$ of $\Z Q$
one may define $M_x=\rho(\pi_x)M$ and let $f_e:M_{t(e)}\to M_{h(e)}$
be the restriction of $\rho(e)$. Note that
$\rho(e)$ maps $M_{t(e)}$ to $M_{h(e)}$ since
$\rho(e)=\rho(e_{h(e)}e)=\rho(ee_{t(e)})$.  
Henceforth we shall freely interchange the symbols $Q$ and $\Z Q$,
writing for example $(Q\dash A)\proj=(\Z Q\dash A)\proj$. 
\subsection{Boundary Links}
We require the following quivers in the study of boundary links (compare
Farber~\cite{Far92},~\cite{Far91})  :
\begin{definition}\index{$P_\mu$|textbf}
Let $P_\mu$ denote the quiver with $\mu$
vertices $x_1,\cdots,x_\mu$ and precisely one edge joining vertex
$x_i$ to vertex $x_j$ for each ordered pair $(i,j)$ with $1\leq
i,j\leq \mu$. 
\end{definition}
\noindent For example, in the case $\mu=2$ the quiver $P_2$ appears as follows:
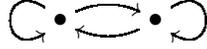
\begin{equation*}
\xymatrix{
\bullet \ar@/^/[r] \ar@(ul,dl)[] & \bullet \ar@/^/[l] \ar@(ur,dr)[]
}
\end{equation*}
As usual, $P_\mu$ will denote both the quiver and its path
ring.

\begin{lemma}
\label{P_is_path_ring}
The path ring $P_\mu$ is isomorphic to the free product of a
polynomial ring $\Z[s]$ by a product of $\mu$ copies of $\Z$. In symbols
\begin{align*}
P_\mu &\cong \Z[s] *_\Z \left(\prod_\mu \Z\right)  \\
&\cong \Z\left\langle s,\pi_1,\cdots,\pi_\mu
\biggm| \sum_{i=1}^\mu \pi_i=1, \pi_i^2=\pi_i, \pi_i\pi_j=0 \text{\
for \mbox{$1\leq i,j\leq \mu$}} \right\rangle.
\end{align*}
\end{lemma}
\begin{proof}
The indeterminate $s$ corresponds to the sum of all the paths of
length one. The idempotents $\pi_i$ correspond to the trivial paths
$\pi_{x_i}$. 
\end{proof}
 A representation of $P_\mu$ is to be thought of as the middle
homology $H_q(V;A)$ of a $\mu$-component embedded Seifert
surface
$V^{2q}$ - see chapter~\ref{chapter:introduction}, section~\ref{section:Flink_seifert_form} above and
lemma~\ref{seifert_and_P} below.
\section{Grothendieck Groups}\index{Grothendieck group}
\subsection{Rings}
Let $A$ be a ring, assumed to be associative and to contain a $1$.
Let $A\proj$\index{$A\proj$} denote the category of projective left $A$-modules. 
\begin{definition}\index{$K_0(A)$}
$K_0(A)$ is the abelian group with one generator $[M]$ for each isomorphism
class of finitely generated projective $A$-modules and one
relation $[M']=[M]+[M'']$ for each identity $M'\cong M\oplus M''$.
\end{definition}
A ring homomorphism $A\to B$ induces an additive functor
\begin{equation*}
A\proj\to B\proj;~M\mapsto B\otimes_A M
\end{equation*}
and therefore a group homomorphism $K_0(A)\to K_0(B)$.

\subsection{Categories}
\label{section:K_0(Cat)}\index{$K_0(\cy{C})$}
More generally, one can define the Grothendieck group $K_0(\cy{C})$ of
any category $\cy{C}$ which has a small skeleton and in which exact
sequences are defined (see Rosenberg~\cite[ch3]{Ros94}). However, the
following is general enough for our purposes: 
\begin{definition}
Suppose $\cy{C}$ is a full subcategory
of an abelian category $\cy{A}$.
The Grothendieck group $K_0(\cy{C})$ is the abelian group with one
generator $[M]$ for each isomorphism class of objects in $\cy{C}$ and one
relation $[M']=[M]+[M'']$ for each exact sequence
$0\to M\to M'\to M''\to0$
in $\cy{C}$.
\end{definition}

$K_0$ is functorial; if $F:\cy{C}\to \cy{D}$ is an exact functor,
i.e.~$F$ preserves exact
sequences, then there is an induced homomorphism 
$F:K_0(\cy{C})\to K_0(\cy{D})$ of abelian groups.

The most important example here is the category $\cy{C}=(R\dash A)\proj$
of representations defined in section~\ref{section:representations}:
A sequence of objects and morphisms
\begin{equation}
\label{equation:repn_ses}
0\to (M,\rho) \mapright{\theta} (M',\rho') \mapright{\theta'}
(M'',\rho'')\to 0
\end{equation}
is said to be exact if the underlying sequence $0\to
M\mapright{\theta} M' \mapright{\theta'} M''\to 0$ is exact.
$(R\dash A)\proj$ is a full subcategory of the abelian category $(R\dash A)\mod$
of representations of $R$ by arbitrary $A$-modules. We write
\begin{equation*}
K_0(R\dash A)=K_0((R\dash A)\proj).\index{$K_0(R\dash A)$}
\end{equation*}
For example, since an exact sequence of projective modules splits, we have
\begin{equation*}
K_0(\Z\dash A)=K_0(A\proj)=K_0(A).
\end{equation*}
On the other hand, if $R$ is any ring
and $A=k$ is a field then $(R\dash k)\proj$ is an abelian category and the
Jordan-H\"older theorem implies that $K_0(R\dash k)$ is a free abelian
group with one generator for each isomorphism class of simple objects. 

If $j:A\to B$ is a ring homomorphism then there is induced an exact functor
$(R\dash A)\proj\to (R\dash B)\proj$
(see section~\ref{section:representations}); for in any exact
sequence~(\ref{equation:repn_ses}) the 
underlying exact sequence of projective modules is split. Thus $j$
induces a group homomorphism $K_0(R\dash A)\to K_0(R\dash B)$.
\section{Witt Groups}\index{Witt group!of a ring with involution}
\subsection{Hermitian Forms over a Ring}
Let $A$ be an associative ring. 
\begin{definition}\index{$A^o$}
The opposite ring $A^o$ is identical to $A$ as an additive group but
multiplication in $A^o$ is reversed: $a\circ a'=a'a$. 
\end{definition}
\begin{definition}\index{Involution}
An involution on $A$ is an isomorphism $A\to A^o$, usually denoted $a\mapsto
\overline a$, such that $a=\overline{\overline a}$ for all $a\in A$.
\end{definition}
An involution on $A$ induces an equivalence of categories 
$\rproj A \to A\proj $; a right $A$-module $M$ becomes a left
$A$-module with $a.x=x\overline a$ for all $a\in A$ and $x\in M$.
In particular, if $M\in A\proj$ then
the dual module $M^*=\Hom_A(M,A)$ is again a finitely generated
projective left $A$-module with
$(a.\xi)(x)=\xi(x)\overline{a}$ for $a\in A$, $\xi\in
M^*$ and $x\in M$.  
We identify $M$ with $M^{**}$ via the natural isomorphism $M\to
M^{**};x\mapsto(\xi \mapsto \xi(x))$. 

Let $\epsilon=1$ or $-1$. An $\epsilon$-hermitian form\index{epsilon hermitian
form@$\epsilon$-hermitian form}\index{Hermitian!form} is a pair
$(M,\phi)$\index{$(M,\phi)$} where $M \in A\proj$ and $\phi:M\to M^*$ satisfies
$\phi^*=\epsilon\phi$. If $\phi$ is an
isomorphism then $(M,\phi)$ is said to be {\it non-singular}. The
category of non-singular $\epsilon$-hermitian forms is denoted
$H^\epsilon(A)$;\index{$H^\epsilon(A)$} a morphism $f:(M,\phi)\to
(M',\phi')$ in $H^\epsilon(A)$ is a map $f:M\to M'$ such that
$\phi=f^*\phi' f$.

If $(M,\phi)$ is an $\epsilon$-hermitian form and $j:L\hookrightarrow
M$ is the inclusion of a summand, one defines
\begin{equation*}
L^\perp~:=~\Ker\left(j^*\phi:M\to L^*\right).
\end{equation*} 
If $L=L^\perp$ then $L$ is called a {\it lagrangian} or {\it
metabolizer} and $(M,\phi)$ is called {\it metabolic}.
Equivalently, $(M,\phi)$ is metabolic if and only if  
\begin{equation*}\index{Metabolizer (=Lagrangian)}
(M,\phi)~\cong~
\left( L\oplus L^*,
\left(\begin{matrix}
0 & 1 \\ 
\epsilon & b
\end{matrix}\right)
:
L\oplus L^* \to L^*\oplus L
\right)
\end{equation*}
for some $L\in A\proj$ and some $\epsilon$-hermitian $b$.

\begin{definition} 
\label{Witt_group_of_ring}
The Witt group $W^\epsilon(A)$\index{$W^\epsilon(A)$} is the abelian
 group with one generator $[M,\phi]$ for each isomorphism class of
 non-singular $\epsilon$-hermitian forms in $H^\epsilon(A)$ subject to
 relations 
\begin{equation*}
\begin{cases}
{[M',\phi']=[M,\phi]+[M'',\phi'']}, 
&\mbox{if $(M',\phi')\cong
(M,\phi)\oplus (M'',\phi'')$} \\
{[M,\phi]=0}, &\mbox{if $(M,\phi)$ is metabolic}.
\end{cases}
\end{equation*}
\end{definition}
Thus two forms represent the same Witt class $[M,\phi]=[M',\phi']$ if
and only if there exist metabolic forms  $(H,\eta)$ and $(H',\eta')$
such that 
\begin{equation*}
(M\oplus H,\phi\oplus\eta)\cong (M'\oplus H', \phi'\oplus\eta').
\end{equation*}  
\begin{remark}
An $\epsilon$-hermitian form $(M,\phi)$ is {\it hyperbolic} if there
exists $L\in A\proj$ such that 
\begin{equation*}
(M,\phi)~\cong~ \left(L\oplus L^*, 
\left(\begin{matrix} 
0 & 1 \\
\epsilon & 0
\end{matrix}\right):L\oplus L^* \to (L\oplus L^*)^*\right)~.
\end{equation*}
Although metabolic forms are not in general hyperbolic, the same group
$W^\epsilon(A)$ is obtained if one substitutes the word `hyperbolic' for
`metabolic' in definition~\ref{Witt_group_of_ring}~\cite[p23]{Scha85}.
\end{remark}
\begin{remark}
The Witt group $W^{\epsilon}(A)$ is isomorphic to a symmetric
$L$-group $L^0_p(A,\epsilon)$ (Ranicki~\cite{Ran80}).
If there exists a central element $a\in A$ such that $a+\overline a=1$ 
then $W^{\epsilon}(A)$ is also isomorphic to the quadratic $L$-group
$L^p_{2q}(A,\epsilon(-1)^q)$ for all $q$.\index{$L_n(A)$}
\end{remark}
\begin{example}
\label{example:interchangerings}
For any ring $A$ the product $A\times A^o$ admits a transposition involution
$\overline{(a,a')}=(a',a)$, for which $W^\epsilon(A\times A^o)=0$.
\end{example}
\begin{proof}
Any $\epsilon$-hermitian form $(M,\phi)$ over $A\times A^o$ has
metabolizer $(1,0)M$.
\end{proof}
\begin{example}\index{$C^-$@$\C^-$, $\C^+$}
Let $\C^-$ denote the field of complex numbers with involution given
by complex conjugation. Let $\C^+$ denote the same field with trivial
involution. Then $W^\epsilon(\C^-) \cong \Z$ and $W^\epsilon(\C^+) \cong \Z/{2\Z}$.
\end{example}
\begin{remark}
If $A$ is a commutative ring then the group $W^1(A)$ is also a
commutative ring with multiplication given by the tensor product 
\begin{equation*}
(M,\phi)\otimes(M',\phi')=(M\otimes_A M',\phi\otimes\phi')
\end{equation*}
 where
$(\phi\otimes\phi')(x_1\otimes x_1')(x_2\otimes
x_2')=\phi(x_1)(x_2)\phi'(x_1')(x_2')$
for $x_1,x_2\in M$ and
$x_1',x_2'\in M'$. 
\end{remark}
\subsection{Change of Rings}
\label{section:Wittgroups_Changeofrings}
A homomorphism $j:A\rightarrow B$ of rings with involution
induces a functor 
\begin{align*}
H^\epsilon(A) &\to H^\epsilon(B) \\
(M,\phi) &\mapsto (B\otimes_A M,\phi_B)
\end{align*}
where  
\begin{align*}
\phi_B : B\otimes_A M &\rightarrow \Hom_B(B\otimes_A M, B); \\
b\otimes x &\mapsto (b'\otimes x' \mapsto b'j(\phi(x)(x'))\overline b). 
\end{align*}
for all $b,b'\in B$ and all $x,x'\in M$. 

If $(M,\phi)$ is metabolic then $(B\otimes_A M,\phi_B)$ is again metabolic so
$j$ induces a group homomorphism $W(A)\to W(B)$. 
If $A$ and $B$ are commutative then $j:W(A)\to W(B)$ is a ring
homomorphism.
\subsection{Hermitian Categories}\index{Hermitian!category}
Quebbemann, Scharlau and Schulte formulated a more general theory of
 quadratic and hermitian forms over a wider class of
 categories~\cite{QSS79}. We summarize next the definitions we require.

Suppose $\cy{C}$ is an additive category (see Bass~\cite[pp12-20]{Bas68}).
\begin{definition}
\label{defn:addcat}\index{Duality functor|textbf}
A {\it duality functor\ }  $*:\cy{C}\rightarrow \cy{C}$
is an additive contravariant functor together with a
natural isomorphism $(i_M)_{M\in\cy{C}}:id\rightarrow **$ such that 
$i^*_Mi_{M^*}=id_{M^*}$ for all $M\in\cy{C}$. We often omit $i$ identifying $M$
with~$M^{**}$. A triple~$(\cy{C},*,i)$ is called a {\it hermitian
category}, and is usually written $\cy{C}$ for brevity. 
\end{definition}
Our main examples of hermitian categories are $A\proj$ and $(R\dash A)\proj$:
\begin{example}
\label{examples_of_hermitian_categories}
The category  $A\proj$ admits the duality functor
 $M\mapsto M^*=\Hom_A(M,A)$. The category $(R\dash A)\proj$ of
 representations of a ring $R$ by f.g.~projective $A$-modules admits
 the duality functor $(M,\rho)\mapsto (M,\rho)^*= (M^*,\rho^*)$
where 
\begin{equation*}
\rho^*(r)(\xi)=(x\mapsto \xi(\overline r.x))
\end{equation*}
for $r\in R$, $\xi\in M^*$ and $x\in M$.
\end{example}
\begin{definition}
Let $\epsilon=+1$ or $-1$.
An $\epsilon$-hermitian form over $\cy{C}$
is\index{epsilon hermitian form@$\epsilon$-hermitian form}
by definition a pair $(M,\phi)$ where $M$ is an object of $\cy{C}$ and
$\phi:M\rightarrow M^*$ satisfies $\phi^*i_M=\epsilon\phi$. $(M,\phi)$ is
{\it non-singular\ } if $\phi$ is an isomorphism. 
\end{definition}
\begin{definition}
An object $M$ is called {\it self-dual} if $M\cong M^*$ and is called
$\epsilon$-self-dual if there exists a non-singular\index{epsilon self
dual@$\epsilon$-self-dual}\index{Self-dual}
$\epsilon$-hermitian form $(M,\phi)$.
\end{definition}
In many cases of interest, self-dual objects are both $1$-self-dual
and $(-1)$-self-dual - see lemmas~\ref{sd=epsilon_sd} and~\ref{sd=-sd} below.

The category of non-singular $\epsilon$-hermitian forms is denoted
$H^\epsilon(\cy{C})$.\index{$H^\epsilon(\cy{C})$}

\subsection{Change of Hermitian Category}
Hermitian categories are themselves the objects of a category. We
define next the morphisms between hermitian categories $\cy{C}$ and $\cy{D}$.
\begin{definition}
\label{duality_preserving_functor}\index{Duality preserving functor|textbf}
A {\it duality preserving functor} $\cy{C}\to \cy{D}$ is a triple
$(F,\Phi,\eta)$ where $F:\cy{C}\to\cy{D}$ is an additive functor,
$\{\Phi_M\}_{M\in\cy{C}}: F(\functor^*) \to F(\functor)^*$ 
is a natural isomorphism, $\eta=1$ or $-1$ and
\begin{equation}
\label{duality_functor}
\Phi_M^*i_{F(M)}= \eta \Phi_{M^*}F(i_M)
: F(M)\to F(M^*)^*
\end{equation}
for all $M\in \cy{C}$.
\end{definition}
A duality preserving functor $(F,\Phi,\eta)$ \index{$(F,\Phi,\eta)$}
is an equivalence of hermitian categories if and only if $F$ is an
equivalence of categories - see proposition~\ref{cat_eq_is_hermitian_eq} of
appendix~\ref{chapter:hermitian_cat}. Precise definitions of
composition and equivalence of duality preserving functors are also given
in appendix~\ref{chapter:hermitian_cat}.
\begin{example}
\label{ring_change_cat_example}
 Suppose as in section~\ref{section:Wittgroups_Changeofrings}
above that $j:A\to B$ is a ring homomorphism and
$B\otimes_A\functor:A\proj \to B\proj$ is the induced functor.
There is a duality preserving functor
$(B\otimes_A\functor,\Phi,1):A\proj\to B\proj$, where $\Phi$ is the
natural isomorphism  
\begin{align*}
\Phi_M:B\otimes_A M^* = B\otimes_A \Hom_A(M,A) &\to \Hom_B(B\otimes_A M,
B) = (B\otimes_A M)^* \\
 b\otimes \xi &\mapsto (b'\otimes x \mapsto b'\xi(x)\overline{b}).
\end{align*}

In the slightly more general setting of representation categories,
essentially the same triple $(B\otimes_A\functor,\Phi,1):(R\dash A)\proj\to
(R\dash B)\proj$ is a duality preserving functor.
\end{example}
\noindent An example in which $\eta$ can be $-1$ will appear in the next
chapter which concerns hermitian Morita equivalence. 

\begin{lemma}
\label{functor_of_forms}
A duality preserving functor $(F,\Phi,\eta):\cy{C}\to \cy{D}$ induces
an additive functor between categories of non-singular hermitian forms:
\begin{align*}
H^\epsilon(\cy{C}) &\to H^{\epsilon\eta}(\cy{D}) \\ 
(M,\phi) &\mapsto (F(M),\Phi_MF(\phi)).
\end{align*}

\end{lemma}
\begin{proof}
Suppose $(M,\phi)\in H^\epsilon(\cy{C})$ is an
$\epsilon$-hermitian form, so that $\phi^*i_M=\epsilon\phi$. 
To show $(F(M),\Phi_MF(\phi))\in H^{\epsilon\eta}(\cy{D})$ we
need $(\Phi_MF(\phi))^*i_{F(M)}=\epsilon\eta\Phi_MF(\phi)$. Indeed,
\begin{align*}
(\Phi_MF(\phi))^*i_{F(M)} &= F(\phi)^*\Phi_M^*i_{F(M)} \\
			  &= \eta F(\phi)^*\Phi_{M^*}F(i_M) \hspace{3mm}
			  \text{by equation~(\ref{duality_functor})} \\
			  &= \eta \Phi_M F(\phi^*)F(i_M) \hspace{3mm}
			  \text{by naturality of $\Phi$} \\
			  &=\eta \Phi_MF(\phi^*i_M) \\
			  &=\epsilon\eta\Phi_MF(\phi).\qedhere
\end{align*}
\end{proof}
\subsection{The Witt Group of a Category}\index{Witt group!of a
			  hermitian category|(} 
Suppose now that $\cy{C}$ is a hermitian category which is a full
subcategory of an abelian category~$\cy{A}$. Suppose further that $\cy{C}$ is
admissible in~$\cy{A}$, i.e.~that if $0\to M\to M'\to M''\to 0$ is an exact
sequence in $\cy{A}$ and $M'$ and $M''$ are in $\cy{C}$ then $M\in\cy{C}$.

A subobject $j: L\hookrightarrow M$ is called
admissible\index{Admissible subobject} in $\cy{C}$, 
if $\Coker(j:L\hookrightarrow M)\in\cy{C}$. In the case $\cy{C}=A\proj$ of
section~\ref{Witt_group_of_ring}, admissible subobjects are precisely
direct summands.

If $L$ is admissible one defines
\begin{equation*}
L^\perp~=~\Ker\left(j^*\phi:M\to L^*\right).
\end{equation*} 
which is an object of $\cy{C}$. We denote by $j^\perp$ the
inclusion~$L^\perp\hookrightarrow M$.  
If $L\subset L^\perp$ then $L$ is called a {\it sublagrangian\ }\index{Sublagrangian} of 
$(M,\phi)$ while if $L=L^\perp$ then one says $L$ is a {\it lagrangian\ } or 
{\it metabolizer\ }\index{Metabolizer (=Lagrangian)|textbf}
for $(M,\phi)$ and $(M,\phi)$ is {\it metabolic}.\index{Metabolic|textbf}

The definition of the Witt group of $\cy{C}$ is very similar to the
Witt group of a ring (see definition~\ref{Witt_group_of_ring} above):
\begin{definition}
\label{defn:Witt}
The Witt group $W^\epsilon(\cy{C})$\index{$W^\epsilon(\cy{C})$, $W^\epsilon(R\dash A)$} is the abelian group
 with one generator $[M,\phi]$ for each isomorphism class of non-singular
 $\epsilon$-hermitian forms $(M,\phi)\in H^\epsilon(\cy{C})$ subject
to relations 
\begin{equation*}
\begin{cases}
{[M',\phi']=[M,\phi]+[M'',\phi'']}, 
&\mbox{if $(M',\phi')\cong
(M,\phi)\oplus (M'',\phi'')$} \\
{[M,\phi]=0}, &\mbox{if $(M,\phi)$ is metabolic}.
\end{cases}
\end{equation*}
Two forms represent the same Witt class $[M,\phi]=[M',\phi']$ if
 and only if there exist metabolic forms  $(H,\eta)$ and $(H',\eta')$
such that 
\begin{equation*}
(M\oplus H,\phi\oplus\eta)\cong (M'\oplus H', \phi'\oplus\eta').
\end{equation*}  
\end{definition}
\begin{definition}
Let $W^\epsilon(R\dash A)$
denote the Witt group of the category of
representations:
\begin{equation*}
W^\epsilon(R\dash A):=W^\epsilon((R\dash A)\proj)
\end{equation*}
\end{definition}  
\noindent In particular, $W^\epsilon(\Z\dash A)=W^\epsilon(A)$.
\begin{remark}
If $A$ is a commutative ring then $W^\epsilon(R\dash A)$ is a module over
$W^1(A)$.
\end{remark}
\begin{lemma}
\label{induced_Witt_map}
An exact duality preserving functor $(F,\Phi,\eta):\cy{C}\to \cy{D}$ induces a
homomorphism of Witt groups:
\begin{align*}
W^\epsilon(\cy{C})&\to W^{\epsilon\eta}(\cy{D}) \\
[M,\phi] &\mapsto [F(M),\Phi_MF(\phi)]
\end{align*}
\end{lemma}
\begin{proof}
By lemma~\ref{functor_of_forms}, $(F,\Phi,\eta)$ induces an
additive functor $H^\epsilon(\cy{C})\to H^{\epsilon\eta}(\cy{D})$. We
need only prove that
if $L$ metabolizes $(M,\phi)$ then $F(L)$ metabolizes $(F(M),\Phi_MF(\phi))$.

Suppose $L=L^\perp$ so there is an exact sequence
\begin{equation*}
0\to L \mapright{j} M \mapright{j^*\phi} L^* \to 0~.
\end{equation*}
Since $F$ is exact the sequence
\begin{equation*}
0\to F(L)\mapright{F(j)} F(M)\mapright{F(j^*\phi)} F(L^*) \to 0
\end{equation*}
is also exact. Observing that 
\begin{equation*}
\Phi_L F(j^*\phi) = \Phi_L F(j^*)F(\phi) = F(j)^*\Phi_M F(\phi):F(M)\to
F(L)^*,
\end{equation*}
 the sequence
\begin{equation*}
0\to F(L)\mapright{F(j)} F(M)\mapright{F(j)\Phi_MF(\phi)} F(L)^* \to 0
\end{equation*}
is exact so $F(L)=F(L)^\perp$.
\end{proof}
For example, the duality preserving functor
$(B\otimes_A\functor,\Phi,1):(R\dash A)\proj\to (R\dash B)\proj$ of 
example~\ref{ring_change_cat_example} induces a homomorphism
$W^\epsilon(R\dash A)\to W^\epsilon(R\dash B)$. Setting $R=\Z$, one recovers the
homomorphism $W^\epsilon(A)\to W^\epsilon(B)$ defined in
section~\ref{section:Wittgroups_Changeofrings}.\index{Witt group!of a
			  hermitian category|)}
\smallskip

In the Witt group computations of chapter~\ref{chapter:devissage} one
needs the following two lemmas: 
\begin{lemma}
\label{lemma:sublagrangian}
If $(M,\phi)\in H^\epsilon(\cy{C})$ and $j:L\hookrightarrow M$ is a
sublagrangian, i.e.
\begin{equation*}
\phi_L=j^*\phi j=0:L\to L^*
\end{equation*}
then
\begin{equation*}
[M,\phi]=\left[\frac{L^\perp}{L},\overline\phi_{L^\perp}\right] \in
W^\epsilon(\cy{C}) 
\end{equation*}
where $\overline\phi_{L^\perp}:\frac{L^\perp}{L}\to \left(\frac{L^\perp}{L}\right)^*$ is induced from
the restriction $\phi_{L^\perp}=\left(j^\perp\right)^*\phi
j^\perp:L^\perp\to (L^\perp)^*$. 
\end{lemma}
\begin{proof}
The image of the diagonal map $\Delta:L^\perp\rightarrow M\oplus
L^\perp/L$ is a metabolizer for $(M,-\phi)\oplus
(L^\perp/L,\overline\phi_{L^\perp})$.
\end{proof}
\begin{lemma}
\label{orthogonal_decomposition}
If $(M,\phi)\in H^\epsilon(\cy{C})$ and $j:M'\hookrightarrow M$ is a
subobject such that the restriction $(M',\phi_{M'})$ is non-singular,
i.e.~$\phi_{M'}=j^*\phi j:M'\to {M'}^*$ is an isomorphism, then
\begin{equation*}
(M,\phi)~\cong~(M',\phi_{M'})\oplus ({M'}^\perp,\phi_{{M'}^\perp}).
\end{equation*}
\end{lemma}
\begin{proof}
The inclusion $j:M'\hookrightarrow M$ is split by $(j^*\phi
j)^{-1}j^*\phi: M\twoheadrightarrow M'$ the kernel of which is
${M'}^\perp$ so $M=M'\oplus {M'}^\perp$ as required. 
\end{proof}
\section{Seifert Forms}
\label{section:define_seifert_form}
As we outlined in sections~\ref{section:introduce_seifert_form}
and~\ref{section:algebraic_Flink_cobordism} of chapter~\ref{chapter:introduction}, the
$F_\mu$-link cobordism group $C_n(F_\mu)$ is isomorphic to a Witt
group of Seifert forms 
$G^{\epsilon,\mu}(\Z)$. We define this Witt group next.
\begin{definition}\index{Seifert!form|textbf}
Let $\epsilon=+1$ or $-1$ and let $A$ be any ring with involution. A
($\mu$-component) Seifert form over $A$ is a $(\mu+2)$-tuple
$(M,\pi_1,\cdots, \pi_\mu,\lambda)$ where $M$ is a finitely generated
projective $A$-module,
$\pi_1,\cdots,\pi_\mu$ is a set of orthogonal 
idempotents in $\End_A(M)$ and $\lambda:M\to M^*$ is an $A$-module
homomorphism such that $\lambda+\epsilon \lambda^*$ is an isomorphism
which commutes with each $\pi_i$.
The category of Seifert forms will be denoted $S^\epsilon(A,\mu)$.
We usually suppress the component structure writing a Seifert form
as $(M,\lambda)$.
\end{definition}
A Seifert form $(M,\lambda)$ is said to be
metabolic\index{Metabolic} if there exists a
metabolizer $L$ for the $\epsilon$-hermitian form $\phi:=
\lambda+\epsilon\lambda : M\to M^*$ such that
\begin{equation*}
\lambda(L)(L)=0 \quad \mbox{and}\quad L=L\cap \pi_1M \oplus \cdots \oplus L\cap \pi_\mu M.
\end{equation*}
\begin{lemma}
A summand $L$ of a Seifert form $(M,\lambda)$ is a metabolizer if and only if 
\begin{equation*}
L=L^{\perp_\lambda}=L^{\perp_\phi}\quad \mbox{and}\quad 
L=L\cap \pi_1M \oplus \cdots \oplus L\cap \pi_\mu M.
\end{equation*}
where $\phi=\lambda+\epsilon\lambda^*$ and $L^{\perp_\lambda}=\{x\in
M\ |\ \lambda(x)(y)=\lambda(y)(x)=0\ \forall y\in L\}$. 
\end{lemma}
\begin{proof}
We must show that if $L=L^{\perp_\phi}$ then $\lambda(L)(L)=0$
if and only if $L=L^{\perp_\lambda}$. 
Certainly $L=L^{\perp_\lambda}$ implies $\lambda(L)(L)=0$. Conversely,
if we have $\lambda(L)(L)=0$ then 
$L\subset L^{\perp_\lambda}\subset L^{\perp_\phi}=L$.
\end{proof}
\begin{definition}\index{$G^{\epsilon,\mu}(A)$|textbf}
The Witt group $G^{\epsilon,\mu}(A)$ of ($\mu$-component) Seifert forms is the
abelian group generated by isomorphism classes
$[M,\lambda]$ of $\mu$-component Seifert forms
subject to relations 
\begin{equation*}
\begin{cases}
{[M',\lambda']=[M,\lambda]+[M'',\lambda'']}, 
&\mbox{if $(M',\lambda')\cong
(M,\lambda)\oplus (M'',\lambda'')$} \\
{[M,\lambda]=0}, &\mbox{if $(M,\lambda)$ is metabolic}.
\end{cases}
\end{equation*}
\end{definition}
Recall from section~\ref{section:quiver} that $P_\mu$ denotes the
ring\index{$P_\mu$}
\begin{equation*}
P_\mu\cong \Z\left\langle s,\pi_1,\cdots,\pi_\mu
\biggm| \sum_{i=1}^\mu \pi_i=1, \pi_i^2=\pi_i, \pi_i\pi_j=0 \text{\
for \mbox{$1\leq i,j\leq \mu$}} \right\rangle.
\end{equation*}
which is the path ring of a quiver with $\mu$ vertices and $\mu^2$ arrows.
If one endows $P_\mu$ with the involution 
\begin{equation}
\label{Pmuinvolution}
\overline{\pi_i}=\pi_i; \hspace{10mm} \overline{s}=1-s
\end{equation}
then the Witt group $G^{\epsilon,\mu}(A)$ is in fact isomorphic to the
Witt group $W^\epsilon(P_\mu\dash A)$ of the representation category:
\begin{lemma}
\label{seifert_and_P}
The category $S^\epsilon(A,\mu)$ of Seifert forms
is equivalent to the category $H^\epsilon(P_\mu\dash A)$ of forms
over representations of $P_\mu$. Consequently, there is an isomorphism
\begin{equation*}
\kappa:G^{\epsilon,\mu}(A)~\xrightarrow{\cong}~W^\epsilon(P_\mu\dash A)
\end{equation*}
which is natural with respect to $A$.
\end{lemma}
\begin{proof}
Suppose $(M,\lambda)$ is a $\mu$-component Seifert form. 
We define a corresponding representation~$(M,\rho:P_\mu\to \End_AM)$ 
and a form $\phi:M\to M^*$ as follows:
\begin{align*}
\phi &=\lambda+\epsilon\lambda^* \\ 
\rho(s) &=\phi^{-1}\lambda \\
\rho(\pi_i^2) &=\rho(\pi_i):M\twoheadrightarrow M_i\rightarrowtail M.
\end{align*}
By definition $\rho(\pi_i)$ projects $M$ onto its $i$th component $M_i$.

Inversely, one can recover a $\mu$-component Seifert form from a triple 
$(M,\rho,\phi)$ by setting $\lambda=\phi\rho(s)$. For 
\begin{align*}
\lambda+\epsilon\lambda^* &=\phi\rho(s)+\epsilon(\phi\rho(s))^* \\
&= (1-\rho(s)^*)\phi + \rho(s)^*\phi \\
&= \phi.
\end{align*}

Metabolic Seifert forms correspond to metabolic
$(P_\mu\dash A)$-forms which implies the last sentence of the lemma.
\end{proof}

\chapter{Morita Equivalence}\index{Morita equivalence|(textbf}
\label{chapter:morita_equivalence}
This expository chapter concerns certain equivalences between
categories of modules and hermitian forms. 
\begin{example}
Suppose $M$ is a simple representation of a ring $R$ over a field
$k$. The endomorphism ring $E=\End_{(R\dash k)}M$ is a division algebra,
and there is a correspondence 
\begin{equation*}
{\overbrace{M\oplus\cdots\oplus M}^{r}=M^{\oplus r}} ~\longleftrightarrow~
E^r
\quad(r\in \N)
\end{equation*} 
between direct sums of copies of $M$ and (right) $E$-modules. 
\end{example}

If $R$ and $k$ are endowed with involutions and $M$ is a simple
representation which is $\epsilon$-self-dual, there
is a Morita equivalence 
between the category of $\epsilon$-hermitian forms 
$M^{\oplus r}\to \left(M^{\oplus r}\right)^*$ and the
category of (symmetric or) hermitian forms $E^{\oplus r}\to
\left(E^{\oplus r}\right)^*$.   

In chapter~\ref{chapter:devissage} the Witt groups
$W^\epsilon(R\dash k)$ of forms over representations of $R$ will be reduced
by a process of devissage to a direct sum of Witt groups of forms over 
isotypic representations $M^{\oplus r}$. Morita equivalence allows one
to pass to the Witt groups of the division algebra of endomorphisms of~$M$.
This equivalence is a generalization of the trace construction of
Milnor~(lemma 1.1 in~\cite{Mil69}) which played a crucial role in his
analysis of isometries of inner product spaces and hence in the computation of
high-dimensional knot cobordism.

In the present chapter, we proceed in greater generality. 
\begin{notation}
Let $\cy{C}$ be an additive category~\cite[pp12-20]{Bas68}.
Each object $M$ of $\cy{C}$ generates a full 
subcategory $\cy{C}|_M$ of objects isomorphic to a direct summand of
some direct sum  $M^{\oplus d}$ of copies of $M$. 
\end{notation}
This restricted category $\cy{C}|_M$ is Morita equivalent to the category of
projective right modules over the endomorphism ring of $M$:
\begin{equation*}
\cy{C}|_M \longleftrightarrow \rproj\End_\cy{C}(M)
\end{equation*} 
A similar result (theorem~\ref{hermitian_Morita_equivalence} below) applies
to the corresponding categories of hermitian forms.

References for this chapter include
Curtis and Reiner~\cite[Vol I,\S3D]{CurRei81} 
for the linear theory and Quebbemann Scharlau and
Schulte~\cite[p271]{QSS79}, Scharlau~\cite[\S7.4]{Scha85} or
Knus~\cite[\S I.9,ch.II]{Knu91} for the hermitian theory.
\section{Linear Morita Equivalence}
\begin{definition}
An additive category $\cy{C}$ is {\it idempotent complete} if each idempotent
$p^2=p:N\to N$ has a splitting $N\mapright{i}N'\mapright{j}N$ such
that $ji=p$ and $ij=\id$.
\end{definition}
All the categories we shall need, such as the category $(R\dash A)\proj$ of
representations of $R$ by finitely generated projective $A$-modules,
are idempotent complete.
\begin{theorem}[Linear Morita Equivalence] 
\label{linear_Morita_equivalence}
Suppose that $M$ is an object in an additive category $\cy{C}$ and
$\cy{C}|_M$ is idempotent complete. Let $E=\End_\cy{C}M$.
Then $\cy{C}|_M$ is equivalent
to the category $\rproj E$ of finitely generated projective
right $E$-modules.
\end{theorem}
\begin{proof}
$M$ admits an action on the left by $E$ 
and if $N\in\cy{C}$ then $\Hom_\cy{C}(M,N)$
is a right $E$-module under composition. 
The functor
\begin{equation*}
\Hom_\cy{C}(M,\functor):\cy{C}|_M  \rightarrow  \rproj E
\end{equation*}
has inverse $\functor\otimes_E M$.
\end{proof}
\begin{corollary}
If $\cy{C}$ is an exact category and $\cy{C}|_M$ is idempotent
complete then 
\begin{equation*}
K_0(\cy{C}|_M) \cong K_0(E^o).
\end{equation*}
\end{corollary} 
\begin{remark}[Naturality of Morita Equivalence]
If $F:\cy{C}\to \cy{D}$ is an additive (covariant) functor then $F$
induces a ring homomorphism $F:\End_\cy{C}M \to \End_\cy{D}F(M)$, and
the following square commutes up to natural isomorphism:
\begin{equation*}
\xymatrix{
{\cy{C}|_M} \ar[r] \ar[d]_\simeq & {\cy{D}|_{F(M)}} \ar[d]^\simeq \\
{\rproj\End_\cy{C}(M)} \ar[r] & {\rproj\End_\cy{D}F(M)}.
}
\end{equation*}
In particular, there is a commutative square
\begin{equation*}
\xymatrix{
K_0(\cy{C}|_M) \ar[r] \ar[d]_\cong & K_0(\cy{D}|_{F(M)}) \ar[d]^\cong \\
K_0({\End_\cy{C}(M)^o}) \ar[r] & K_0(\End_\cy{D}{F(M)^o}).
}
\end{equation*}
\end{remark}
\begin{proof}
The map
\begin{align*}
\Hom_\cy{C}(M,N)\otimes_{\End_\cy{C}M} \End_\cy{D}F(M)
&\to \Hom_\cy{D}(F(M),F(N)) \\
\phi\otimes f  &\mapsto F(\phi)f
\end{align*}
is natural for $N\in \cy{C}|_M$ and is an isomorphism because
\begin{equation*}
\End_\cy{C}(M)\otimes_{\End_\cy{C}M} \End_\cy{D}F(M) \cong
\End_\cy{D}F(M). \qedhere
\end{equation*} 
\end{proof}

The functors we shall be considering are those induced by a change
of ground ring:
\begin{align*}
F: (R\dash A)\proj &\to (R\dash B)\proj \\
(M,\rho) &\mapsto (B\otimes_A M,\rho_B)
\end{align*}
where $\rho_B(r)(b\otimes x)= b\otimes \rho(r)x$ for all $r\in R$,
$b\in B$ and $x\in M$.
There is a commutative diagram
\begin{equation*}
\xymatrix{
K_0\left((R\dash A)\proj|_M\right) \ar[r] \ar@{<->}[d]_\cong  & 
K_0\left((R\dash B)\proj|_{B\otimes M}\right) \ar@{<->}[d]^\cong \\
K_0\left((\End_{(R\dash A)}M)^o\right) \ar[r] &
K_0\left((\End_{(R\dash B)}B\otimes M)^o\right)~. 
}
\end{equation*}
In particular, if $A$ and $B$ are commutative and if $B$ is flat as an
$A$-module then there is an isomorphism 
\begin{align*}
B\otimes_A\End_{(R\dash A)}M &\cong  \End_{(R\dash B)}(B\otimes M) \\
b\otimes f  &\mapsto  (b'\otimes x \mapsto bb' \otimes f(x))~.
\end{align*}
\section{Hermitian Morita equivalence}
Suppose $\cy{C}$ is a hermitian category and $M\cong M^*$ is a
self-dual object in $\cy{C}$. The full subcategory $\cy{C}|_M$
inherits a hermitian structure. 
\begin{theorem}[Hermitian Morita Equivalence]
\label{hermitian_Morita_equivalence}
Let $\eta=+1$ or $-1$. Suppose that $(M,b)$ is a non-singular
$\eta$-hermitian form in a hermitian category $\cy{C}$, and assume
further that the hermitian subcategory $\cy{C}|_M$ is idempotent
complete. Let $E=\End_\cy{C}M$
be endowed with the adjoint involution $f\mapsto \overline f=b^{-1}f^*b$. 
Then there is  an equivalence of hermitian categories
\begin{equation*}
\Theta_{M,b}=(\Hom(M,\functor), \Phi, \eta) : \cy{C}|_M\to E\proj
\end{equation*}
\end{theorem}
The precise definition~\ref{equivalence_of_hermitian_categories} of 
equivalence of hermitian categories is given in
appendix~\ref{chapter:hermitian_cat}.
\begin{proof}[Proof of Theorem~\ref{hermitian_Morita_equivalence}]
Identifying $\rproj E$ with $E\proj$ there is a functor
\begin{equation*}
\Hom(M,\functor):\cy{C}|_M\to E\proj
\end{equation*}
which is an equivalence of categories by
theorem~\ref{linear_Morita_equivalence} above. We aim to 
extend this functor to an equivalence of hermitian categories.
By proposition~\ref{cat_eq_is_hermitian_eq} of
appendix~\ref{chapter:hermitian_cat} any 
duality preserving functor $(\Hom(M,\functor),\Phi,\eta)$ is
automatically an equivalence of hermitian categories.

One can define $\Phi:\Hom(M,\functor^*) \to
\Hom(M,\functor)^*$ 
by asserting that $\Phi_N$ should be the composite of the following natural
isomorphisms:
\begin{align}
\label{natural_isomorphism_Phi}
\Hom_\cy{C}(M,N^*)&\to\Hom_\cy{C}(N,M)\to \Hom_E(\Hom_\cy{C}(M,N),E) \\
\gamma &\mapsto b^{-1}\gamma^*;\hspace{8mm}\delta\mapsto
(\alpha\mapsto \delta\alpha). \notag\qedhere
\end{align}
\end{proof}
\begin{corollary}
If $\cy{C}$ is a hermitian admissible subcategory of an abelian
 category, $(M,b)$ is a non-singular $\eta$-hermitian form and 
$\cy{C}|_M$ is idempotent complete then there is an isomorphism
\begin{align*}
\Theta_{M,b}: W^\epsilon(\cy{C}|_M) &\xrightarrow{\cong}~
W^{\epsilon\eta}(\End_{\cy{C}}M); \\
[N,\phi] &\mapsto [\Hom(M,N),\Phi_N\phi_*]
\end{align*}
where $\phi_*:\Hom(M,N)\to\Hom(M,N^*)$ and 
$(\Phi_N\phi_*)(\alpha)(\beta)=\epsilon
b^{-1}\alpha^*\phi\beta \in \End_{\cy{C}}M$ for all $\alpha,\beta\in
\Hom(M,N)$.
\end{corollary}
\begin{proof}
The corollary follows from lemma~\ref{induced_Witt_map} because
$\Theta_{M,b}$ is an exact functor.
\end{proof}
\begin{Kexample}
A simple self-dual representation over~$\Q$ of~$P_1=\Z[s]$
may be written
\begin{equation*}
K=\frac{\Q[s]}{p(s)}\cong \frac{\Q[s]}{p(1-s)} = \Hom_\Q(K,\Q)
\end{equation*}
where $p\in\Q[s]$ is an irreducible polynomial. 

Milnor's trace construction (lemma 1.1 in~\cite{Mil69}) is a special
case of hermitian Morita equivalence. The trace construction turns a
symmetric form $\phi:K^r\to \Hom_\Q(K^r,\Q)$ which
respect the action of $s$ into a hermitian form over the field
$K=\frac{\Q[s]}{p}$ with involution $s\mapsto 1-s$. The key to the
construction is the following isomorphism
\begin{align*}
\label{Milnor_trace_isomorphism}
\Hom_K(K^r,K) &\xrightarrow{\cong}\Hom_\Q(K^r,\Q) \\
\delta &\mapsto \Trace_{K/\Q}\circ\,\delta. \notag
\end{align*}
This isomorphism is inverse to~(\ref{natural_isomorphism_Phi}) above
if one defines $b:K\to\Hom_\Q(K,\Q)$ by
$b(x)(y)=\Trace(x.y)$ for all $x,y\in K$.

The endomorphism ring of $K$ is just $K$ itself.
Since $K$ is commutative any choice of isomorphism 
$b:K\to \Hom_\Q(K,\Q)$ 
which respects the representations - i.e.~$b(sm)(m')=b(m)((1-s)m')$ -
will induce the same adjoint involution 
\begin{equation*}
K \to K; \quad s\mapsto b^{-1}s^*b = 1-s. 
\end{equation*}
\end{Kexample} 
\medskip
\begin{remark}[Naturality]
\label{naturality_of_hermitian_Morita_equivalence}
If $(F,\Phi,\eta):\cy{C}\to\cy{D}$ is a duality preserving functor
and $(M,b)\in H^\eta(\cy{C})$ then $(F(M),\Phi_MF(b))\in
H^\eta(\cy{D})$ and the following diagram of hermitian Morita
equivalences
\begin{equation*} 
\xymatrix{
{\cy{C}|_M} \ar[r] \ar@{<->}[d] & {\cy{D}|_{F(M)}} \ar@{<->}[d] \\
{\End_\cy{C}(M)\proj} \ar[r] & {\End_\cy{D}F(M)\proj}.
}
\end{equation*}
commutes up to natural isomorphism (of duality preserving functors).
If $F$ is exact then the following square also commutes:
\begin{equation*}
\xymatrix{
W^\epsilon(\cy{C}|_M) \ar[r] \ar@{<->}[d]_\cong &
W^\epsilon(\cy{D}|_{F(M)}) \ar@{<->}[d]^\cong \\
W^{\epsilon\eta}({\End_\cy{C}(M)}) \ar[r] & W^{\epsilon\eta}(\End_\cy{D}{F(M)}).
}
\end{equation*}
\end{remark}
In particular, a homomorphism $A\to B$ of rings with
involution induces a commutative diagram
\begin{equation*}
\xymatrix{
W^\epsilon((R\dash A)\proj|_M) \ar[r] \ar@{<->}[d]_\cong  & 
W^\epsilon((R\dash B)\proj|_{B\otimes M}) \ar@{<->}[d]^\cong \\
W^{\epsilon\eta}(\End_{(R\dash A)}M) \ar[r] &
W^{\epsilon\eta}(\End_{(R\dash B)}B\otimes M). 
}
\end{equation*}
As we remarked above, if $A\to B$ is a flat homomorphism of
commutative rings then $\End_{(R\dash B)}B\otimes M$ is naturally
isomorphic to $B\otimes_A\End_{(R\dash A)}M$.\index{Morita equivalence|)textbf}
\chapter{Devissage}\index{Devissage|(textbf}
\label{chapter:devissage}
In this chapter we assume that $\cy{C}$ is an abelian
category. Roughly speaking, an abelian category is a category in which every 
morphism has a kernel, an image and a cokernel. The standard reference
for a precise definition is the book of H.Bass~\cite[p21]{Bas68}.

We shall assume further that each object in $\cy{C}$ has finite
length. We discuss the Jordan-H\"older
theorem and the corresponding structure theorem in the hermitian
setting. The example to keep in mind is the representation category
$\cy{C}=(R\dash k)\proj$ where $k$ is a field. 

Recall from section~\ref{section:K_0(Cat)} of
chapter~\ref{chapter:preliminaries} that $K_0(\cy{C})$ is the
abelian group generated by the isomorphism classes of $\cy{C}$ with relations
corresponding to exact sequences in $\cy{C}$.
\begin{theorem}[Jordan-H\"older]\index{Jordan-H\"older theorem|textbf}
\label{linear_devissage}
Suppose $\cy{C}$ is an abelian category with ascending and descending
chain conditions. There is a canonical isomorphism
\begin{equation*}
K_0(\cy{C})\cong \bigoplus_{M\in\cy{M}^s(\cy{C})} K_0(\cy{C}|_M) 
\end{equation*}
where $\cy{M}^s(\cy{C})$ is the set of isomorphism classes of
simple objects $M$ in $\cy{C}$.
\end{theorem}
\begin{proof}
See, for example Lang~\cite[p22,157]{Lan93}.
\end{proof}
Every infinite ordered set contains either an infinite ascending chain 
$i_1<i_2<\cdots$ or an infinite descending chain $i_1>i_2>\cdots$. The
hypothesis that $\cy{C}$ has ascending and descending chain conditions
implies that every object $N$ of $\cy{C}$ has a finite composition
series
\begin{equation*}
0=N_0\subset N_1 \subset \cdots \subset N_l=N
\end{equation*}
such that each subfactor $N_i/N_{i-1}$ is simple. The content of
theorem~\ref{linear_devissage} is that these subfactors are uniquely
determined up to reordering. The natural number $l$ is also uniquely
determined and is called the composition length of $N$.

\begin{corollary}
$K_0(\cy{C})$ is isomorphic to a direct sum of copies of $\Z$, one for
each simple isomorphism class:
\begin{equation*}
K_0(\cy{C}) \cong \Z^{\oplus \cy{M}^s(\cy{C})}~\left(\cong
\bigoplus_{M\in\cy{M}^s(\cy{C})}K_0(\End_{\cy{C}}(M)^o)\right).
\end{equation*}
\end{corollary}

The corresponding theorem for Witt groups is the following:
\begin{theorem}[Hermitian Devissage]\index{Devissage!hermitian}
\label{hermitian_devissage}
Suppose $\cy{C}$ is an abelian hermitian category with ascending and
descending chain conditions.  There is an isomorphism of Witt groups
\begin{equation*}
W^\epsilon(\cy{C})~\cong \bigoplus_{M\in \ssdM(\cy{C},\epsilon)} W^\epsilon(\cy{C}|_M)
\end{equation*}
where $\ssdM(\cy{C},\epsilon)$\index{$M(C,epsilon)$@$\ssdM(\cy{C},\epsilon)$} denotes
the set of isomorphism class of simple $\epsilon$-self-dual objects in
$\cy{C}$.
\end{theorem}
\begin{corollary}
\label{general_devissage_and_ME}
Under the hypotheses of theorem~\ref{hermitian_devissage} The Witt
group of $\cy{C}$ is isomorphic to a direct sum of Witt groups of
division rings:
\begin{equation*}
W^\epsilon(\cy{C})~\cong
\bigoplus\limits_{M\in\ssdM(\cy{C},\epsilon)} W^1(\End_\cy{C}M)~.
\end{equation*}
\end{corollary}
\begin{proof}[Proof of corollary]
For each simple $\epsilon$-self-dual object $M$, choose an $\epsilon$-hermitian
form $(M,b)$ and apply theorem~\ref{hermitian_Morita_equivalence}.
\end{proof}

When $M$ is simple there is often little distinction between the
adjectives self-dual and $\epsilon$-self-dual (e.g.~\cite[p255]{Scha85}):
\begin{lemma}
\label{sd=epsilon_sd}
Suppose $M$ is a simple self-dual object in an abelian category $\cy{C}$ and
$2.\id:M\to M$ is an isomorphism. Then $M$ is either $1$-self-dual or
$(-1)$-self-dual or both.
\end{lemma}
\begin{proof}
Let $f:M\to M^*$ be an isomorphism. Since
\begin{equation*}
f=\frac{1}{2}(f+f^*) + \frac{1}{2}(f-f^*)
\end{equation*}
the forms $(f+f^*)$ and $(f-f^*)$ cannot both be zero, so one or other
is an isomorphism.
\end{proof}
Suppose next that $M$ is any $\epsilon$-self-dual object. If
$2.\id_M=0$ then there is no distinction between $\epsilon$-self-dual
and $(-\epsilon)$-self-dual. On the other hand, if $2.\id_M\neq0$ the
next lemma gives a criterion for $M$ to be $(-\epsilon)$-self
dual. Recall that an $\epsilon$-hermitian form $b:M\to M^*$ induces an
involution
\begin{align*}
I_b:\End(M) &\to \End(M) \\
f &\mapsto b^{-1}f^*b.
\end{align*}
\begin{lemma}
\label{sd=-sd}
An $\epsilon$-self-dual object $M$ {\rm fails} to be
$-\epsilon$-self-dual if and only if some (and therefore every)
$\epsilon$-hermitian form $b:M\to M^*$ induces the identity involution
$I_b=\id$ on $\End_\cy{C}M$.
\end{lemma}
\begin{proof}
Suppose $\epsilon b^*=b:M\to M^*$ and $f\neq I_b(f)$ for some
endomorphism $f\in \End_{\cy{C}}M$. Let $c=f-I_b(f)$ and observe that
$I_b(c)=-c$. The composite $b'=bc$ is a $(-\epsilon)$-hermitian form.

Conversely, suppose $b$ is an $\epsilon$-hermitian form and $b'$ is
a $(-\epsilon)$-hermitian form. Setting $f=b^{-1}b'$ we have
\begin{equation*}
b^{-1}f^*b=b^{-1}(b')^*b^{-*}b=b^{-1}(-\epsilon)b'\epsilon
b^{-1}b=-b^{-1}b'=-f.\qedhere
\end{equation*}
\end{proof}

In particular, lemmas~\ref{sd=epsilon_sd} and~\ref{sd=-sd} imply that
every self-dual simple complex representation is both $1$-self-dual and
$(-1)$-self-dual.
\begin{corollary}
$W^\epsilon(R\dash \C^-)$ is isomorphic to a direct sum of copies of $\Z$, one for
each isomorphism class of simple self-dual complex representations of $R$:
\begin{equation*}
W^\epsilon(R\dash \C^-)\cong \Z^{\oplus\ssdM(R\dash \C)}~.
\end{equation*}
\end{corollary}
\begin{proof}
Set $\cy{C}=(R\dash \C^-)\proj$ and observe that each of the endomorphism
rings $\End_\cy{C}M$ is isomorphic to $\C^-$. 
\end{proof} 
We have proved a part of proposition~\ref{Flink_signatures}:
\begin{corollary}
$G^{\epsilon,\mu}(\C^-)\cong \Z^{\oplus\infty}$.
\end{corollary}
A generalization of Pfister's theorem (see
chapter~\ref{chapter:pfister} below) will imply that the composite
\begin{equation*}
C_{2q-1}(F_\mu)\cong G^{\epsilon,\mu}(\Z)\to
G^{\epsilon,\mu}(\C^-)\cong \Z^{\oplus \infty}
\end{equation*}
with $q\geq2$ and $\epsilon=(-1)^q$) defines a {\it complete} set of
signature invariants.

We have also proved a part of proposition~\ref{complete_invariants_outline}:
\begin{corollary}
\label{devissage_over_Q}
There is an isomorphism
\begin{equation*}
G^{\epsilon,\mu}(\Q)\cong \bigoplus_M W^1(\End(M))
\end{equation*}
with one summand for each isomorphism class of $\epsilon$-self-dual
simple rational representation of~$P_\mu$.
\end{corollary}

To prove theorem~\ref{hermitian_devissage} we need the following lemma
\begin{lemma}
\label{lemma:unstable}
Suppose $\cy{C}$ is a hermitian abelian category. A stably
metabolic $\epsilon$-hermitian form over $\cy{C}$ is metabolic. In other
words, $(M,\phi)$ represents the zero element of
$W^\epsilon(\cy{C})$ if and only if $(M,\phi)$ is metabolic.
\end{lemma}
\begin{proof} (e.g.~\cite[Cor 6.4]{QSS79})
Suppose $L$ metabolizes $(H,\eta)$ and $L'$ metabolizes $(M\oplus
H,\phi\oplus\eta)$. Let $L''= L + L^\perp\cap L'\subset M\oplus H$ and observe that
\begin{align*}
{L''}^\perp= L^\perp \cap (L^\perp \cap L')^\perp = L^\perp \cap
(L + {L'}^\perp) &= L^\perp \cap (L + L') \\
&= L+L^\perp\cap L' \hspace{3mm} \text{since $L\subset L^\perp$} \\
&= L''
\end{align*}
so $L''$ metabolizes $(M\oplus H,\phi\oplus\eta)$. Since $L''$
contains $L$ we have $L''=L\oplus (L''\cap M)$ whence $L''\cap M$
metabolizes $(M,\phi)$.
\end{proof}
\begin{proof}[Proof of theorem~\ref{hermitian_devissage}]
A general element of $\bigoplus W^\epsilon(\cy{C}|_M)$ can be expressed
as a formal sum $\sum_{i=1}^r [N_i,\phi_i]$ where each form
$(N_i,\phi_i)$ represents an element of $W^\epsilon(\cy{C}|_{M_i})$.
Let
\begin{align*}
j:\bigoplus_M W^\epsilon(\cy{C}|_M) &\rightarrow W^\epsilon(\cy{C}) \\
\sum_{i=1}^r [N_i,\phi_i] &\mapsto \left[\bigoplus_{i=1}^r N_i, \bigoplus_{i=1}^r \phi_i \right].
\end{align*}
We check first that $j$ is well-defined and injective.
Indeed, if each $(N_i,\phi_i)$ is metabolic  with
metabolizer $L_i$ say then $L=\bigoplus L_i$ metabolizes $(N,\phi)=\left(\bigoplus_{i=1}^r N_i, \bigoplus_{i=1}^r \phi_i \right)$.
Conversely, suppose $L$ is a metabolizer for $(N,\phi)$. Since $N$ is
semisimple, $L$ is also semisimple. Moreover, each simple subobject of~$L$ is
isomorphic to one of the $M_i$ and is contained in the corresponding
$M_i$-isotypic summand $N_i\subset N$. 
Thus $L=\bigoplus_i L\cap N_i$ and $L\cap N_i$ metabolizes $(N_i,\phi_i)$.

It remains to show that $j$ is surjective. Let $(N,\phi)$ be a
representative of a Witt class in~$W^\epsilon(\cy{C})$ and let
$l=l(N)$ be the composition length of~$N$.
We argue by induction on~$l$. If $l$=1 then $N$ is simple and there is
nothing to prove.  If $l(N)>1$, let
$j:M_i\hookrightarrow N$
be the inclusion of a simple subobject. The composite
\begin{equation*}
M_i\hookrightarrow N\mapright{\phi} N^* \twoheadrightarrow M_i^*
\end{equation*}
is either the zero map or an isomorphism. 
In the former case
\begin{equation*}
[N,\phi] = \left[\frac{M_i^\perp}{M_i},\overline\phi_{M_i^\perp}\right]
\end{equation*}
by lemma~\ref{lemma:sublagrangian}. In the latter case the restriction
$(M_i,\phi_{M_i})$ is non-singular so
$(W,\phi)\cong(M_i,\phi_{M_i})\oplus (M_i^\perp, \phi_{M_i^\perp})$ by
lemma~\ref{orthogonal_decomposition}.
Either way, the inductive hypothesis applies.
\end{proof}
\begin{remark}[Naturality]
\label{naturality_of_devissage}
If $\cy{C}$ and $\cy{D}$ satisfy the hypotheses of
theorem~\ref{linear_devissage} and $F:\cy{C}\to\cy{D}$ is an exact
functor, then there is a commutative diagram
\begin{equation*}
\xymatrix{
K_0(\cy{C}) \ar@{<->}[d]_\cong \ar[r] & K_0(\cy{D}) \\
{\bigoplus_M K_0(\cy{C}|_M)} \ar[r] & {\bigoplus_M K_0(\cy{D}|_{F(M)})} \ar[u]
}
\end{equation*}
with both direct sums indexed over the simple isomorphism classes in
$\cy{C}$. 

Analogously, if $\cy{C}$ and $\cy{D}$ satisfy the hypotheses of
theorem~\ref{hermitian_devissage} and $(F,\Phi,\eta):\cy{C}\to \cy{D}$
is an exact duality preserving functor then there is a commutative
diagram
\begin{equation*}
\xymatrix{
W^\epsilon(\cy{C}) \ar[r] \ar@{<->}[d]_\cong & W^{\epsilon\eta}(\cy{D}) \\
{\bigoplus_M W^\epsilon(\cy{C}|_M)} \ar[r] & {\bigoplus_M
W^{\epsilon\eta}(\cy{D}|_{F(M)})}. \ar[u]
}
\end{equation*}
\end{remark}\index{Devissage|)textbf}
\chapter{Varieties of Representations}
\label{chapter:variety}\index{Variety of representations|(textbf}
The aim of this chapter is to give some geometric structure to the set of
signature invariants of Seifert forms over $\C^-$. In
chapter~\ref{chapter:rationality_of_representations} characters will
be used to identify the `algebraically integral' representations
relevant to Seifert forms over $\Z$ and hence to boundary
link cobordism.

\begin{notation}\index{$M(R)$@$\cy{M}(R)=\cy{M}(R\dash
\C)$}\index{$M(R)$, $M^s(R)$@$\sdM(R)$, $\cy{M}^s(R)$, $\ssdM(R)$}
Let $\cy{M}(R)=\cy{M}(R\dash \C)$ denote the set of
isomorphism classes of semisimple representations of a
ring $R$ over $\C$.
We denote the subset of self-dual
representations~$\sdM(R) \subset \cy{M}(R)$. 
Let $\cy{M}^s(R)\subset \cy{M}(R)$ be the subset of simple
representations and let $\ssdM(R)\subset \sdM(R)$ denote the simple
self-dual representations. 
\end{notation}
In this chapter we focus on the case $R=P_\mu$.

\begin{Kexample} 
When $\mu=1$ we have $P_1=\Z[s]$ and there is a
correspondence 
\begin{align*}
\C &\leftrightarrow  \cy{M}^s(\Z[s]) \\
\nu &\leftrightarrow  \C[s]/(s-\nu)
\end{align*}
which parameterizes the isomorphism classes of simple representations
and identifies $\cy{M}^s(\Z[s])$ with one-dimensional affine space. 

The self-dual representations $\sdM(\Z[s])$
correspond to the points $\nu\in\C$ which satisfy the equation
$\overline\nu = 1-\nu$. This equation defines a one-dimensional real variety
\begin{equation*}
\{1/2+ib\ |\ b\in \R\}\subset \C
\end{equation*}
(compare Ranicki~\cite[p609]{Ran98}). 
We recall the notion of a real variety in
definition~\ref{defn:realvariety} below.
\end{Kexample}

The aim here is to give some analogous algebraic structure to
$\cy{M}^s(P_\mu)$ and $\ssdM(P_\mu)$, when $\mu\geq 2$, in the 
framework of Mumford's geometric invariant theory~\cite{Mum94}. 
In chapter~\ref{chapter:smalldim}, explicit computations are given for two
particular dimension vectors $\alpha$.

Since $P_\mu$ is the path ring of a quiver 
(see lemma~\ref{P_is_path_ring}) a description of $\cy{M}^s(P_\mu)$ can be 
read off from the work of Le~Bruyn and Procesi~\cite{LeBPro90} (see
also~Procesi~\cite{Pro76}) who applied the \'Etale slice machinery of 
D. Luna~(\cite{Lun73},~\cite{Lun75}) to study the semisimple
representations of quivers.  

Given a quiver $Q$ and a fixed dimension vector $\alpha:Q_0\to
\Z_{\geq0}$ the isomorphism classes of semisimple 
representations of $Q$ with dimension vector~$\alpha$ correspond to the points
in an affine algebraic variety (which is not smooth in general).
The semisimple isomorphism classes $\cy{M}(P_\mu)$ are 
therefore a disjoint union of affine varieties $\cy{M}(P_\mu,\alpha)$, one for each dimension vector~$\alpha$. 
 
Each variety `admits a finite stratification into locally closed
smooth irreducible subvarieties corresponding to the different types
of semisimple decompositions of dimension vector $\alpha$. Moreover,
one stratum lies in the closure of another if the corresponding
representations are deformations' (\cite[p586]{LeBPro90}). 
In particular, if there exist simple representations with dimension vector
$\alpha$ then they form a dense (Zariski) open smooth subvariety of the
variety of semisimple representations.

A duality functor on $(Q\dash \C^-)\proj$ induces an involution on each
variety of semisimple representations. The subset invariant under the
involution is then a real algebraic variety.
\section{Existence of Simple Representations}
The first question is whether any self-dual simple representations
exist with dimension vector~$\alpha$.
When $\mu\geq2$ the answer is usually yes: 
\begin{lemma}
\label{nonempty}
Let $\alpha:\{x_1,\cdots,x_\mu\}\to \Z_{\geq0};~~x_i\mapsto \alpha_i$
be a dimension vector for the quiver $P_\mu$. The following are equivalent:
\begin{enumerate}
\item$\ssdM(P_\mu,\alpha)\neq \emptyset$.
\label{ssdne}
\item $\cy{M}^s(P_\mu,\alpha)\neq\emptyset$.
\label{sdne}
\item Either $|\Support(\alpha)|\geq 2$ or $\alpha=\delta^i$ for some $i$, where
$\delta^i_j=1$ if $j=i$ and $\delta^i_j=0$ if $j\neq i$.
\label{descdim}
\end{enumerate}
\end{lemma}
\begin{proof}
\ref{ssdne}$\Rightarrow$\ref{sdne}: Immediate.\\ \noindent
\ref{sdne}$\Rightarrow$\ref{descdim}: If ${\rm Support}(\alpha)=\{i\}$
then a representation of $P_\mu$ of dimension vector~$\alpha$ is essentially
a representation of a polynomial ring $\C[s_{ii}]$ which by the
fundamental theorem of algebra is not simple unless $\alpha_i=1$.\\ \noindent
\ref{descdim}$\Rightarrow$\ref{ssdne}:
The case where $\Support(\alpha)=\{i\}$ is easy so assume
$|\Support(\alpha)|\geq2$. We shall construct a simple self-dual
representation $(M,\rho)$ with dimension vector $\alpha$.

Recall that $s$ is the sum of the paths of length one in the quiver
$P_\mu$. We denote by $s_{ij}$ the arrow from vertex $x_j$ to vertex $x_i$.
For each $i\in\Support(\alpha)$ let
\begin{equation*}
\rho(s_{ii})=\left(
\begin{matrix}
\nu_1 & 0 &  \cdots & 0 \\
0 & \nu_2 &  \cdots & 0 \\
\vdots & \vdots & \ddots & \vdots \\
0 & 0  & \cdots & \nu_{\alpha_i}
\end{matrix}
\right)
\end{equation*}
where $\nu_1,\cdots, \nu_{\alpha_i}$ are distinct elements of
the set $\{1/2+bi\ |\ b\in\R\}$. When $i,j\in\Support(\alpha)$ and $i\neq j$, let
\begin{equation*}
\rho(s_{ij})={\rm sign}(i-j)\left(
\begin{matrix}
1 & 1 & \cdots & 1 \\
1 & 1 & \cdots & 1 \\
\vdots & \vdots & \ddots & \vdots \\
1 & 1 & \cdots & 1 
\end{matrix}\right).
\end{equation*}
The total matrix $\rho(s)_{ij}=\rho(s_{ij})$ satisfies
$\rho(s)=1-\overline{\rho(s)}^t$ and so $(M,\rho)$ is self-dual.

We must check that $M$ is simple. Regarding $\pi_iM$ as a
representation of a polynomial ring $\Z[s_{ii}]$ for each $i\in
\Support(\alpha)$, we have
\begin{equation*}
\pi_iM \cong \frac{\C[s_{ii}]}{(s_{ii}-\nu_1)}\oplus\cdots\oplus\frac{\C[s_{ii}]}{(s_{ii}-\nu_{\alpha_i})},
\end{equation*}
a direct sum of simple representations of $\Z[s_{ii}]$ no two of which are
isomorphic. Given a subrepresentation $M'\subset M$, it suffices to show that 
$\pi_iM'=0$ or $\pi_iM$ for each $i$.

Suppose that $\pi_iM'\neq0$. Since
$|\Support(\alpha)|\geq 2$ we can choose
an element $j\in\Support(\alpha)$ with $j\neq i$ and we find
\begin{equation*}
\langle(1,1,\cdots,1)\rangle\in \rho(s_{ij})\rho(s_{ji})(\pi_iM')\subset \pi_iM'
\end{equation*}
which implies that $\pi_iM'=\pi_iM$ as required.
\end{proof}
\section{Semisimple Representations of Quivers}
Let us fix a quiver~$Q$, a dimension vector $\alpha:Q_0\to \Z_{\geq0}$
and a family of
vector spaces $\{M_x=\C^{\alpha_x}\}_{x\in Q_0}$.
As usual, let
$m=\sum_{x\in Q_0}\alpha_x$.
Dimension vector $\alpha$ representations of~$Q$ correspond to families of
linear maps 
\begin{equation*}
\{f_e\}_{e\in Q_1} \in R(Q,\alpha)=\bigoplus_{e\in Q_1}
\Hom(\C^{\alpha_{h(e)}},\C^{\alpha_{t(e)}}),
\end{equation*}
points in an affine space of dimension
$\sum_{e\in Q_1} \alpha_{h(e)}\alpha_{t(e)}$. The coordinate
ring of this affine space is a polynomial ring $\C[X]$, where~$X$ is a
set of (commuting) indeterminates. 

Isomorphism classes of representations correspond to orbits for the
action of the reductive group
\begin{equation*}\index{$GL(alpha)$@$\GL(\alpha)$}
\GL(\alpha)=\prod_{i=1}^\mu \GL(\alpha_i)
\end{equation*}
by conjugation on $R(Q,\alpha)$, i.e.~$g.f_e=g_{t(e)}f_eg^{-1}_{h(e)}$
where $e\in Q_1$ and $g=\{g_x\}_{x\in Q_0}$. By Mumford's theory the
closed orbits correspond to semisimple isomorphism classes. 

The ring of invariants $\C[X]^{\GL(\alpha)}$ is the coordinate ring of
a variety 
\begin{equation*}\index{$M(R,\alpha)$@$\cy{M}(R,\alpha)$, $\sdM(R,\alpha)$}
\cy{M}(Q,\alpha)=R(Q,\alpha)//\GL(\alpha)
\end{equation*}
which parameterizes these semisimple isomorphism classes; the inclusion of
$\C[X]^{\GL(\alpha)}$ in $\C[X]$ induces a categorical quotient 
$R(Q,\alpha)\twoheadrightarrow \cy{M}(Q,\alpha)$,~the universal
$\GL(\alpha)$-invariant morphism. The Zariski topology on
$\cy{M}(Q,\alpha)$ coincides with the quotient topology so
$\cy{M}(Q,\alpha)$ is irreducible.

\begin{theorem}[Procesi and Le~Bruyn]
\label{oriented_cycles_generate}
$\C[X]^{\GL(\alpha)}$ is generated by traces of oriented cycles in the
quiver of length at most $m^2=(\sum\alpha_i)^2$. 
\end{theorem}
\begin{proof} 
See~\cite[Theorem 1]{LeBPro90}. 
\end{proof}
It follows from further work of Procesi~\cite{Pro87} that all the
relations between these traces can be deduced from Cayley-Hamilton
relations of $m\times m$ matrices.
\medskip

As we remarked above, there is a Luna stratification\index{Luna stratum} of
$\cy{M}(Q,\alpha)$ into smooth locally closed subvarieties, one for
each representation type:
\begin{definition}
A semisimple representation
\begin{equation*}
M=M_1^{\oplus d_1}\oplus\cdots\oplus M_k^{\oplus d_k}
\end{equation*}
is said to be of {\it representation type}\index{Representation!type}
$\tau=(d_1,\beta_1;\cdots;d_k,\beta_k)$ if each $M_i$ is simple and
$\beta_i$ is the dimension vector of $M_i$.
\end{definition}
The stratum of representations of type~$\tau'$ is contained in the closure
of the type~$\tau$ stratum if and only if $\tau'$ can be obtained
from~$\tau$ by a sequence of deformations each decomposing some
vector~$\beta_i$ as a sum of smaller dimension vectors 
$\beta_i=\beta^{(1)}_i+\beta^{(2)}_i$ (and perhaps regrouping if
$\beta^{(j)}_i=\beta_k$ for some~$j$ and~$k$). In particular, 
the stratum of simple isomorphism classes, if it exists, contains all other
strata in its closure; it is a dense open subset of $\cy{M}(Q,\alpha)$.
\begin{remark}
\label{dimension_of_varieties}
The dimension of $\cy{M}(Q,\alpha)$ is easy to calculate. We are
concerned here with those $\alpha$ for which the generic representation $M$
is simple so that $\dim(\Aut_{P_\mu}(M))=\dim(\End_{P_\mu}(M))=1$.

By the orbit-stabilizer theorem 
\begin{align*}
\dim\left(\cy{M}(Q,\alpha)\right) 
&=\dim\left(R(Q,\alpha)\right) -
\dim\left(\GL(\alpha)\right) + \dim(\GL(\alpha)_M) \\
\intertext{(where $\GL(\alpha)_M$ is the stabilizer of $M$)}
&= \sum_{e\in Q_1} \alpha_{h(e)}\alpha_{t(e)} - \sum_{x\in Q_0}
\alpha_x^2 +1.
\end{align*}
In particular when $Q=P_\mu$,
\begin{equation}
\label{dimension_of_variety}
\dim\left(\cy{M}(P_\mu,\alpha)\right) = \sum_{i,j=1}^\mu
\alpha_i\alpha_j - \sum_{i=1}^\mu\alpha_i^2 +1 = 1+ \sum_{1\leq
i<j\leq\mu}2\alpha_i\alpha_j
\end{equation}
\end{remark}
\section{Self-Dual Representations}
We proceed to give the self-dual isomorphism classes
$\sdM(P_\mu,\alpha)$\index{$M(R,\alpha)$@$\cy{M}(R,\alpha)$, $\sdM(R,\alpha)$} a real algebraic variety structure. Let us first
recall some of the basic theory of real algebraic geometry; for
further details the reader is referred to Bochnak,
Coste and Roy~\cite{BCR98} or Lam~\cite{Lam84}.

\subsection{Real Algebraic Sets}
Let $\R[X] =\R[X_1,\cdots,X_m]$ denote a commutative polynomial ring over the real numbers.
If $U\vartriangleleft \R[X]$ then one can define
\begin{equation*}
V(U)=\bigcap_{f\in U}\{x\in\C^m\mid f(x)=0\}.
\end{equation*}
The intersection $V_\R(U)=V(U)\cap\R^m$ is called a real
algebraic set;\index{Real!algebraic set} these generate the closed sets in the
Zariski topology for $\R^m$.  
In fact the points in $V_\R(U)$ correspond
to the maximal ideals of $\R[X]$ which both contain $U$ and are real in the
following (algebraic) sense (lemma~\ref{realideals} below):
\begin{definition}\index{Real!ring (=Formally real ring)}
\label{formally_real}
A commutative ring $A$ is {\it real} (also called {\it formally
real}) if it has the property
\begin{equation*}
a_1^2+\cdots +a_l^2=0 ~\Rightarrow~ a_1=\cdots=a_l=0~
\end{equation*}
for all $a_1,\cdots, a_l\in A$.
An ideal $U$ of a commutative ring $A$ is said to be real if and only if
the quotient ring $A/U$ is real. 
\end{definition}

On the other hand, if one is given a subset $S$ of~$\R^m$ one can define
\begin{equation*}
U(S)=\bigcap_{x\in S}\{f\in\R[X]\mid f(x)=0\}.
\end{equation*}
To explain the extent to which $V_\R$ and $U$ are mutually inverse, one
must introduce the real radical:
\begin{definition}\index{Real!radical}
The {\it real radical} $\sqrt[r]{U}$ of an ideal $U\vartriangleleft A$
is by definition 
\begin{equation*}
\{a\in A \mid a^{2m}+a_1^2+\cdots +a_l^2\in U \mbox{\ for some\ } m\geq
0 \mbox{\ and\ } a_1,\cdots,a_n\in A\}.
\end{equation*}
In fact $\sqrt[r]{U}$ is the intersection of the real ideals containing $U$.
\end{definition}
\begin{theorem}[Dubois-Risler Real Nullstellensatz] \hspace{2mm}\index{Real!nullstellensatz}
\smallskip \\ \noindent
i) If $U\vartriangleleft \R[X]$ then $U(V_\R(U))=\sqrt[r]{U}$. \smallskip \\ \noindent
ii) If $S\subset \R^m$ then $V_R(U(S))= \overline S$
where $\overline S$ denotes the closure of $S$ for the Zariski topology
on $\R^m$.
\end{theorem}

The dimension of a real algebraic set $V_\R(U)$ is by definition
the (Krull) dimension of $\R[X]/\sqrt[r]{U}$. Dubois and
Efroymson~\cite{Dub74} proved that given a real variety $V_\R(\mf{p})$ one
can find a chain of real varieties 
\begin{equation*}
V_\R(\mf{p})=V_\R(\mf{p}_d)\supset
V_\R(\mf{p}_{d-1})\supset\cdots\supset V_\R(\mf{p}_0) 
\end{equation*}
where $V_\R(\mf{p}_i)$ has dimension $i$; this result confirms that
the algebraic definition of dimension is reflected in the real geometry.  

The decomposition theory of semi-algebraic sets provides further
reassurance. A semi-algebraic set
is a subset of $\R^m$ defined by polynomial inequalities
in addition to, or in place of, polynomial equations. Any semi-algebraic set
decomposes~\cite[Thm2.3.6, \S2.8]{BCR98}
as a disjoint union of semi-algebraic sets each semi-algebraically 
homeomorphic to an open hypercube $(0,1)^d$. The largest occurring
value of~$d$ coincides with the Krull dimension defined above.  

\subsection{Real Varieties}
\begin{definition}\index{Real!variety}
\label{defn:realvariety}
An (affine) real algebraic variety $V$ is a subset of $\R^m$ of the
form $V=V_\R(\mf{p})$ where $\mf{p}$ is a real prime ideal of $\R[X]$.
\end{definition}
 Defined thus, a real variety is not only irreducible in the Zariski topology,
but absolutely irreducible (i.e.~$V(\mf{p})$ is also irreducible). Moreover
$V_\R(\mf{p})$ is Zariski dense in $V(\mf{p})$.

\subsection{Smoothness}
A point $\nu$ in an irreducible real algebraic set $V=V_\R(U)$ is
said to be {\it regular} if the localized ring 
\begin{equation*}
T=\left(\frac{\R[X]}{\sqrt[r]{U}}\right)_{(X_1-\nu_1,\cdots,X_n-\nu_n)}
\end{equation*}
is regular. 
The tangent space to $V$ at $\nu$ 
\begin{equation*}
\bigcap_{f\in U}\left\{x\in\R^m \biggm|
\sum_{i=1}^m\frac{\partial f}{\partial X_i}(\nu)x_i=0\right\}
\end{equation*}
coincides in dimension with $V$ if and only if $\nu$ is regular~\cite[Prop3.3.6]{BCR98} so a regular point is also described as
{\it non-singular} or {\it smooth}.

Smooth points give a useful test for reality:
\begin{proposition}
\label{simpleptcriterion}
(e.g.~Lam~\cite[p797]{Lam84}) Let $\mf{p}$ be a prime ideal in
$\R[X]$. Then $\mf{p}$ is real if and only if $V_\R(\mf{p})$ contains
a smooth point.  
\end{proposition}
\subsection{Varieties with Involution}
\label{section:variety_with_involution}
Returning to representations of $P_\mu$, the duality functor
$(\functor)^*$ on $(P_\mu\dash \C^-)\proj$ induces an involution on
$\cy{M}(P_\mu)$. 
To be more explicit, the involution~(\ref{Pmuinvolution})
on $P_\mu$ induces an involution
$I(\theta)=1-\overline{\theta}^{\raisebox{-1mm}{\scriptsize{$*$}}}$ on each
component $R(P_\mu,\alpha)$, and hence induces an involution on the
coordinate ring $\C[X]$. Here, $X=\{X_{ij}\}_{1\leq i,j\leq n}$ denotes
a set of $m^2$ commuting indeterminates and the induced involution
\begin{equation*}
I(X_{ij})=
\begin{cases}
-X_{ji} &\text{if $j\neq i$} \\
1-X_{ii} &\text{if $j=i$}
\end{cases}  
\end{equation*}
may be briefly written $I(f)(X)=\overline{f}(1-X^t)$. 

Let us combine the action of $\GL(\alpha)$ on $\C[X]$ with the involution.
\begin{definition}
Let $\GL(\alpha)\rtimes\frac{\Z}{2\Z}$ denote the semidirect product,
where $\frac{\Z}{2\Z}$ acts on $\GL(\alpha)$ by the formula
$I(g)=(g^{-1})^t$.
\end{definition}
\begin{lemma}
\label{semidirectaction}
There is an action of $\GL(\alpha)\rtimes\frac{\Z}{2\Z}$ on $\C[X]$
which extends both the conjugation action of $\GL(\alpha)$ and the
involution.
\end{lemma}
\begin{proof}
If $f\in\C[X]$ and $g\in\GL(\alpha)$ then 
\begin{align*}
(g.I(f))(X) = I(f)(g^{-1}X g) &= 
\overline f(1-(\overline{g}^{-1}X\overline{g})^t) \\ 
 &= (\overline{g}^{-1})^t.\overline{f}(1-X^t) = I((g^{-1})^t.f)(X).\qedhere
\end{align*}
\end{proof}
The following lemma says that $I$ induces an
involution on $\cy{M}(P_\mu,\alpha)=R(P_\mu,\alpha)//\GL(\alpha)$:
\begin{lemma}
The involution $I$ acts on the invariant ring $\C[X]^{\GL(\alpha)}$.
\end{lemma}
\begin{proof}
If $f\in\C[X]^{\GL(\alpha)}$ and $g\in \GL(\alpha)$ then
\begin{equation*}
(g.I(f))(X)=I((\overline{g}^{-1})^t.f)(X)=I(f)(X).\qedhere
\end{equation*}
\end{proof}
The orbits of the action of $I$ correspond to the maximal ideals in the
invariant ring 
\begin{equation*}
\left(\C[X]^{\GL(\alpha)}\right)^{\frac{\Z}{2\Z}}=\C[X]^{\GL(\alpha)\rtimes\frac{\Z}{2\Z}};
\end{equation*}
we are interested in the orbits which contain only one element, for
they correspond to self-dual isomorphism classes.

Let us now assume that there are simple representations of $P_\mu$ of
dimension vector~$\alpha$ (see lemma~\ref{nonempty}). 
\begin{proposition}
\label{realvariety}
The isomorphism classes of self-dual semisimple representations of
$P_\mu$ with dimension vector $\alpha$ correspond to the real maximal ideals of the invariant ring
$\C[X]^{\GL(\alpha)\rtimes\frac{\Z}{2\Z}}$ and hence to the points in a real
algebraic variety of dimension
\begin{equation*}
 1+ \sum_{1\leq i<j\leq\mu}2\alpha_i\alpha_j~.
\end{equation*}
\end{proposition}

\begin{proposition}
\label{simplesubvariety}
The subset $\ssdM(P_\mu,\alpha)$ of self-dual simple isomorphism classes
is a smooth Zariski open subvariety of $\sdM(P_\mu,\alpha)$.
\end{proposition}
\noindent The remainder of this chapter concerns the proof of
propositions~\ref{realvariety} and~\ref{simplesubvariety}.
 
Let $A$ denote any finitely generated
commutative $\C$-algebra with an involution $I$ which restricts to
complex conjugation on $\C$. Let $A_0=A^{\frac{\Z}{2\Z}}\subset A$ be the
involution-invariant $\R$-algebra. 
\begin{lemma}
\label{invariantsubring}
There is a natural isomorphism $A\cong \C^-\otimes_\R A_0$ of rings with
involution. In particular, $A$ is an integral extension of $A_0$
and $A_0$ is a finitely generated $\R$-algebra.
\end{lemma}
\begin{proof}
The natural map $\C^-\otimes_\R A_0 \to A;~\nu\otimes a \mapsto \nu a$
has inverse
\begin{equation*}
A \rightarrow \C\otimes_\R A_0;\hspace{6mm} 
a~\mapsto~1\otimes\frac{1}{2}((a+I(a))\,-\, i\otimes i(a - I(a)).\qedhere
\end{equation*}
\end{proof}

Recall that an ideal $\mf{m}$ of $A$ is said to {\it lie over} an ideal
$\mf{m}_0$ of $A_0$ if and only if $\mf{m}\cap A_0=\mf{m}_0$.
\begin{lemma}
\label{ideals_lying_over}
Over each maximal ideal $\mf{m}_0$ of $A_0$  
lies either one maximal ideal $\mf{m}\vartriangleleft A$ with
$\mf{m}=I(\mf{m})$ or precisely two distinct maximal ideals $\mf{m}$
and $I(\mf{m})$.
\end{lemma}
\begin{proof}
There exists a maximal ideal $\mf{m}\vartriangleleft A$ with $\mf{m}\cap
A_0=\mf{m}_0$ by the `lying over' theorem for integral
extensions. Although $I(\mf{m})$ may or may not be distinct
from~$\mf{m}$, $I(\mf{m})$ certainly lies over $\mf{m}_0$. It remains
to show that if 
$\mf{m}'\vartriangleleft A$ is maximal and $\mf{m}'\cap
A_0=\mf{m}_0$ then $\mf{m}'$ coincides either with $\mf{m}$ or
with $I(\mf{m})$. If we assume the contrary then 
by the Chinese remainder theorem there exists $a\in A$ such that 
\begin{equation*}
a\equiv 0 \pmod{\mf{m}},\quad
a\equiv 1 \pmod{\mf{m}'}\quad \mbox{and}\quad a\equiv 1\pmod{I(\mf{m}')}.
\end{equation*}
The condition
$a\in\mf{m}$ implies that $aI(a)\in \mf{m}\cap A_0=\mf{m}_0$
whereas the other conditions $a\equiv I(a)\equiv 1 \pmod{\mf{m}'}$ have
the contradictory implication $aI(a)\notin \mf{m}'$.
\end{proof}
\begin{lemma}
\label{realideals}
There is precisely one maximal ideal $\mf{m}$ over $\mf{m}_0$ if and only if $\mf{m}_0$ is real.
\end{lemma}
\begin{proof}
By lemma~\ref{invariantsubring} we may write
$A_0=\frac{\R[Y]}{U}$ and $A=\frac{\C[Y]}{\C\otimes U}$
for some finite set $Y=\{Y_1,\cdots,Y_m\}$ of commuting indeterminates and some
ideal $\mf{p}$. Now $\mf{m}$ has the form
$(Y_1-\nu_1,\cdots,Y_m-\nu_m)$ for some $\nu\in\C^m$ and
$\mf{m}=I(\mf{m})$ if and only if $\nu\in\R^m$. If $\nu\in\R^m$ then 
$\R[Y]/(\mf{m}\cap A_0)\cong \R$ is real. On the other
hand, if some $\nu_i\notin\R$ then 
\begin{equation*}
Y_i^2-(\nu_i+\overline{\nu}_i)Y_i +
\nu_i\overline{\nu}_i=(Y_i-\nu_i)(Y_i-\overline{\nu}_i)\in\mf{m}\cap A_0 
\end{equation*}
so in $A_0/(\mf{m}\cap A_0)$ we have
$(Y_i-\frac{1}{2}(\nu_i+\overline\nu_i))^2=\frac{1}{4}(\nu_i-\overline\nu_i)^2 < 0$.
\end{proof}

\begin{proof}[Proof of Proposition~\ref{realvariety}]
To obtain a correspondence between self-dual isomorphism classes of
representations and real maximal
ideals in the invariant ring $\C[X]^{\GL(\alpha)\rtimes\frac{\Z}{2\Z}}$, we apply
lemmas~\ref{ideals_lying_over} and~\ref{realideals}, putting
$A=\C[X]^{\GL(\alpha)}$. 
We must show that $\overline{\cy{M}}(P_\mu,\alpha)$ is a real algebraic
variety and compute its dimension.

Now $A=\C[X]^{\GL(\alpha)}$ is an integral domain so
the invariant ring $A_0=\C[X]^{\GL(\alpha)\rtimes \frac{\Z}{2\Z}}$ is of the
form $\frac{\R[Y]}{\mf{p}}$ where $\mf{p}$ is a prime ideal. 
We must check that $\mf{p}$ is real; by
proposition~\ref{simpleptcriterion} it suffices to prove that
$V_\R(\mf{p})$ contains a smooth point.

Lemma~\ref{nonempty} above assures us that there are self-dual
representations among the simple
dimension vector~$\alpha$ representations, so $\C[X]^{\GL(\alpha)}$ contains a
maximal ideal $\mf{m}$ which is both smooth and involution
invariant. 
The proof of lemma~\ref{invariantsubring} implies that
$\mf{m}=\C\otimes\mf{m}_0$ where $\mf{m}_0=\mf{m}\cap
\C[X]^{\GL(\alpha)\rtimes\frac{\Z}{2\Z}}$ so
\begin{multline*}
\dim_\R\left(\frac{\mf{m}_0}{\mf{m}_0^2}\right)
=
\dim_\C\left(\frac{\mf{m}}{\mf{m}^2}\right)
=
\Krulldim\left(\C[X]^{\GL(\alpha)}\right)
\\
= 
\Krulldim\left(\C\otimes_\R \C[X]^{\GL(\alpha)\rtimes\frac{\Z}{2\Z}}\right)
=
\Krulldim\left(\C[X]^{\GL(\alpha)\rtimes\frac{\Z}{2\Z}}\right).
\end{multline*}
Thus $\mf{m}_0$ is a smooth real point and so $\mf{p}$ is a real prime.
 
The (real) dimension of $\cy{M}(P_\mu,\alpha)$ is equal to the
(complex) dimension of $\cy{M}(P_\mu,\alpha)$ which, by
remark~\ref{dimension_of_varieties} is $1+ \sum_{1\leq
i<j\leq\mu}2\alpha_i\alpha_j$.
\end{proof}

\begin{proof}[Proof of Proposition~\ref{simplesubvariety}]
Note first that the 
involution~$I$ respects the Luna stratification of
$\cy{M}(P_\mu,\alpha)$; in other words $I$ preserves representation type. 
The inclusion of $\C[X]^{\GL(\alpha)\rtimes\frac{\Z}{2\Z}}$ in
$\C[X]^{\GL(\alpha)}$ induces a quotient
\begin{equation*}
q:\cy{M}(P_\mu,\alpha)\to \cy{M}(P_\mu,\alpha)//\frac{\Z}{2\Z}
\end{equation*}
and the Zariski topology on $\sdM(P_\mu,\alpha)$ coincides with the
subspace topology induced by the inclusion of $\sdM(P_\mu,\alpha)$ in
$\cy{M}(P_\mu,\alpha)//\frac{\Z}{2\Z}$.   
The image under $q$ of the Zariski open set
$\cy{M}^s(P_\mu,\alpha)$ is open in $\cy{M}(P_\mu,\alpha)//\frac{\Z}{2\Z}$ and
its intersection $\ssdM(P_\mu,\alpha)=\sdM(P_\mu,\alpha)\cap
q(\cy{M}^s(P_\mu,\alpha))$ is therefore open in $\sdM(P_\mu,\alpha)$.
Smoothness of $\ssdM(P_\mu,\alpha)$ follows (as in the proof of
proposition~\ref{realvariety}) from the smoothness of $\cy{M}^s(P_\mu,\alpha)$.
\end{proof}

\chapter{Generalizing Pfister's Theorem}
\label{chapter:pfister}
This chapter initiates discussion of rationality questions proving
that the varieties of signatures defined in the preceding chapters are
a {\it complete} set of torsion-free invariants for boundary link cobordism.  

To recap, in 
section~\ref{section:Flink_seifert_form} of chapter~\ref{chapter:introduction}
and section~\ref{section:define_seifert_form} of
chapter~\ref{chapter:preliminaries} we identified the odd-dimensional
$F_\mu$-link cobordism group $C_{2q-1}(F_\mu)$ with the Witt group
$W^\epsilon(P_\mu\dash \Z)$ where $\epsilon=(-1)^{2q-1}$
and chapters~\ref{chapter:morita_equivalence}
and~\ref{chapter:devissage} demonstrated an
isomorphism $W^\epsilon(R\dash \C^-)\cong \Z^{\oplus \ssdM(R)}$ for
any ring $R$ with involution. 
We will see in lemma~\ref{rationals_injection} that $W^\epsilon(R\dash \Z)$ injects 
into $W^\epsilon(R\dash \Q)$. In the present chapter we prove the following:
\begin{theorem}
\label{2primary_kernel}
Every element $v$ in the kernel of the natural map
\begin{equation}
\label{generic_Q_into_C}
W^\epsilon(R\dash \Q)\to W^\epsilon(R\dash \C^-)\cong
\Z^{\oplus \ssdM(R\dash \C)} 
\end{equation}
satisfies $2^nv=0$ for some integer $n$.
\end{theorem}
Using number theory one can prove that the kernel
of~(\ref{generic_Q_into_C}) is in fact $8$-torsion - see
corollary~\ref{eight_torsion} below.

In fact we shall prove a generalization of
theorem~\ref{2primary_kernel} substituting for $\Q$  
an arbitrary field $k$ (with trivial involution). 
This generalization is (Morita) equivalent to
a result of Scharlau~\cite[\S5]{Scha70} on Witt groups of semisimple
algebras. Indeed, our proof follows Scharlau's closely employing his
Frobenius reciprocity method.\index{Frobenius reciprocity}  
The special case $R=\Z$ is Pfister's theorem (see
Scharlau~\cite[Thm2.7.3 p56]{Scha85}):
\begin{theorem}[Pfister]\index{Pfister's theorem}
\label{pfister}
A symmetric form $\phi$ over a field $k$ represents a torsion element
of the Witt group $W^1(k)$ if and only if the signature of $\phi$ is
zero with respect to every ordering of $k$.
\end{theorem}  

\section{Artin-Schreier Theory}
It will be convenient to speak not of the orderings of a field, as in
theorem~\ref{pfister}, but of real closures. We
explain briefly the equivalence between the two ideas.
The reader is referred to Bochnak, Coste and Roy~\cite[Chapter
1]{BCR98}, Milnor and Husemoller~\cite[Chapter III \S2]{MilHus73} or
Scharlau~\cite[p113]{Scha85} for further details.
 
Let us first recall a definition of ordering for a field:
\begin{definition}
\label{ordering}\index{Ordered field}
An ordering of a field $k$ is a subset $\omega\subset k$ 
which is closed under addition and multiplication and is such that   
$k$ is the disjoint union:
\begin{equation*}
k~=~-\omega\sqcup \{0\} \sqcup \omega
\end{equation*}
where $-\omega=\{-x\mid x\in \omega\}$.
\end{definition}
Given an ordering $\omega$ one can define a relation $<_\omega$ on $k$
which satisfies all the usual axioms of an ordered set and moreover
respects addition and multiplication. Let
\begin{equation*}
x<_\omega y\quad\mbox{if and only if}\quad y-x \in \omega.
\end{equation*}
Inversely, given a relation $<$ one may set $\omega=\{x\in k\mid x>0\}$.

If a field $k$ admits an ordering, the non-zero squares are certainly
positive with respect to the ordering so $k$ is real. Conversely, in a
real field the subset of elements which can be expressed as a sum of
non-zero squares does not contain zero and can be extended,
using Zorn's lemma, to an ordering $\omega$. Note that $\omega$ is not 
in general unique.
\begin{definition}
A real field $k$ is said to be {\it real closed}\index{Real!closed
field} if no proper algebraic extensions of $k$ are real. A {\it real
closure}\index{Real!closure of a field} of a real field $k$ is an algebraic extension which is real closed.
\end{definition}
A real closed field has much in common with the field $\R$ of real
numbers. In particular, if $k$ is real closed then $k(\sqrt{-1})$ is
algebraically closed and $k$ admits the unique ordering
$\omega=\{x^2\mid x\in k^{\,\bullet}\}$.

A Zorn's lemma argument shows that every ordered field $(k, \omega)$ has a real
closure $k_\omega$ (whose unique ordering extends $\omega$). Conversely,
\begin{theorem}[Artin-Schreier]
\label{Artin-Schreier}
There is precisely one real closure of an ordered field up to
$k$-algebra isomorphism. Thus there is a canonical bijective
correspondence between the orderings $\omega$ of $k$ and the $k$-isomorphism
classes of real closures $k_\omega$ of $k$. 
\end{theorem}
We shall need one other basic fact about real fields:
\begin{lemma}
\label{real_closure}
If a field $k$ is real but not real closed then either there
exists an extension field of odd degree $q>1$ or, 
for some ordering of $k$, there exists a positive non-square element
$a\in k$ (or both). 
\end{lemma}
\begin{proof}
See~\cite[p9]{BCR98} or~\cite[p113]{Scha85}. 
\end{proof}

The main theorem of the chapter says that every 
element of infinite order in $W^\epsilon(R\dash k)$ is non-zero in some
$W^\epsilon(R\dash k_\omega)$:
\begin{theorem}
\label{pfister_generalized}
Let $k$ be a (commutative) field with trivial involution. \smallskip \\
a) If $k$ is not real then $2^nW^\epsilon(R\dash k)=0$ for some $n$.
\smallskip \\
b) Suppose $k$ is real and
$\{k_\omega\}_{\omega\in\Omega}$ is the family of real closures of
$k$. Then every element $v$ in the kernel of the canonical map
\begin{equation*}
W^\epsilon(R\dash k)\rightarrow \bigoplus_\omega W^\epsilon(R\dash k_\omega)
\end{equation*}
satisfies $2^nv=0$ for some $n$.
\smallskip \\
c) Suppose $K=k(\sqrt{-1})$ is a field extension of degree $2$ with
involution fixing precisely $k$. Then every
$v\in\Ker(W^\epsilon(R\dash k)\rightarrow W^\epsilon(R\dash K))$ satisfies $2v=0$.
\end{theorem}

In the particular case $k=\Q$ there is just one real closure, the set
of real algebraic numbers~$\Q_\omega=\R\cap \overline{\Q}$. 
Theorem~\ref{pfister_generalized} implies that the kernels of the natural maps
\begin{equation*}
W^\epsilon(R\dash \Q) \to W^\epsilon(R\dash \R\cap \overline{\Q}^-) \to W^\epsilon(R\dash \overline{\Q}^-)
\end{equation*}
are $2$-primary. Theorem~\ref{2primary_kernel} follows easily, because the natural map 
\begin{equation*}
W^\epsilon(R\dash \overline{\Q}^-) \cong \Z^{\oplus
\ssdM(R\dash \overline{\Q}^-)}\to \Z^{\ssdM(R\dash \C^-)} \cong
W^\epsilon(R\dash \C^-) 
\end{equation*}
is an injection. 

Note that when $R=\Z$ and $\epsilon=1$, the kernel of
$W^1(\Q)\to W^1(\C^-)$ already fails to be finitely generated -
Further discussion of invariants which
detect the torsion part of $W^\epsilon(R\dash \Q)$ can be found in
chapter~\ref{chapter:complete_invariants}.
The remainder of this chapter is devoted to the proof of
theorem~\ref{pfister_generalized}.
\section{Frobenius reciprocity}\index{Frobenius reciprocity} 
\label{section:frobenius}
\begin{notation}
Let $k$ denote any field with involution. For each
involution-invariant element $a=\overline{a}\in k$ there is a
one-dimensional (symmetric or) hermitian form 
\begin{equation*}
\phi:k\to k^*; x\mapsto (y\mapsto
xa\overline{y})
\end{equation*}
which is denoted $\langle a\rangle_k$ or simply
$\langle a\rangle$. 
Any (finite-dimensional) hermitian form over $k$
can be expressed as a direct sum of one-dimensional forms
(diagonalized). One writes
\begin{equation*}
\langle a_1, a_2, \cdots, a_m\rangle := \langle
a_1\rangle\oplus\langle a_2\rangle \oplus\cdots\langle a_m\rangle.
\end{equation*}
We use the same notation to denote the Witt class in $W^1(k)$
represented by a non-singular hermitian form.
\end{notation}
Let $K/k$ be an extension of commutative fields with involution.
We have seen in chapter~\ref{chapter:preliminaries},
section~\ref{section:Wittgroups_Changeofrings} that the inclusion
of $k$ in $K$ induces a homomorphism of Witt rings
$W^1(k)\rightarrow W^1(K)$ and a group
homomorphism $W^\epsilon(R\dash k)\rightarrow W^\epsilon(R\dash K)$. Each
of these maps is denoted $v\mapsto v_K$.

Suppose $s:K\rightarrow k$ is any $k$-linear map which respects the
involutions. Then $s$ induces group homomorphisms 
\begin{equation*}
W^\epsilon(K)\rightarrow W^\epsilon(k)\quad\mbox{and}\quad
W^\epsilon(R\dash K)\rightarrow W^\epsilon(R\dash k)
\end{equation*}
via the equation $s(M,\phi)=(M,s\phi)$ where 
\begin{equation*}
(s\phi)(m_1)(m_2)=s(\phi(m_1)(m_2))
\end{equation*}
for all $m_1,m_2\in M$. 
\begin{lemma}[Frobenius Reciprocity, Scharlau~\cite{Scha70}]
Let $u\in W^1(K)$ and $v\in W^\epsilon(R\dash k)$. Then
\begin{equation*}
s(u.v_K)=su.v
\end{equation*}
In particular, 
\begin{equation*}
s(v_K)=(s\langle1\rangle_K).v
\end{equation*}
\end{lemma}
\begin{proof}
Straightforward.
\end{proof}
\section{Proof of Theorem~\ref{pfister_generalized}}
\label{section:proof_pfister}
a) If $k$ is not formally real then $2^nW^1(k)=0$ for some positive
 integer $n$ (e.g.~see Milnor and Husemoller~\cite[p68 or p76]{MilHus73}). Since
$W^\epsilon(R\dash k)$ is a $W^\epsilon(k)$-module we have $2^nW^\epsilon(R\dash k)=0$. 
\smallskip \\ \noindent
b) Assume the contrary, that there exists $v\in W^\epsilon(R\dash k)$ such
that $2^n v\neq0$ for all $n\geq0$ and the image of $v$ is zero in
each group $W^\epsilon(R\dash k_\omega)$. By Zorn's lemma there exists a
 maximal algebraic extension $K/k$ 
such that $2^n v_K\neq0$ for all $n\geq0$.
Now $K$ is real by a) but is not real closed so 
lemma~\ref{real_closure} says that either there exists an extension
$K'=K(\xi)/K$ of odd degree $d>1$ or for some ordering of $K$ there is
a positive non-square $a\in K$.

In the former case, define $s:K'\rightarrow K$ by $s(1)=1$ and
$s(\xi^i)=0$ for $1\leq i\leq q-1$. Then 
$s\langle1\rangle_{K'}=\langle1\rangle_K\in W^1(K)$ because the elements
$\xi^\frac{q+1}{2},\cdots,\xi^{q-1}$ span a sublagrangian for
$s(1_{K'})$ - see Scharlau~\cite[p49]{Scha85} for further details.  
Frobenius reciprocity gives the equation $s(v_{K'})=v_K$.
It follows from the maximality of $K$ that
$2^n v_{K'}=0$ for some $n$ so $2^n v_K=s(2^nv_{K'})=0$ and we have reached a
contradiction.

The latter case is similar: Suppose $a\in K$ is non-square and
positive in some ordering so $-a$ is also non-square. Let
$K'=K(\sqrt{a})$ and define $s:K'\rightarrow K$ by $s(1)=1$ and
$s(\sqrt{a})=0$. Then $s\langle1\rangle_{K'}=
\langle1, a\rangle_K$ so Frobenius reciprocity gives
$s(v_{K'})=\langle1, a\rangle.v_K$. Hence
\begin{equation}
\label{1_a_is_2_primary}
2^n\langle1, a\rangle.v_K=0~~\mbox{for some $n$}.
\end{equation}
Now $-a$ is another non-square in $K$ so by the same argument 
\begin{equation}
\label{1_-a_is_2_primary}
2^n\langle1, -a\rangle.v_K=0~~\mbox{for some $n$}.
\end{equation}
For suitably large $n$ we may sum equations~(\ref{1_a_is_2_primary})
and~(\ref{1_-a_is_2_primary}) to find that
$2^{n+1}v_K=0\in W^\epsilon(R\dash K)$. Once again we have reached a
contradiction.
\smallskip \\
c) Let $v\in W^\epsilon(R\dash k)$ and suppose $v_K=0$. Define
$s:K\rightarrow k$ by $s(1)=1$ and $s(\sqrt{-1})=0$. Frobenius
reciprocity yields 
\begin{equation*}
0=s(v_K)=\langle1, 1\rangle.v=2v~.
\end{equation*}
\chapter{Characters}\index{Character|(textbf}
\label{chapter:characters}
Having associated cobordism invariants $\sigma_{M,b}$ to simple self-dual
representations of the ring~$P_\mu$ we turn now to the theory of
characters to distinguish such representations. A `trace' invariant
for boundary links was first introduced by Farber~\cite{Far92} 
under the assumption that the ground ring should be an algebraically
closed field. Some of his work was later simplified by Retakh, Reutenauer and
Vaintrob~\cite{RRV99}.

In the present chapter $k$ is any field of 
characteristic zero. All representations will be implicitly assumed to
be finite-dimensional.
The characters of simple
representations of an arbitrary associative ring $R$ over $k$ are
shown to be linearly independent which implies that a semisimple
representation is determined up to isomorphism by its character.   

Recalling from chapter~\ref{chapter:devissage} that $F_\mu$-link signatures
are associated to self-dual representations, it is shown at the
conclusion of the present chapter that a semisimple representation 
is self-dual if and only if the corresponding character
$\chi:R \to k$ respects the involutions of~$R$ and~$k$ (cf~\cite[\S6]{Far92}). 

The books by Serre~\cite{Ser77} and Curtis and Reiner~\cite{CurRei81} on
representation theory of finite groups are used for basic reference.
We do not assume the invariant theory of Procesi and Le Bruyn which we
employed in chapter~\ref{chapter:variety}. 

The linear independence of characters will be exploited again in the
study of rationality questions in
chapter~\ref{chapter:rationality_of_representations} below. 
\section{Artin algebras}
We must first recall the definition and basic theory of Artin algebras
\begin{definition}\index{Algebra!Artin}
A $k$-algebra~$S$ is called {\it Artinian} or an {\it Artin $k$-algebra}
if every descending chain of left ideals
\begin{equation*}
I_1 \supset I_2 \supset \cdots \supset I_n \supset \cdots
\end{equation*}
terminates, i.e.~eventually, $I_n=I_{n+1}=I_{n+2}=\cdots$~.
\end{definition}
\noindent For example, every finite-dimensional $k$-algebra is Artinian.
Recall that a ring~$S$ is by definition {\it simple} if~$0$ and~$S$ are
the only two-sided ideals. 
\begin{theorem}
\label{Artinian_structure}
i) An Artin $k$-algebra~$S$ is simple if and only if it admits a faithful
simple representation $\rho:S\hookrightarrow \End_k(M)$. \\ \noindent
ii) (Skolem-Noether) A simple Artin $k$-algebra admits a unique simple
representation. \\ \noindent 
ii) Every simple Artin algebra is isomorphic to a
matrix ring $M_n(D)$ over a division $k$-algebra $D$.
\end{theorem}
\section{Independence of Characters}\index{Character!independence of}
\begin{definition}
Suppose~$A$ is a commutative ring. The {\it character}
$\chi_M\in\Hom_\Z(R,A)$ of a representation $(M,\rho)$ of~$R$ over~$A$ is 
\begin{align*}
\chi_M:R &\to A \\
r &\mapsto \Trace(\rho(r)).
\end{align*}
\end{definition}
We concentrate here on the case $A=k$ where~$k$ denotes a field
of characteristic zero.

\begin{definition}\index{Representation!simple (=irreducible)|textbf}
A representation~$M$ over~$k$ is called {\it simple} or
{\it irreducible} if there are no subrepresentations other than $0$
and $M$. One says $M$ is {\it semisimple} if it is a direct sum of simple
representations. 
\end{definition}
\begin{proposition}
\label{independence_of_characters}
If  $(M_1,\rho_1)$, $(M_2,\rho_2)$,
$\cdots$, $(M_l,\rho_l)$ are (pairwise) non-isomorphic simple
representations then $\chi_{M_1}$, $\chi_{M_2}$, $\cdots$,
$\chi_{M_l}$ are linearly independent over $k$.
\end{proposition}
\begin{corollary}
\label{character_classification}
Let~$k$ be a (commutative) field of characteristic zero. 
Two semisimple representations $(M,\rho)$ and $(M',\rho')$ of~$R$ over
$k$ are isomorphic if and only if $\chi_M=\chi_{M'}$. 
\end{corollary}
We shall prove proposition~\ref{independence_of_characters} by
studying the kernels of the maps $\rho_i$.
Suppose $A$ is a ring, $M$ is an $A$-module and $T$ is an abelian group.
A homomorphism $f\in\Hom_\Z(T,M)$ determines, and is determined by, the
$A$-module map $Af : A\otimes_\Z T \to M$ given by 
$Af(a\otimes t)=af(t)$ for all $a\in A$ and $t\in T$. There is a
commutative diagram
\begin{equation*}
\xymatrix{
T \ar[r]^f \ar[d]^i & M \\
A\otimes_\Z T \ar[ur]_{Af} & }
\end{equation*}
where~$i$ is the natural map $t\mapsto 1\otimes t$.
In particular, a representation $\rho:R\to \End_k(M)$ and its character
$\chi_M:R\to k$ determine and are determined by $k$-linear maps
$k\rho: k\otimes_\Z R \to \End_k(M)$ and
$k\chi_M:k\otimes_\Z R \to k$ respectively.
\begin{lemma}
\label{maximal_ideals}
Let~$R$ be any ring. Two simple representations $(M,\rho)$ and
$(M',\rho')$ over~$k$ are isomorphic 
if and only if $\Ker(k\rho)=\Ker(k\rho')$.
\end{lemma}
\begin{proof}
If $\theta:(M,\rho)\to(M',\rho')$ is an isomorphism then
\begin{equation*}
k\rho(x)=\theta^{-1}k\rho'(x)\theta
\end{equation*}
 for all $x\in k\otimes_\Z R$ so $\Ker(k\rho)=\Ker(k\rho')$.

Conversely, suppose
$\Ker(k\rho)=\Ker(k\rho')=I$. The quotient $\frac{k\otimes_\Z R}{I}$
has faithful irreducible representations $(M,k\rho)$ and $(M',k\rho')$.
Since $\frac{k\otimes_\Z R}{I}$ is finite-dimensional over~$k$, and hence 
Artinian it follows from theorem~\ref{Artinian_structure} that
$(M,k\rho)\cong(M',k\rho')$, so $(M,\rho)\cong(M',\rho')$. 
\end{proof}
\begin{proof}[Proof of proposition~\ref{independence_of_characters}]
Let $I_i=\Ker(\rho_i)\subset k\otimes_\Z R$ for $1\leq i\leq l$. By
lemma~\ref{maximal_ideals} and theorem~\ref{Artinian_structure} i), $I_1, \cdots, I_l$ are distinct maximal
two-sided ideals of $k\otimes_\Z R$. In particular, $I_i+I_j=R$ for
$i\neq j$ so, fixing $i\in\{1,\cdots,l\}$,
\begin{equation*}
R=\prod_{j\mid j\neq i}(I_i+I_j) \subset I_i+\prod_{j\mid
j\neq i}I_j
\end{equation*}
and there is an equation $1=x'_i+x_i$ with $x'_i\in I_i$ and $x_i\in
\prod_{j\mid j\neq i} I_j \subset \bigcap_{j\mid j\neq i}I_j$. Thus
$x_i$ acts as the identity on $M_i$ and as zero on the other $M_j$
(cf~\cite[ex9 p170]{CurRei81}) so
$(k\chi_{M_i})(x_i)=\dim_kM_i$ and $(k\chi_{M_j})(x_i)=0$ when $j\neq i$.

If there is a linear relation $\sum_{j=1}^la_j\chi_{M_j}=0$ with each
$a_j\in k$ then, evaluating at $x_i$, we have
\begin{equation*}
0=\sum_{j=1}^la_j(k\chi_{M_j})(x_i)=a_i\dim_kM_i
\end{equation*}
so $a_i=0$ for each $i\in\{1,\cdots,l\}$. 
\end{proof} 
\subsection{Self-Dual Characters}
Let~$R$ be a ring with involution and let~$k$ be a field with
involution. Recall from example~\ref{examples_of_hermitian_categories}
the duality functor $(M,\rho)\mapsto(M^*\rho^*)$ on the category
$(R\dash k)\proj$.
\begin{lemma}
\label{detecting_self-dual}
A semisimple representation~$M$ of~$R$ over~$k$ is self dual $(M,\rho)\cong
(M^*,\rho^*)$ if and only if the character $\chi_M$ respects the
involutions, i.e.~$\chi_M(\overline r)=\overline{\chi_M(r)}$
for all~$r\in R$.
\end{lemma}
\begin{proof}
By corollary~\ref{character_classification},
$(M,\rho)\cong(M^*,\rho^*)$ if and only if $\chi_M=\chi_{M^*}$. Now
$\chi_{M^*}(\overline{r})=\Trace(\rho^*(\overline{r}))=\Trace(\rho(r)^*)=\overline{\Trace(\rho(r))}=\overline{\chi_M(r)}$
so $(M,\rho)\cong(M^*,\rho^*)$ if and only if $\chi_M(\overline
r)=\overline{\chi_M(r)}$ for all $r\in R$.
\end{proof}\index{Character|)textbf}
\chapter{Detecting Rationality and Integrality}
\label{chapter:rationality_of_representations}
As discussed in chapter~\ref{chapter:main_results}, and in more detail
in chapter~\ref{chapter:devissage}, the composite
\begin{equation*}
W^{\epsilon}(P_\mu\dash \Z)\to W^{\epsilon}(P_\mu\dash\C^-) \cong \Z^{\bigoplus\infty}
\end{equation*}
associates to a Seifert form one signature invariant $\sigma_{M,b}$ for
each self-dual simple complex representation~$M$ of the quiver
$P_\mu$. 

These are a complete set of signatures, but
some are surplus to requirements. One source of redundancy derives
from the fact that each signature is equal to its complex conjugate
$\sigma_{M,b}=\sigma_{\overline{M},\overline{b}}$. A second source is
that a necessary condition for a 
signature $\sigma_{M,b}$ to be non-zero is
that~$M$ should be {\it algebraically integral}, i.e.~$M$ should be a
summand of some representation induced up from an integral representation. 

The present chapter shows how the character can be used to identify the
algebraically integral representations among a variety of
complex representations.
The chapter is arranged in two sections dubbed rationality and integrality. 
In the rationality section we study an arbitrary field extension
$k_0\subset k$ in characteristic zero showing that a semisimple representation
$M$ of an arbitrary associative ring $R$ over~$k$ is a summand of a representation induced up from $k_0$
if and only if the character~$\chi_M$ takes values in a finite
 extension field of~$k_0$.  
 
Turning to questions of integrality, we prove that a
representation~$M$ over~$\C$ is a summand of an integral
representation if and only 
if $\chi_M$ takes values in the ring of algebraic integers of some finite
extension of~$\Q$. 
In particular, when $R$ is a path ring of a quiver, one need 
only check that the traces of oriented cycles lie in such an algebraic
number ring.

\section{Rationality}
In this section $k_0\subset k$ denotes a (possibly infinite) extension
of characteristic zero fields. 
\subsection{Induction}
If $M_0$ is a representation of~$R$ over $k_0$ then, as in
the third paragraph of chapter~\ref{chapter:preliminaries},
section~\ref{section:representations}, $k\otimes_{k_0}M$ is a
representation of~$R$ over $k$. A representation over~$k$ which is
isomorphic to
$k\otimes_{k_0}M_0$ for some $M_0$ is called a {\it $k_0$-induced
representation}, or simply an {\it induced
representation}\index{Representation!induced} when the identity of
$k_0$ is clear.
\begin{lemma}[Induction Lemma]
\label{induction_lemma}
Suppose $(M_0,\rho_0)$ and $(M'_0,\rho'_0)$ are representations of~$R$
over $k_0$. \\ \noindent
i) $M_0$ is semisimple if and only if $k\otimes_{k_0}M_0$ is
semisimple. \smallskip \\ \noindent 
ii) If $M_0$ and $M'_0$ are semisimple then 
\begin{equation*}
M_0\cong M'_0~\Longleftrightarrow~k\otimes_{k_0}M_0\cong
k\otimes_{k_0}M'_0.
\end{equation*}
iii) If $M_0$ and $M'_0$ are non-isomorphic simple representations, no
summand of $k\otimes_{k_0}M_0$ is isomorphic to a summand of
$k\otimes_{k_0}M'_0$.
\end{lemma}
\begin{proof}
i) See for example Curtis and Reiner~\cite[3.56(iii), 7.5]{CurRei81} or
Bourbaki~\cite[Theorem 2 p87]{Bou58}. \\ \noindent
ii) Follows from proposition~\ref{character_classification} since
$\chi_{k\otimes_{k_0}M_0}=\chi_{M_0}\in\Hom_\Z(R,k_0)$. \\ \noindent 
iii) As in the proof of proposition \ref{character_classification}
there exists an element $x\in k_0\otimes_\Z R$ such that~$x$ acts as
the identity on $M_0$ and as zero on $M'_0$. 
Now~$x$ also acts as the
identity on $k\otimes_{k_0}M_0$ and as zero on $k\otimes_{k_0}M'_0$.
\end{proof}

In the following lemma, we assume that $R$, $k_0$ and $k$ are endowed
with involutions.
\begin{lemma}
\label{sd_induction_lemma}
Suppose $M_0$ is simple and
$M$ is a simple summand of the induced representation
$k\otimes_{k_0}M_0$. If~$M$ is self-dual then $M_0$ is also self-dual.
\end{lemma}
\begin{proof}
Observing that $M^*$ is a summand of $k\otimes_{k_0}M_0^*$, if
$M\cong M^*$ then, lemma~\ref{induction_lemma} iii) implies that $M_0\cong
M_0^*$. 
\end{proof}
The converse to lemma~\ref{sd_induction_lemma} is not true in
general; a self-dual simple representation $M_0$ over $k_0$ decomposes over
$k$ into summands among which there may be dual pairs~$M\oplus M^*$.

\subsection{Restriction}
If~$k$ is a finite-dimensional extension field of $k_0$ and
$(M,\rho)$ is a representation of~$R$ over~$k$ then~$M$ can be
regarded as a representation over $k_0$:
\begin{definition}\index{Representation!restriction of}
The {\it restriction} $\res^k_{k_0}(M,\rho)$ of~$M$ to $k_0$ is the
representation $(\res^k_{k_0}M,F\circ\rho)$ where $\res^k_{k_0}M$ is
$M$ regarded as a vector space over $k_0$ and $F:\End_kM\to\End_{k_0}M$ is the forgetful map.
\end{definition}
\begin{lemma}[Restriction Lemma]
\label{restriction_lemma}
Suppose~$k$ is a finite-dimensional extension of $k_0$ and~$M$
is a representation of~$R$ over $k$. Then~$M$ is semisimple if and only 
if the restriction $\res_{k_0}^kM$ is semisimple.
\end{lemma}
\begin{proof}
Suppose first that~$M$ is simple and let $N_0$ be any simple submodule of
$\res_{k_0}^kM$. If $a_1, \cdots, a_m$ is a basis for~$k$
over $k_0$ then for each $i$, $a_iN_0$ is a simple representation of
$R$ over $k_0$ and $\sum_{i=1}^m a_iN_0=\res_{k_0}^kM$. Hence
$\res_{k_0}^kM$ is semisimple. It follows immediately that if~$M$ is
semisimple then $\res_{k_0}^kM$ is semisimple.

Conversely, suppose $\res_{k_0}^kM$ is semisimple. There is a natural
surjection $k\otimes_{k_0}\res_{k_0}^kM \twoheadrightarrow M$ and by
lemma~\ref{induction_lemma}~i) above $k\otimes_{k_0}\res_{k_0}^kM$ is
semisimple. 
\end{proof}
\subsection{Criteria for Rationality}
\begin{proposition}
\label{character_in_subfield}
Suppose $k_0\subset k$ and 
$(M,\rho)$ is a semisimple representation of~$R$ over $k$. The
following are equivalent: 
\begin{enumerate}
\item There exists a positive integer~$d$ and a semisimple
representation $M_0=(M_0,\rho_0)$ such that $M^{\oplus d}\cong
k\otimes_{k_0}M_0$. 
\item $\det(x-\rho(r))\in k_0[x]$ for all $r\in R$.
\item The character $\chi_M$ is $k_0$-valued, i.e.~$\Trace(\rho(r))\in
k_0$ for all $r\in R$. 
\end{enumerate}
\end{proposition}
\begin{remark}
\label{induction_uniqueness}
Case a) of the proof below, together with lemma~\ref{induction_lemma} iii),
demonstrates that if~$M$ is a simple representation over~$k$
and $\chi_M$ is $k_0$-valued then $M_0$ can also be chosen, in a
unique way, to be a simple representation. The positive integer
$d=d(M,k/k_0)$ is then the `relative Schur index'.
\end{remark}
\begin{proof}[Proof of Proposition~\ref{character_in_subfield}]
$2\Leftrightarrow 3$: 
Let $f=\det(1-x\rho(r))\in k[x]$ be the
`reverse characteristic polynomial' of $\rho(r)$ which has the property $f\in
k_0[x]$ if and only if $\det(x-\rho(r))\in k_0[x]$.
The exponential trace formula 
\begin{equation*}
\left(-\frac{d}{dx}\log f = \right)-f^{-1}\frac{df}{dx}=\sum_{n\geq1}\trace(\rho(r^n))x^{n-1}
\end{equation*}
in the ring $k[[x]]$ of formal power series implies that
$f$ has coefficients in $k_0$ if and only if
$\trace(\rho(r^i))\in k_0$ for all $i\geq1$. \smallskip \\ 
$1\Rightarrow3$: $\chi_M=\frac{1}{d}\chi_{M^{\oplus
d}}=\frac{1}{d}\chi_{k\otimes_{k_0}M_0}=\frac{1}{d}\chi_{M_0}$  so
$\chi_M$ is $k_0$-valued. \smallskip \\ \noindent
The rest of the proof deduces statement~$1$ from statements~$2$ and
$3$. We proceed in stages assuming that a)~$M$ is simple and b)~$M$ is
semisimple. In case b) we first consider finite extensions $k/k_0$ before
addressing arbitrary field extensions (in characteristic
zero). \smallskip \\ \noindent
a) Suppose that $\chi_M$ is $k_0$-valued and~$M$ is simple. Let~$R$ act
on $\End_kM$ by left multiplication $r.\alpha=\rho(r)\alpha$ so that,
as a representation,
$\End_kM$ is isomorphic to $M^{\oplus m}$ where $m=\dim_kM$. Similarly,
$k\rho(R)$ and $k_0\rho(R)$ can be regarded as representations of~$R$
over~$k$ and $k_0$ respectively.

If we can prove that the natural surjection
\begin{equation}
\label{nat_on_induction}
k\otimes_{k_0}k_0\rho(R) \to k\rho(R) \subset \End_kM\cong M^{\oplus n}
\end{equation}
is an isomorphism, i.e.~that $\dim_{k_0}k_0\rho(R)=\dim_kk\rho(R)$,
then any simple subrepresentation $M_0$ of $k_0\rho(R)$ has the
property $k\otimes_{k_0}M_0\cong M^{\oplus d}$ for some $d$.

To show that (\ref{nat_on_induction}) is an isomorphism, 
let $S\subset R$ be a finite subset with the property
that $\rho(S)$ is a basis for $k\rho(R)$ over $k$. 
Note that~$M$ is a faithful simple representation of
$k\rho(R)$ so $k\rho(R)$ is a simple
$k$-algebra by theorem~\ref{Artinian_structure} i). It follows that 
\begin{align*}
\Phi:k\rho(R) &\to k^S \\
\alpha &\mapsto (s\mapsto \trace(\rho(s)\alpha))
\end{align*}
is injective, for if $\trace(\rho(s)\alpha)=0$ for all $s\in S$ then
$\trace(k\rho(R)\alpha)=0$ which implies that
$\trace(k\rho(R)\alpha k\rho(R))=0$
and hence that $\alpha=0$. Since $\dim(k\rho(R))=\dim(k^S)$, $\Phi$ is an isomorphism of
$k$-vector spaces and the restriction $\Phi|:k_0\rho(R)\to k_0^S$
is an isomorphism of $k_0$-vector spaces. It follows
that $\dim_{k_0}k_0\rho(R)=|S|=\dim_kk\rho(R)$
so~(\ref{nat_on_induction}) is an isomorphism as required. \smallskip
\\ \noindent
b) Suppose~$k$ is finite-dimensional over $k_0$.  Let
$M_0=\res^k_{k_0}M$ (cf~Serre~\cite[lemma 12 p92]{Ser77}).
Now 
\begin{equation*}
\chi_{k\otimes_{k_0}M_0}=\chi_{M_0}=\Tr_{k/k_0}(\chi_M)=[k:k_0]\chi_M
\end{equation*}
since $\chi_M$ is $k_0$-valued. By
proposition~\ref{character_classification} $k\otimes_{k_0}M_0\cong
M^{\oplus [k:k_0]}$.

If $k/{k_0}$ is an arbitrary extension, the following lemma reduces
the problem to the finite case:
\begin{lemma}
\label{finite_character_sum}
Let $k/k_0$ be an arbitrary field extension (in characteristic
zero). Suppose $\chi_1$, $\chi_2$, $\cdots$, $\chi_l$ are characters of
simple representations over $k$, and $\sum_{i=1}^l t_i\chi_i$ is
$k_0$-valued with each $t_i\in\Z$. Then there 
exists a finite extension $k_1/k_0$ such that $\chi_i$ is
$k_1$-valued for $1\leq i\leq l$. 
\end{lemma}
Suppose~$M$ is a semisimple representation over~$k$ and $\chi_M$ is
$k_0$-valued. Let us deduce from lemma~\ref{finite_character_sum} that
$M^{\oplus d}\cong k\otimes_{k_0}M_0$.
Writing~$M$ as a direct sum of simple
representations
\begin{equation*}
M=(M^{(1)})^{\oplus
t_1} \oplus \cdots \oplus (M^{(l)})^{\oplus t_l}
\end{equation*}
we obtain an equation $\chi_M=\sum_{i=1}^lt_i\chi_{M^{(i)}}$.
By lemma~\ref{finite_character_sum}, all the simple characters
$\chi_{M^{(i)}}$ take values in some finite extension $k_1$ of $k_0$ and   
by case a) above, there exists a positive integer
$d^{(i)}$ and a representation $M^{(i)}_1$ over $k_1$ such that
$(M^{(i)})^{\oplus d^{(i)}}\cong k\otimes_{k_1} M^{(i)}_1$. 
Writing $d_1={\rm lcm}(d^{(1)},\cdots, d^{(l)})$ we have $M^{\oplus d_1}\cong
k\otimes_{k_1}M_1$ for some semisimple representation $M_1$ over
$k_1$.
As in b) above  $M_1^{\oplus [k_1:k_0]}\cong k_1\otimes_{k_0}M_0$ where
$M_0=\Res^{k_1}_{k_0}M_1$ so $M^{\oplus d_1[k_1:k_0]} \cong k\otimes_{k_0}
M_0$.
\begin{proof}[Proof of Lemma~\ref{finite_character_sum}]
The proof is divided into three stages, which apply to increasingly
general classes of extensions $k/{k_0}$. We consider
i) Galois extensions; 
ii) Algebraic extensions; 
iii) Arbitrary extensions. \smallskip \\ \noindent
i) The Galois group $G=\Gal(k/k_0)$ acts on $\Hom_\Z(R,k)$ by the equation
$\chi^g(r)=g\chi(r)$ where $g\in G$, $\chi\in\Hom_\Z(R,k)$ and $r\in
R$.  

If $\chi=t_1\chi_1+\cdots+t_l\chi_l$ is $k_0$-valued then 
\begin{equation*}
t_1\chi_1+\cdots +t_l\chi_l=\chi=\chi^g=t_1\chi_1^g +
\cdots t_l\chi_l^g
\end{equation*}
so, by the linear independence of characters
(proposition~\ref{independence_of_characters}),~$G$ permutes the set 
$\{\chi_1,\cdots, \chi_l\}$.

The kernel of this action is a normal subgroup $H\leq G$ of finite
index. Let $k_1=k^H$ be the intermediate field $k_0\subset k_1\subset k$
of elements fixed by $H$. 
Each character $\chi_i$ is fixed by~$H$ and is therefore $k_1$-valued so
it remains to note that $[k_1:k_0]\leq [G:H]<\infty$ by the
following standard Galois theory argument: 
If $\alpha\in k_1$ then the orbit $G.\alpha\subset k_1$ is
finite with cardinality at most $[G:H]$, so $\alpha$ is a root of a
polynomial $\prod_{\lambda\in G\alpha}(x-\lambda)\in k_0[x]$ of
degree at most $[G:H]$. It follows by the primitive element
theorem~\ref{bog_standard_primitive_element_theorem}
that $k_1$ is a finite extension of $k_0$ with $[k_1:k_0]\leq [G:H]$. 
(In fact $[k_1:k_0]=[G:\overline H]$ where $\overline H$ is the closure
of~$H$ in the Krull topology - see for example~\cite[Theorem
2.11.3]{RibZal00}). \smallskip \\ \noindent
ii) Given an algebraic extension $k/k_0$, let $K/k_0$ be a Galois
extension such that $K\supset k$ (e.g.~let~$K$ be the algebraic closure
$\overline{k_0}$). By case i) there exists a finite field extension
$K_1$ of $k_0$ such that
$k_0\subset K_1\subset K$ and each $\chi_i$ is
$K_1$-valued. Since $\chi_i$ is certainly $k$-valued we may set
$k_1=k\cap K_1$. \smallskip \\ \noindent 
iii) Arguing as in ii) we can assume that~$k$ is algebraically
closed and therefore contains the algebraic closure $\overline{k_0}$
of $k_0$. 
Let $(M^{(i)},\rho^{(i)})$ be a simple representation over $k$ with
character $\chi_i$ and let $M=(M^{(1)})^{\oplus
t_1} \oplus \cdots \oplus (M^{(l)})^{\oplus t_l}$. Suppose $\chi_M$ is
$k_0$-valued. For each $r\in R$ the characteristic polynomial of the
endomorphism $\rho^{(1)}(r)^{\oplus t_1} \oplus \cdots \oplus
\rho^{(l)}(r)^{\oplus t_l}$ is in $k_0[x]$ so the eigenvalues are in $\overline{k_0}$. It follows that the eigenvalues of each $\rho^{(i)}(r)$ 
are in $\overline{k_0}$ so $\chi_i$ is $\overline{k_0}$-valued for all
$i$ and we have reduced the problem to case i).
\end{proof}
This completes the proof of proposition~\ref{character_in_subfield}.
\end{proof}
\begin{remark} 
For the purposes of our application to boundary links one could
add the hypothesis that~$R$ is finitely generated, which
simplifies the proof of lemma~\ref{finite_character_sum} as follows: 
If $k/{k_0}$ is an algebraic extension and $\chi$ is the character of
a representation~$M$ over~$k$ then, fixing a basis for $M$, a finite
generating set for~$R$ acts via 
(finitely many) square matrices, the entries in which generate a
finite extension field $k_1$ of $k_0$. Plainly $\chi$ is $k_1$-valued.
\end{remark} 
\begin{corollary}
\label{character_in_fin_ext}
Under the hypotheses of proposition~\ref{character_in_subfield}, the
following are equivalent: 
\begin{enumerate}
\item There exists a semisimple representation $(M_0,\rho_0)$ 
over $k_0$ such that~$M$ is a summand of $k\otimes_{k_0}M_0$. 
\item There exists a finite extension $k_1/k_0$, a positive integer
$d$ and a semisimple representation~$(M_1,\rho_1)$ over $k_1$ such
that $M^{\oplus d} \cong k\otimes_{k_1}M_1$.
\item There exists a finite extension $k_1/k_0$ such that
$\det(x-\rho(r))\in k_1[x]$ for all $r\in R$.
\item There exists a finite extension $k_1/k_0$ such that
the character $\chi_M$ is $k_1$-valued. 
\end{enumerate}
\end{corollary}
\begin{proof}
Statements $2$,~$3$ and~$4$ are equivalent by
proposition~\ref{character_in_subfield}. \\ \noindent
$1\Rightarrow3$: Suppose $M\oplus M'\cong k\otimes_{k_0}M_0$. Let
$\chi_1$, $\cdots$, $\chi_l$ denote the characters of the distinct
isomorphism classes of simple summands of $M\oplus M'$.
By lemma~\ref{finite_character_sum} above there is a finite
extension~$k_1$ of~$k_0$ such that all the characters $\chi_i$ are
$k_1$-valued. Hence $\chi_M$ is $k_1$-valued. \smallskip \\ \noindent
$2\Rightarrow1$: It suffices to show that there exists a semisimple
representation $M_0$ over $k_0$ such that $M_1$ is a summand of
$k_1\otimes_{k_0}M_0$. Set $M_0=\Res_{k_0}^{k_1}M_1$. By
lemma~\ref{restriction_lemma} $M_0$ is semisimple and there is a
natural surjection $k_1\otimes_{k_0}M_0\twoheadrightarrow M_1$.
\end{proof} 
\begin{definition}
A semisimple complex representation~$M$ of~$R$ is {\it
algebraic}\index{Representation!algebraic} if~$M$ is a
summand of a $\Q$-induced representation.
\end{definition}
\begin{corollary}
A semisimple complex representation~$M$ is algebraic if and only
if $\chi_M$ takes values in an algebraic number field. 
\end{corollary}
\begin{proof}
This is a special case of
corollary~\ref{character_in_fin_ext}.
\end{proof}
\begin{definition}
Two algebraic representations~$M$ and $M'$ are {\it conjugate}\index{Representation!conjugate} if
both~$M$ and $M'$ are summands of the same $\Q$-induced representation.
\end{definition}
\begin{corollary}
i) Two simple algebraic representations~$M$ and $M'$ are conjugate if
and only if $\chi_M$ and $\chi_{M'}$ lie in the same $\Gal(\overline
\Q/\Q)$-orbit. \smallskip \\ \noindent
ii) If~$R$ is the path ring of a quiver $Q$, then
conjugate simple algebraic representations have the same dimension vector. 
\end{corollary}
\begin{proof}
i) Let $G=\Gal(\overline{\Q}/\Q)$.
If~$O$ denotes the $G$-orbit of $\chi_M$ then $\sum_{\chi\in O}\chi$
is $G$-invariant and therefore $\Q$-valued. Note that $O$ is finite by
corollary~\ref{character_in_fin_ext}. 
By proposition~\ref{character_in_subfield}, there is a positive integer~$d$
such that $d\sum_{\chi\in O}\chi$ is the character of a rational
representation $M_0$. Every rational subrepresentation has
a $G$-invariant character which, by
proposition~\ref{independence_of_characters}, must be $d'\sum_{\chi\in
O}\chi$ for some integer $d'\leq d$.
If we choose~$d$ to be minimal then $M_0$ is
simple 
Lemma~\ref{induction_lemma} iii) implies that
$M'$ is conjugate to~$M$ if and only if $\chi_{M'}\in O$. \smallskip
\\ \noindent 
ii) Recall that $e_x$ denotes the idempotent in $\Z Q$ corresponding
to the trivial path at $x\in Q_0$. 
After i) it suffices to observe that for any semisimple complex
representation~$M$ of~$Q$ and any vertex $x\in Q_0$
\begin{equation*}
\dim_\C (M_x) = \chi(e_x) \in \Z_{\geq0}
\end{equation*} 
and is invariant under $\Gal(\overline \Q/\Q)$.
\end{proof}
\section{Integrality}
In place of a field extension, we consider in this section the
inclusion of~$\Z$ in~$\Q$. 
The general theory of integral representations is well-known to be
far more subtle. In particular, $(R\dash \Z)\proj$ is not an abelian
category. Nevertheless, one can give criteria for a rational
representation to be induced from some integral representation:
\begin{proposition}
\label{character_in_Z}
Suppose $(M,\rho)$ is a semisimple representation of~$R$ over~$\Q$.
The following are equivalent: 
\begin{enumerate}
\item There exists a representation $(M_0,\rho_0)$ of~$R$ over $\Z$
with the property that $M\cong \Q\otimes_\Z M_0$. 
\item $\det(x-\rho(r))\in \Z[x]$ for all $r\in R$.
\item The character $\chi_M$ is $\Z$-valued.
\end{enumerate}
\end{proposition}
\begin{caveat}
In contrast to remark~\ref{induction_uniqueness}, even if one assumes
that~$M$ is simple, $M_0$ is not in general unique.
\end{caveat}
\begin{proof}[Proof of proposition~\ref{character_in_Z}]
$1\Rightarrow2$: Immediate. \\ \noindent 
$2\Rightarrow3$: $\Trace(\rho(r))$ is a coefficient of
$\det(x-\rho(r))$. \\ \noindent
$3\Rightarrow2$:  We assume $\chi_M$ is $\Z$-valued and aim to prove that
$f=\det(1-x\rho(r))$ is $\Z$-valued. Once again we use the exponential trace
formula
\begin{equation*}
-f^{-1}\frac{df}{dx}=\sum_{i\geq 1}\trace(\rho(r^i)x^{i-1}) \in \Z[[x]].
\end{equation*}
If $f=1+a_1x + a_2x^2 + \cdots + a_nx^n$, let~$b$ denote the smallest
positive integer such that $ba_i\in\Z$ for all $i$, i.e.~$b$ is the
least common multiple of the denominators of 
the $a_i$. We aim to prove that $b=1$. Now $f=1+g/b$ where
$g=ba_1x+\cdots +ba_nx^n\in A[x]$ and $\hcf(b,g)=1$ so
\begin{equation*}
f^{-1}=1-\frac{g}{b} + \frac{g^2}{b^2} - \frac{g^3}{b^3} + \cdots 
\end{equation*}
If~$b$ has a prime factor~$p$ then plainly $f^{-1}\frac{df}{dx} \notin
\Z[[x]]$. Thus $b=1$ and $f\in \Z[x]$ as required. \smallskip \\ \noindent
$2,3\Rightarrow1$: 
Suppose first that~$M$ is a simple
representation so that $\Q\rho(R)$ is a simple algebra.
Choose a minimal finite subset $S\subset R$ such that $\rho(S)$ is a
$\Q$-basis for $\Q\rho(R)$. The vector space isomorphism
\begin{align*}
\Q\rho(R) &\to \Q^S \\
\alpha &\mapsto (s\mapsto \trace(\alpha\rho(s)))
\end{align*}
restricts to an injection $\rho(R)\to \Z^S$, so $\rho(R)$ is finitely
generated as a $\Z$-module. If $v_1, \cdots, v_n$ is a $\Q$-basis for
$M$ then we may set $M_0=\sum_{i=1}^n \rho(R)v_i$ and the natural
map $\Q\otimes_\Z M_0\to \Q M_0=M$ is an isomorphism.

To extend the result to the semisimple case note that a direct sum 
$(M,\rho)\oplus(M',\rho')$ satisfies condition 2 if and only if
$(M,\rho)$ and $(M',\rho')$ both satisfy condition 2; for 
\begin{equation*}
\det\left(x-(\rho(r)\oplus\rho'(r))\right)=\det\left(x-\rho(r)\right)\det\left(x-\rho'(r)\right) 
\end{equation*}
is in $\Z[x]$ if and only if $\det(x-\rho(r))\in\Z[x]$ and
$\det(x-\rho'(r))\in\Z[x]$.
\end{proof}
\begin{definition}
A complex representation~$M$ of~$R$ is {\it algebraically
integral}\index{Representation!algebraically integral|textbf} if
$M$ is a direct summand of a $\Z$-induced representation.
\end{definition}
Let $\cy{O}$\index{$O$@$\cy{O}$, $\cy{O}_K$} denote the ring of algebraic
integers. If~$K$ is an algebraic number field, then
$\cy{O}_K=\cy{O}\cap K$ denotes the ring of algebraic integers in $K$. 
\begin{corollary}
\label{character_in_number_ring}
Suppose $M=(M,\rho)$ is a semisimple representation of~$R$ over
$\C$. Then the following are equivalent:
\begin{enumerate}
\item~$M$ is algebraically integral.
\item There exists an algebraic number field~$K$
such that $\det(x-\rho(r))\in \cy{O}_K[x]$ for all $r\in R$. 
\item There exists an algebraic number field~$K$ such that $\chi_M$ is
$\cy{O}_K$-valued.  
\end{enumerate}
\end{corollary}
\begin{proof}
$1\Rightarrow2$: The eigenvalues of each endomorphism $\rho(r)$ are
algebraic integers so each characteristic polynomial $\det(x-\rho(r))$
has coefficients in $\cy{O}$.
Since~$M$ is a summand of $\C\otimes_\Q\Q\otimes_\Z M_0$,
corollary~\ref{character_in_fin_ext} implies that there is a finite extension
$K$ of $\Q$ such that each polynomial $\det(x-\rho(r))$ has
coefficients in~$\cy{O}_K$. \\ \noindent 
$2\Rightarrow3$: Trivial. \\ \noindent
$3\Rightarrow1$: Let $\Sigma$ be the set of embeddings $K\hookrightarrow\C$.
If $\chi_M$ is $\cy{O}_K$-valued then~$M$ is a summand of a
complex representation $M'$ with the $\Z$-valued character
$\bigoplus_{\sigma\in\Sigma}\chi_M^\sigma$.
By proposition~\ref{character_in_subfield}~$M$ is a summand of a 
$\Q$-induced representation $(M')^{\oplus d}$. 
By proposition~\ref{character_in_Z}, $(M')^{\oplus d}$ is $\Z$-induced.
\end{proof}
\begin{corollary}
\label{oriented_cycles_in_number_ring}
Let $M=(M,\rho)$ be a semisimple representation of a quiver~$Q$. The
following are equivalent:
\begin{enumerate}
\item $M$ is algebraically integral.
\item There is an algebraic number field $K$ such that the
characteristic polynomial of every oriented cycle in the quiver has
coefficients in $\cy{O}_K$, i.e.~$\det(x-\rho(r))\in \cy{O}_K[x]$.
\item There is an algebraic number field $K$ such
that~$\Trace(\rho(r))$ lies in~$\cy{O}_K$ for
every oriented cycle~$r\in Q$.
\end{enumerate}
\end{corollary}
\begin{proof}
It follows immediately from corollary~\ref{character_in_number_ring}
that $1\Rightarrow2$ and $2\Rightarrow3$.

To show that $3\Rightarrow1$, observe that if $r\in Q$ is a
non-trivial path but not a cycle then $\Trace(\rho(r))=0$. 
Thus condition 3 of corollary~\ref{oriented_cycles_in_number_ring}
implies condition 3~of corollary~\ref{character_in_number_ring}.
\end{proof}
%
%
\chapter[Example Varieties]{Representation Varieties: Two Examples}
\label{chapter:smalldim}
In this chapter we describe the
varieties of isomorphism classes of semisimple representations of
$P_2$ for dimension vectors
$\alpha=(1,1)$ and $\alpha=(1,2)$. In each case 
we compute the open subvariety of simple
representations, the real variety of self-dual
representations and the subset of algebraically integral
representations. I am grateful to Raf Bocklandt for a helpful
discussion of the case $\alpha=(1,2)$.

By theorem~\ref{oriented_cycles_generate}
the coordinate ring~$\C[X]^{\GL(\alpha)}$ for the variety of
semisimple isomorphism classes with dimension vector~$\alpha$ is generated
by a finite number of traces of
oriented cycles in the quiver $P_\mu$. Moreover the relations between
these generators can all be deduced from `Cayley-Hamilton relations'
(for further explanation see Procesi~\cite{Pro87}). 
While explicit sets of generators and relations 
can in principle be computed for any dimension vector~$\alpha$ a
general formula is not available at the time of writing.  
\medskip \\
A particular representation $(M,\rho)$ of $P_2$ can be displayed as follows:
\begin{equation}
\label{matrix_repn_of_P_2}
\left(
\begin{array}{c|c}
a_{11} & a_{12} \\
\hline
a_{21} & a_{22}
\end{array}
\right)
~=~~~
\xymatrix{
\bullet \ar@/^/[r]^{a_{21}} \ar@(ul,dl)[]_{a_{11}} & \bullet
\ar@/^/[l]^{a_{12}} \ar@(ur,dr)[]^{a_{22}} 
}
\end{equation}
If $e_{ij}$ denotes the arrow from vertex $j$ to vertex $i$ in the
quiver $P_\mu$ then $a_{ij}$ is by definition $\rho(e_{ij}):\pi_jM\to \pi_iM$.

\section{Dimension~$(1,1)$}
\begin{proposition}
\label{the_case_11}
The isomorphism classes of algebraically integral self-dual simple
representations of~$P_2$ with dimension vector $(1,1)$ correspond to triples of algebraic integers 
\begin{equation*}
\left(\frac{1}{2}+r_1i~,~\frac{1}{2}+r_2i~,~r_3\right)
\end{equation*}
with $r_1,r_2,r_3\in\R$ and $r_3\neq0$.
\end{proposition}
\subsection{Semisimple and Simple Representations}
\label{section:dimension_1_1}
If $\alpha=(1,1)$ then each $a_{ij}$ is a complex number and
the representation space 
\begin{equation*}
R(P_2,(1,1))=\bigoplus_{1\leq i,j\leq2}\Hom(\C,\C)
\end{equation*}
is $4$-dimensional affine
space. The coordinate ring is denoted
\begin{equation*}
\C[X]=\C[X_{11},X_{12},X_{21},X_{22}]
\end{equation*}
and the ring of invariant
polynomials for the conjugation action of $\GL(1,1)\cong
\C^{\,\bullet}\times\C^{\,\bullet}$ is generated by three algebraically
independent (traces of) oriented cycles:
\begin{equation*}
\C[X]^{\GL(1,1)}=\C[x_1,x_2,x_3]
\end{equation*}
with $x_1=X_{11}$, $x_2=X_{22}$ and~$x_3=X_{12}X_{21}$. 
The variety of semisimple dimension vector~$(1,1)$ isomorphism classes
of representations is
therefore $3$-dimensional affine space. A
representation~(\ref{matrix_repn_of_P_2}) is simple if 
and only if the complex number $a_{12}a_{21}$ is non-zero, so the
open subvariety of simple representations is defined by $x_3\neq 0$.

\subsection{Self-Dual Representations}
A representation is
self-dual if and only if its character $\chi:P_2\to \C^-$ respects the
involutions (lemma~\ref{detecting_self-dual}), i.e. 
\begin{equation*}
a_{11}=1-\overline{a_{11}},\quad
a_{22}=1-\overline{a_{22}}\quad\mbox{and}\quad
a_{12}a_{21}=\overline{a_{12}a_{21}}.
\end{equation*}
Thus the isomorphism classes of self-dual semisimple representations
correspond to the points in a
$3$-dimensional real affine space
\begin{equation}
\label{real_variety(1_1)}
\overline{\cy{M}}(P_2,(1,1))=
L\times L\times\R
\end{equation}
where $L=\left\{\frac{1}{2}+bi \mid b\in
\R\right\}\cong \R$.

Let us rewrite this computation from a dual point of view.
The duality functor on $(P_2\dash \C^-)\proj$ induces an involution on 
$\C^-[x_1,x_2,x_3]$ given by
\begin{equation*}
x_1\mapsto 1-x_1; \hspace{5mm} x_2\mapsto 1-x_2;\hspace{5mm}
x_3\mapsto x_3
\end{equation*}
(see chapter~\ref{chapter:variety},
section~\ref{section:variety_with_involution}). 
The elements $i(x_1-\frac{1}{2})$,
$i(x_2-\frac{1}{2})$ and~$x_3$ are plainly involution invariant and the
natural map  
\begin{equation*}
\C^-\otimes
\R\left[i\left(x_1-\frac{1}{2}\right)~,~i\left(x_2-\frac{1}{2}\right)~,~x_3\right]
\to \C^-[x_1,x_2,x_3]
\end{equation*}
is an isomorphism so the
involution-invariant subring is
\begin{equation*}
\C^-[X]^{\GL(1,1)\rtimes\frac{\Z}{2\Z}}
=\R\left[i\left(x_1-\frac{1}{2}\right)~,~i\left(x_2-\frac{1}{2}\right)~,~x_3\right]~.
\end{equation*} 
Thus real points correspond to maximal ideals
\begin{equation*}
\left(i\left(x_1-\frac{1}{2}\right) -
r_1~,~i\left(x_2-\frac{1}{2}\right)-r_2~,~x_3-r_3\right)
\end{equation*}
with $r_1,r_2,r_3\in\R$. The maximal ideals of
$\C^-[x_1,x_2,x_3]$ lying over these reals ideals correspond precisely
to the elements of $L\times L\times\R$ which
confirms~(\ref{real_variety(1_1)}) above.
\subsection{Algebraically Integral Representations}

Recall from chapter~\ref{chapter:rationality_of_representations} that
a complex representation is said 
to be algebraically integral if it is a summand of a representation
which is induced from an integral representation.
By corollary~\ref{oriented_cycles_in_number_ring}, a
semisimple complex representation of a quiver is algebraically
integral if and only if all the traces of oriented cycles 
evaluate to algebraic integers. When $\alpha=(1,1)$ the trace of an
oriented cycle is just some product of the generators $x_1$, $x_2$ and
$x_3$, so the algebraically integral representations correspond to
triples of algebraic integers in $L\times L\times\R$. This completes
the proof of proposition~\ref{the_case_11}.

\section{Dimension~$(2,1)$}
This section is devoted to the proof of the following proposition:
\begin{proposition}
\label{summary_dim_2_1}
The isomorphism classes of algebraically integral self-dual simple
representations of~$P_2$
with dimension vector~$(2,1)$ correspond to quintuples of algebraic integers
\begin{equation}
\left(1+2r_1i~,~r_2+r_1i~,~2r_3~,~r_3+r_4i~,~\frac{1}{2}+r_5i\right) 
\end{equation}
such that $r_1,\cdots,r_5\in\R$ and 
$r_3^2-2r_1r_4r_3+r_4^2\neq 4r_2r_3^2$.
\end{proposition}
\subsection{Semisimple Representations}
Fixing the dimension vector $\alpha=(2,1)$, the representation space $R(P_2,(2,1))$ is $9$-dimensional and has
coordinate ring $\C[X]$ where $X=\{X_{ij}\}_{1\leq i,j\leq 3}$. Let $Y_{ij}$ denote the matrix of indeterminates corresponding to the
linear map $a_{ij}$ so
\begin{equation*}
\left(\begin{array}{c|c}
Y_{11} & Y_{12} \\
\hline
Y_{21} & Y_{22}
\end{array}\right)
=
\left(
\begin{array}{cc|c}
X_{11} & X_{12} & X_{13} \\
X_{21} & X_{22} & X_{23} \\
\hline
X_{31} & X_{32} & X_{33}
\end{array}
\right).
\end{equation*}
\begin{proposition}
\label{dim_2_1_generators}
The invariant ring $\C[X]^{\GL(2,1)}$ is generated by the following
five polynomials:
\begin{align*}
z_1 &= \Trace(Y_{11}) = X_{11}+X_{22} \\ 
z_2 &= \Det(Y_{11}) = X_{11}X_{22}-X_{12}X_{21} \\
z_3 &= \Trace(Y_{12}Y_{21})= X_{13}X_{31}+X_{23}X_{32} \\
z_4 &= \Trace(Y_{11}Y_{12}Y_{21})= 
\left(\begin{matrix}
X_{31} & X_{32}
\end{matrix}\right)
\left(\begin{matrix}
X_{11} & X_{12} \\
X_{21} & X_{22}
\end{matrix}\right)
\left(\begin{matrix}
X_{13} \\
X_{23} 
\end{matrix}\right)
\\
z_5 &= \Trace(Y_{22})=Y_{22}=X_{33}
\end{align*}
Moreover these polynomials are algebraically independent which implies that
$\cy{M}(P_2,(2,1))$ is five-dimensional affine space.
\end{proposition}
\begin{proof}
It is easy to see that each $z_i$ is $\GL(\alpha)$-invariant,
i.e.
\begin{equation*}
\C[z_1,\cdots,z_5]\subset \C[X]^{\GL(\alpha)}.
\end{equation*}
To prove that $z_1$, $\cdots$, $z_5$ generate $\C[X]^{\GL(\alpha)}$
as a $\C$-algebra it suffices to show that the trace of every oriented
cycle lies in $\C[z_1,\cdots,z_5]$. In fact we prove slightly more:
\begin{lemma}
\label{traces_generated_over_Z}
The trace of every oriented cycle lies in $\Z[z_1,\cdots,z_5]$.
\end{lemma} 
\begin{proof}
Let $Z$ be a non-trivial oriented cycle.
$Y_{22}$ commutes with other oriented cycles at vertex~$2$ so we may
assume that there are no occurrences of $Y_{22}$ in $Z$ and aim to show
that $\Trace(Z)\in\Z[z_1,\cdots,z_4]$. Furthermore,
we may assume that the cycle $Z$ begins and ends at vertex $1$, if necessary by
applying some equation of the form
\begin{equation*}
\Trace(Y_{21}Z'Y_{12})=\Trace(Z'Y_{12}Y_{21}).
\end{equation*}

Let $c_1=Y_{11}$ and $c_2=Y_{12}Y_{21}$ denote the two basic cycles at
vertex $1$. An arbitrary cycle $Z$ is a word in the alphabet
$\{c_1,c_2\}$ and we have
\begin{equation*}
z_1=\Trace(c_1);\quad z_2=\Det(c_1);\quad z_3=\Trace(c_2);\quad
z_4=\Trace(c_1c_2).
\end{equation*}
Noting that $\det(c_2)=0$, the cycle $Z$ may be expressed in
terms of shorter cycles by means of the following Cayley-Hamilton identities:
\begin{align}
\label{2_1_Cayley-Hamilton}
c_1^2 &=\Trace(c_1)c_1 - \Det(c_1) = z_1c_1 - z_2 \\
c_2^2 &= \Trace(c_2)c_2 = z_3c_2 \\
c_1c_2c_1c_2 &= \Trace(c_1c_2)c_1c_2 = z_4c_1c_2.
\end{align}
Indeed, every sufficiently long word in $c_1$ and $c_2$ contains one of the
three subwords $c_1^2$, $c_2^2$ or $c_1c_2c_1c_2$; 
the only words which do not contain one of these three are the following: 
\begin{equation}
\label{smallcycles}
c_1,\: c_2,\: c_1c_2,\: c_2c_1,\: c_1c_2c_1,\: c_2c_1c_2.
\end{equation}
So $\Trace(Z)$ is a linear combination (with coefficients in
$\Z[z_1,\cdots,z_4]$) of the traces of these cycles~(\ref{smallcycles}).  
Now $\Trace(c_1c_2c_1)$ and $\Trace(c_2c_1c_2)$ each simplify further
because
\begin{equation*}
\Trace(c_1c_2c_1)=\Trace(c_1^2c_2)\quad\mbox{and}\quad 
\Trace(c_2c_1c_2)=\Trace(c_2^2c_1)
\end{equation*}
and $\Trace(c_2c_1)=\Trace(c_1c_2)=z_4$ so~$\Trace(Z)$ lies in
$\Z[z_1,\cdots, z_4]$.
\end{proof}

To prove the last sentence of proposition~\ref{dim_2_1_generators}, it
suffices to note that by remark~\ref{dimension_of_varieties} the
dimension of the variety $\cy{M}(P_2,(2,1))$ of isomorphism classes of
semisimple representations is $5$.
\end{proof}
\subsection{Simple Representations}
\label{section:simple_2_1}
\begin{proposition}
\label{2_1_open_subvariety}
The isomorphism classes of semisimple representations which are not
simple correspond to the points on the subvariety defined by
\begin{equation}
\label{2_1_nonsimple}
z_1z_3z_4-z_4^2 =z_2z_3^2.
\end{equation}
\end{proposition}
\begin{proof}
If $M$ is a semisimple representation of dimension~$(2,1)$ but is not
simple then there must be a (simple) summand of dimension $(1,0)$; for if $M$
decomposes as a direct sum of representations of dimension $(0,1)$ and $(2,0)$
respectively then the latter representation must decompose further. 
We may therefore express $M$ as a direct sum of two representations
whose dimensions are $(1,1)$ and $(1,0)$.

Our aim is to show that the image of the direct sum morphism
\begin{equation}
\label{direct_sum_morphism}
\operatorname{\oplus}:\cy{M}(P_2,(1,1))\times \cy{M}(P_2,(1,0)) \to \cy{M}(P_2,(2,1)).
\end{equation}
is the variety given by equation~(\ref{2_1_nonsimple}).

It follows from the definitions of the $z_i$ in
proposition~\ref{dim_2_1_generators}
that the morphism~\ref{direct_sum_morphism} is dual to the map of
coordinate rings
\begin{align*}
\theta:\C[z_1,\cdots,z_5] &\to \C[x_1,x_2,x_3,y] \\
z_1 &\mapsto x_1+y \\
z_2 &\mapsto x_1y \\
z_3 &\mapsto x_3 \\
z_4 &\mapsto x_1x_3 \\
z_5 &\mapsto x_2
\end{align*}
where $\C[x_1,x_2,x_3]$ is the coordinate ring of
$\cy{M}(P_2,(1,1))$ as in section~\ref{section:dimension_1_1} and 
$\C[y]$ denotes the coordinate ring of $\cy{M}(P_2,(1,0))$. 
We must show that the kernel of~$\theta$ is the ideal generated by 
$p=z_1z_3z_4-z_4^2 - z_2z_3^2$. 

It is straightforward to check that $p$ is
contained in $\Ker(\theta)$. To prove the converse note that 
the restriction of $\oplus$ to
$\cy{M}^s(P_2,(1,1))\times \cy{M}(P_2,(0,1))$ is injective and so the
image of $\oplus$ is four-dimensional and $\Ker(\theta)$ is a prime of
height one. Therefore we need only check that $p$ is
irreducible. Indeed, $p$ is linear in $z_1$ and the coefficient
$z_3z_4$ is coprime to $-z_4^2-z_2z_3^2$.
This completes the proof of proposition~\ref{2_1_open_subvariety}.
\end{proof}
\subsection{Self-Dual Representations}
As we discussed in chapter~\ref{chapter:variety},
section~\ref{section:variety_with_involution} the
duality functor on $\cy{M}(P_2,(2,1))$ induces the involution
\begin{equation*}
I(X_{ij})=
\begin{cases}
-X_{ji} &\text{if $j\neq i$} \\
1-X_{ii} &\text{if $j=i$}
\end{cases}  
\end{equation*}
on $\C^-[X]$. It follows from the definition of $z_1,\cdots, z_5$  (in
proposition~\ref{dim_2_1_generators}) 
that the restriction of this involution to the $\GL(2,1)$ invariant
ring $\C^-[X]^{\GL(2,1)}=\C^-[z_1,\cdots,z_5]$ is given by:
\begin{align*}
z_1 &\mapsto 2-z_1 \\
z_2 &\mapsto z_2-z_1+1 \\
z_3 &\mapsto z_3 \\
z_4 &\mapsto -z_4 + z_3 \\
z_5 &\mapsto 1-z_5.
\end{align*}

The elements $i(z_1-1)$, $2z_2-z_1$, $z_3$, $i(z_3-2z_4)$ and
$i\left(z_5-\frac{1}{2}\right)$ are fixed by the involution and the natural map
\begin{equation*}
\C^-\otimes_\R \R\left[i(z_1-1)~,~2z_2-z_1~,~z_3~,~i(z_3-2z_4)~,~i\left(z_5-\frac{1}{2}\right)\right]
~\to~\C^-[z_1,\cdots, z_5]
\end{equation*}
is an isomorphism so by lemma~\ref{invariantsubring}
\begin{equation*}
\C^-[X]^{\GL(2,1)\rtimes \frac{\Z}{2\Z}}~=~\R\left[i(z_1-1)~,~2z_2-z_1~,~z_3~,~i(z_3-2z_4)~,~i\left(z_5-\frac{1}{2}\right)\right]. 
\end{equation*}
Thus there is one real
maximal ideal and therefore one self-dual semisimple isomorphism
class of representations
for each quintuple
\begin{equation}
\label{self_dual_values}
\left(1+2r_1i~,~r_2+r_1i~,~2r_3~,~r_3+r_4i~,~\frac{1}{2}+r_5i\right) 
\end{equation}
with $r_1,\cdots,r_5\in\R$.
\subsection{Self-Dual Simple Representations}
In section~\ref{section:simple_2_1} above we showed that the
non-simple dimension $(2,1)$ isomorphism classes of representations
lie on the hypersurface defined by equation~(\ref{2_1_nonsimple}). We
must compute the intersection of this hypersurface with the real
variety of self-dual representations. Substituting the
quintuple~(\ref{self_dual_values}) for $(z_1,\cdots,z_5)$
in~(\ref{2_1_nonsimple}) and simplifying we obtain
\begin{equation}
\label{self-dual_2_1_nonsimple}
r_3^2-2r_1r_4r_3+r_4^2 = 4r_2r_3^2
\end{equation}
Self-dual simple representations correspond to
quintuples~(\ref{self_dual_values}) which do not satisfy
equation~(\ref{self-dual_2_1_nonsimple}) 
\subsection{Algebraically Integral Representations}
\begin{lemma}
A semisimple dimension~$(2,1)$ complex representation $M$ is
algebraically integral if and only if the invariants $z_1$, $\cdots$,
$z_5$ take algebraic integer values at $M$.
\end{lemma}
\begin{proof}
Apply lemmas~\ref{oriented_cycles_in_number_ring}
and~\ref{traces_generated_over_Z}. 
\end{proof}
This completes the proof of proposition~\ref{summary_dim_2_1}.
\chapter{Number Theory Invariants}
\label{chapter:complete_invariants}
In this chapter we use the theory of
symmetric and hermitian forms over division algebras to obtain
invariants which distinguish elements of finite order in 
$W^\epsilon(P_\mu\dash \Q)$ and hence, by the following lemma, to
distinguish such elements in $W^\epsilon(P_\mu\dash \Z)$.
\begin{lemma}
\label{rationals_injection}
Given any ring $R$ with involution, the natural map
\begin{equation*}
W^\epsilon(R\dash \Z)\to W^\epsilon(R\dash \Q)
\end{equation*}
is injective.
\end{lemma}
\begin{proof}
Suppose $(M,\rho,\phi)\in H^\epsilon(R\dash \Z)$ is a non-singular
$\epsilon$-hermitian form which represents the zero class in
$W^\epsilon(R\dash \Q)$. By lemma~\ref{lemma:unstable}, the induced form
 $\Q\otimes_\Z (M,\rho,\phi)$ is metabolic with metabolizer
$L\subset\Q\otimes M$ say. It follows that $L\cap M$ metabolizes
$(M,\rho,\phi)$.
\end{proof}
\noindent We enhance lemma~\ref{rationals_injection} in
section~\ref{section:localization_exact_sequence} below.

Devissage followed by hermitian Morita equivalence decompose
$W^\epsilon(P_\mu\dash \Q)$
as in corollary~\ref{general_devissage_and_ME}
\begin{equation*}
W^\epsilon(P_\mu\dash \Q) \cong \bigoplus_M W^1(\End_{(P_\mu\dash \Q)}M)
\end{equation*}
with one summand for each isomorphism class of $\epsilon$-self-dual
simple representations $M$ of 
$P_\mu$ over $\Q$. 
Each endomorphism ring $\End(M)$ is a finite $\Q$-dimensional
division algebra with involution. Sections~\ref{section:division_algebras_overQ} and~\ref{section:number_theorm_Witt_invariants} of the
present chapter are therefore a summary of results in the Witt theory
of such division algebras; fortunately, complete Witt invariants are
available.
We do not provide proofs in these sections but the relevant theory
can be found in the books of Albert~\cite{Alb39}, Lam~\cite{Lam73} and
Scharlau~\cite[Chapters 6,8,10]{Scha85} and in papers of
Lewis~\cite{Lew82',Lew82'',Lew82}. Closely
related $L$-theory computations have also been performed by
I.Hambleton and I.Madsen~\cite{HamMad93}.

There are five distinct classes of algebras with involution to be
considered; in four of these 
a local-global principle applies. In chapter~\ref{chapter:endrings}
below we prove that if $\mu\geq2$ {\it every} finite $\Q$-dimensional
division algebra with involution is the endomorphism ring of some
simple (integral) representation of $P_\mu$ over $\Q$, so all the five
classes are germane. Recall, however, that in the case $\mu=1$ of knot
theory, or in the case of a split $F_\mu$-link, one need only consider
algebraic number fields with non-trivial involution, class 2a in
theorem~\ref{table_of_invariants}. \vspace{1ex} 

It is natural to ask whether a complete set of torsion invariants of
$W^\epsilon(R\dash\Q)$ can be defined by completing $\Q$ at finite
primes, since all the torsion-free invariants can be obtained by
completing $\Q$ at the (unique) real prime (see
part b) of theorem~\ref{pfister_generalized}).
The answer is negative; since there is a class of
division algebra with involution whose Witt group does not satisfy a
local-global principle it follows from Morita equivalence that the
`Hasse-Minkowski' map 
\begin{equation*}\index{Hasse-Minskowski map}
W^\epsilon(P_\mu\dash \Q)\to \prod_p W^\epsilon(P_\mu\dash \Q_p)
\end{equation*}
is not injective - in fact the kernel is isomorphic to
$\left(\frac{\Z}{2\Z}\right)^{\oplus\infty}$. Nonetheless, the vast
majority of invariants of $W^\epsilon(P_\mu\dash \Q)$ can be defined locally.
\section{Division Algebras over $\Q$}
\label{section:division_algebras_overQ}
The basic structure theorem (Albert~\cite[p149]{Alb39}) is that
every finite-dimensional division algebra (indeed every central simple
algebra) over an algebraic number field is a cyclic
algebra.\index{Algebra!cyclic} 
A cyclic algebra $E$ over a field $K$ may be defined by the following
construction (see~\cite[p316-320]{Scha85} for example):

Let $L/K$ be a finite Galois extension, say of degree $m$, with
cyclic Galois group $\Gal(L/K)$. Let $\sigma$ be a generator of
$\Gal(L/K)$ and let $a_0\in K^{\,\bullet}$ be any non-zero
element. The cyclic algebra 
$E=E(L/K,\sigma,a_0)$ is the $m$-dimensional $L$-vector space
\begin{equation*}
E:=L.1\oplus L.e\oplus L.e^2\oplus \cdots \oplus L.e^{m-1}
\end{equation*}
with multiplication defined by:
$e^i.e^j:=e^{i+j}$, $e^m:= a_0 \in L.1$ and $ea=\sigma(a)e$ for all
$a\in L$. 
The center of $E$ is $K$ and the dimension of $E$ over $K$ is $m^2$.

Let us now describe an important class of cyclic algebras,
namely the quaternion algebras.\index{Algebra!quaternion|textbf}
Suppose $\alpha$ and $\beta$ are non-zero elements in $K$ and $\alpha$
is non-square. Let $L=K(\sqrt\alpha)$, let $\sigma$ be the non-trivial Galois
automorphism $\sqrt\alpha \mapsto -\sqrt\alpha$ and let
$a_0=\beta$. Then the cyclic algebra $E(L/K,\sigma,\beta)$ is the
quaternion algebra:
\begin{align*}
E= (\alpha,\beta)_K &= K\langle i,j \bigm| i^2=\alpha; j^2=\beta;
                       ij=-ji\rangle \\ 
                    &= K.1\oplus K.i\oplus K.j \oplus K.k 
\end{align*}
with multiplication defined by $i^2=\alpha$, $j^2=\beta$, $k^2=-\alpha\beta$ 
and $k=ij=-ji$.
$E$ is a non-commutative four-dimensional $K$-algebra with $K$-basis
$1$, $i$, $j$ and $ij$.

Of course, not every cyclic algebra is a division ring. In particular
a quaternion algebra $(\alpha,\beta)$ fails to be a division ring if and only
if the norm form $\langle 1,-\alpha,-\beta,\alpha\beta\rangle$ is
isotropic, i.e.~if and only if there is a non-trivial equation
$0=x_1^2-\alpha x_2^2 - \beta x_3^2 + \alpha\beta x_4^2$ with
$x_1,x_2,x_3,x_4\in K$ (see Scharlau~\cite[p76]{Scha85}). 

\subsection{Kinds of Involution}
One must distinguish two kinds of involution.
Suppose $E$ is an algebra with center a field $K$ and let $I:E\to E^o$ be an
involution. The restriction $I|_K:K\to K$ is an automorphism satisfying
$I|_K^2=\id_K$.   
\begin{definition}\index{Involution!of the first kind}\index{Involution!of the second kind|textbf}
The involution $I$ is said to be {\it of the first kind} if its
restriction to the center is the identity automorphism
$I|_K=\id_K$. Otherwise, $I$ is said to be {\it of the second kind}.
\end{definition}
\begin{remark}
If $M$ is an $\epsilon$-self-dual simple representation then the involution
\begin{equation*}
\End(M)\to\End(M);\ f\mapsto b^{-1}f^*b
\end{equation*}
depends in general on a
choice of hermitian (or skew-hermitian) form $b:M\to M^*$. However,
distinct choices
yield involutions which are conjugate and, in particular, are of the
same kind. They are not, however, isomorphic in general.

In more detail, 
if $b$,$b': M\to M^*$ are non-singular hermitian forms then we may
write $b'=bc$ for some invertible $c\in\End(M)$ so
\begin{equation*}
{b'}^{-1}f^*b'= (bc)^{-1}f^*(bc)= c^{-1}b^{-1}f^*bc
\end{equation*}
for all $f\in\End(M)$. Thus we may speak of simple self-dual
representations of the first and second kind.
\end{remark}
\subsection{Involutions of the First Kind}
The following theorem of Albert~\cite[p161]{Alb39} can also be found in Scharlau~\cite[p306,354]{Scha85}:
\begin{theorem}
\label{first_kind_classification}
i) A central simple algebra $A$ over a field $K$ admits an involution of the
first kind if and only if $A\cong A^o$. \\
ii) If $K$ is an algebraic number field and $A$ is a central simple
algebra with center $K$ such that $A\cong A^o$ then either $A=K$ or
$A$ is a quaternion algebra.
\end{theorem}
Henceforth we assume that $K$ is a number field. It follows quickly
from theorem~\ref{first_kind_classification} that
the simple algebras with involution fall into three classes as follows:
\begin{enumerate}
\item[1a)] $A=K$ is a number field with trivial involution; 
\item[1b)] $A=K.1 \oplus K.i\oplus K.j\oplus K.k$ is a quaternion
algebra with the `standard' involution $a\mapsto \overline a$ defined by
\begin{equation*}\index{Involution!standard}
\overline{i}=-i,\hspace{4mm} \overline{j}= -j,\hspace{4mm},\overline{k}=-k
\end{equation*}
which fixes precisely $K$;
\item[1c)] $A=K.1 \oplus K.i \oplus K.j\oplus K.k$ is a quaternion algebra with
`non-standard' involution\index{Involution!non-standard|textbf} $a\mapsto \widehat{a}$ defined by
\begin{equation*}
\widehat{i}= -i,\hspace{4mm} \widehat{j}=j,\hspace{4mm} \widehat{k}=k
\end{equation*}
which fixes a three-dimensional $K$-vector space.
\end{enumerate}
Since $\widehat{a}=i^{-1}\overline{a}i$ for all $a\in A$,
multiplication by $i$ transforms an hermitian form over a quaternion
algebra with non-standard involution 
$a\mapsto\widehat{a}$ into a skew-hermitian form over the
same algebra with standard involution $a\mapsto\overline{a}$ so
\begin{equation*}
W^1(A,\widehat{\ })\cong W^{-1}(A,\bar{\ }).
\end{equation*}
It is more difficult to define complete Witt invariants in case
1c) than in the other cases because the local-global principle
fails. One requires a secondary invariant which is defined if all 
local invariants vanish - see section~\ref{section:local-global} below.

\subsection{Involutions of the Second Kind}\index{Involution!of the second kind}
Suppose we are given an involution
$\bar{\ }:K\to K$ which fixes an index~$2$ subfield $k\subset
K$ say. We distinguish just two classes of algebras with involutions of
the second kind. Either 
\begin{enumerate}
\item[2a)] $A=K$ or
\item[2b)] $A$ is a cyclic algebra over $K$.
\end{enumerate}
There is no need to draw any distinction between hermitian and
skew-hermitian forms, for if $K=k(\sqrt\alpha)$ then 
$\overline{\sqrt\alpha}=-\sqrt\alpha$ and the transformation
$\phi\mapsto \phi\alpha$ converts hermitian forms into skew-hermitian
forms (and vice versa) so $\displaystyle{W^1(K)\cong W^{-1}(K)}$
\section{Witt Invariants} 
\label{section:number_theorm_Witt_invariants}
The simplest Witt invariant of a symmetric or hermitian form is the
rank $m$ modulo~$2$. Of course, if one or more signature invariants
are defined then this rank invariant can often be deduced from the
signatures.  
An exhaustive multi-signature invariant $\sigma$ has already been
defined in earlier chapters. We will shortly recall the definitions of
two further local invariants - the discriminant $\Delta$ and the
Hasse-Witt invariant $c$. In the case of quaternionic algebras
with non-standard involution, one requires an extra relative invariant
such as the invariant $\theta$ which was introduced by D.Lewis. 

Among many possible sources, the reader is referred
to~Landherr~\cite{Lan36}, Lewis~\cite{Lew82}, Milnor and
Husemoller~\cite{MilHus73}, O'Meara~\cite{O'M63} and
Scharlau~\cite{Scha85} for a proof of the following theorem:
\begin{theorem}
\label{table_of_invariants}
The following table indicates a sufficient set of invariants
to distinguish Witt classes of forms over each of the five
classes of finite-dimensional division $\Q$-algebra:
\begin{equation*}
\begin{array}{|c|l|l|l|}
\hline
{\it Class} & {\it Division Algebra} & {\it Involution} & {\it Invariants} \\
\hline
1a & \mbox{Commutative} & \mbox{Trivial} & \mbox{$m\; (2)$,
$\sigma$, $\Delta$, $c$} \\
\hline
1b & \mbox{Quaternionic} & \mbox{Standard} & \mbox{$m\; (2)$,
$\sigma$}  \\
\hline
1c & \mbox{Quaternionic} & \mbox{Non-standard} & 
\mbox{$m\; (2)$, $\sigma$, $\Delta$, $\theta$} \\
\hline
2a & \mbox{Commutative} & \mbox{Non-trivial} & \mbox{$m\; 
(2)$, $\sigma$, $\Delta$} \\
\hline
2b & \mbox{Cyclic} & \mbox{Second Kind} & \mbox{$m\; (2)$, $\sigma$,
$\Delta$} \\
\hline
\end{array}
\end{equation*}
\end{theorem}
\noindent Note that in the case $1c$ the invariant $\theta$ is defined only if
all the other invariants vanish.
\begin{corollary}
\label{eight_torsion}
Let $R$ be any ring with involution. The subgroup
$8W^\epsilon(R\dash \Q)\subset W^\epsilon(R\dash \Q)$ is torsion-free.
\end{corollary}
\begin{proof}
By corollary~\ref{general_devissage_and_ME} $W^\epsilon(R\dash \Q)$ is
isomorphic to a direct sum of Witt groups of division algebras
$W^1(\End(M))$. It follows from theorem~\ref{table_of_invariants},
that the exponent of each summand $W^1(\End(M))$ divides eight (compare
caveats~\ref{addition_for_hermitian_fields}
and~\ref{addition_in_witt_ring_of_number_field} below.)
\end{proof}
\subsection{Discriminant $\Delta$}\index{Discriminant|textbf}
\label{define_discriminant}
\paragraph*{1a}
The {\it determinant} of a non-singular symmetric form $\phi:K^m\to
(K^m)^*$ over any field $K$ (with trivial involution) is the determinant
of the $m\times m$ matrix $A$ representing $\phi$ with respect to some
choice of basis. Changing the basis one has
$\det(PAP^t)=\det(A)\det(P)^2$ so $\det(\phi)$ is a well-defined
element of $\frac{K^{\,\bullet}}{\left(K^{\,\bullet}\right)^2}$ where
$\left(K^{\,\bullet}\right)^2=\{x^2\mid x\in K^{\,\bullet}\}$.

The {\it discriminant} of $\phi$ is, by definition,
\begin{equation*}
\Delta(\phi):=(-1)^{\frac{m(m-1)}{2}}\det(\phi) \in 
\frac{K^{\,\bullet}}{\left(K^{\,\bullet}\right)^2}.
\end{equation*}
The discriminant vanishes on hyperbolic forms and is therefore a
well-defined Witt invariant:
\begin{equation*}
\Delta: W^1(K)\to \frac{K^{\,\bullet}}{\left(K^{\,\bullet}\right)^2}~.
\end{equation*}
\begin{caveat}
\label{addition_for_hermitian_fields}
$\Delta$ is not a group homomorphism. On the other hand, the rank
modulo~$2$ and the discriminant together with a group homomorphism
\begin{equation*}
W^1(K)\to
\frac{\Z}{2\Z}\times\frac{K^{\,\bullet}}{\left(K^{\,\bullet}\right)^2}.  
\end{equation*}
where the group operation on the right-hand side is
`$\frac{\Z}{2\Z}$-graded':
\begin{align*}
(0,a)+(0,b) &=(0,ab) \\
(0,a)+(1,b) &=(1,ab) \\
(1,a)+(1,b) &=(0,-ab)~.
\end{align*}
\end{caveat}

\paragraph*{1b} Suppose next that $K$ is a field with a non-trivial involution
$a\mapsto\overline{a}$ which fixes a subfield $k$ of index~$2$.
The discriminant of a hermitian form $\phi$ over $K$
is defined as above, but the value group is slightly different: 
\begin{align*}
\Delta: W^1(K) &\to
\frac{k^{\,\bullet}}{K^{\,\bullet}\overline{K^{\,\bullet}}} \\
[\phi] &\mapsto (-1)^{\frac{m(m-1)}{2}}\det(\phi).
\end{align*}
To explain the notation,
$K^{\,\bullet}\overline{K^{\,\bullet}}=\{x\overline x\mid x\in
K^{\,\bullet}\}$.

\paragraph*{1c and 2b} Finally, if $E$ is a division ring of dimension $d^2$ over its center
$K$ and $\phi:E^m\to (E^m)^*$ is an hermitian form then the
determinant of $\phi$ is by definition the reduced norm of a matrix
representing $\phi$. 
It can be computed as follows:
Choose a Galois extension $L$ of
$K$ which splits $E$, i.e.~such that $L\otimes_K E\cong M_d(L)$. Then
$\phi$ induces a form $\phi_{L\otimes E}$ which is represented by a
$dn\times dn$ matrix over $L$ and whose determinant is $\det(\phi)$. 
To obtain a Witt invariant one defines the discriminant
\begin{equation*}
\Delta(\phi)=(-1)^{\frac{m(m-1)}{2}d}\det(\phi)
\end{equation*}
which takes values in
$\frac{K^{\,\bullet}}{\left(K^{\,\bullet}\right)^2}$ if the involution
is of the first kind or in
$\frac{k^{\,\bullet}}{K^{\,\bullet}\overline{K^{\,\bullet}}}$ if the
involution is of the second kind. Note that if $d$ is even, such as in the case
1c, then we have $\Delta(\phi)=\det(\phi)$. 
\subsection{Hasse-Witt Invariant $c$}\index{Hasse-Witt invariant|textbf}
Suppose $K$ is an algebraic number field with trivial involution and
let $\varPsi$ denote the set of primes.
For each prime $\mf{p}\in\varPsi$ and each pair of units $a,b\in
K^{\,\bullet}$ there is defined a Hilbert symbol 
\begin{equation*} 
(a,b)_{\mf{p}} = 
\begin{cases} 
\phantom{-}1 &\mbox{if there exist $x,y\in K_{\mf{p}}$ such that $ax^2+by^2=1$}
\\ 
-1 &\mbox{otherwise}
\end{cases}
\end{equation*}
Equivalently, $(a,b)_{\mf{p}}=1$ if and only if the quaternion algebra 
$(a,b)_{K_{\mf{p}}}$ is isomorphic to the matrix algebra
$M_2(K_{\mf{p}})$. 
The number of primes at which $(a,b)_{\mf{p}}=-1$ is finite (and,
according to Hilbert reciprocity, even) so, considering all primes together,
we may write 
\begin{equation*}
(a,b)\in \bigoplus_\varPsi\{+1,-1\}.
\end{equation*}

Given a non-singular symmetric form $\phi=\langle
a_1,a_2,\cdots,a_m\rangle$ over $K$ 
let
\begin{equation*}
s(\phi)=\prod_{i<j}(a_i,a_j)~\in~\bigoplus_\varPsi\{+1,-1\}.
\end{equation*}
It turns out that $s$ is independent of the diagonalization.
Some adjustment is needed to obtain a well-defined Witt
invariant $c$~\cite[p81]{Scha85}: 
\begin{equation*}
c(\phi) :=
\begin{cases}
s(\phi) &\mbox{if $m\cong 1,2\pmod 8$} \\
(-1,-\det(\phi))s(\phi) &\mbox{if $m\cong 3,4\pmod 8$} \\
(-1,-1)s(\phi) &\mbox{if $m\cong 5,6\pmod 8$} \\
(-1,\det(\phi))s(\phi) &\mbox{if $m\cong 7,8\pmod 8$} \\
\end{cases}.
\end{equation*}
\begin{caveat}
\label{addition_in_witt_ring_of_number_field}
The Hasse-Witt invariant $\displaystyle{c:W^1(K)\to \bigoplus_\varPsi\{+1,-1\}}$ is not a 
homomorphism. However, the rank modulo $2$, the discriminant and the
Hasse-Witt invariant together give a homomorphism
\begin{equation*}
W^1(K)\to \frac{\Z}{2\Z}\times \frac{K^{\,\bullet}}{\{a^2\mid a\in
K^{\,\bullet}\}} \times \left(\bigoplus_\varPsi\{+1,-1\}\right)
\end{equation*}
where the group law on the right hand side is given by
\begin{align*}
(0,d,c)+(0,d',c') &=(0,dd',(d,d')cc') \\
(0,d,c)+(1,d',c') &=(1,dd',(d,-d')cc') \\
(1,d,c)+(0,d',c') &=(1,dd',(-d,d')cc') \\
(1,d,c)+(1,d',c') &=(0,-dd',(d,d')cc').
\end{align*}
\end{caveat}
\subsection{Local-Global Principle}
\label{section:local-global}
Naturality of devissage and hermitian Morita
equivalence (remarks~\ref{naturality_of_devissage}
and~\ref{naturality_of_hermitian_Morita_equivalence}) give rise to a
commutative diagram: 
\begin{equation*}
\xymatrix{
 W^\epsilon(R\dash \Q) \ar[rr] \ar@{<->}[d] & & {\prod_{p}
W^\epsilon(R\dash \Q_p)}  \\
{\bigoplus_M W^\epsilon_M(R\dash \Q)} \ar[r] \ar@{<->}[d]_{\cong} 
& {\bigoplus_M\prod_p W^\epsilon_{\Q_p\otimes M}(R\dash \Q_p)}
\ar@{>->}[r] \ar@{<->}[d]_{\cong}  
& {\prod_p \bigoplus_M W^\epsilon_{\Q_p\otimes_\Q M} (R\dash \Q_p)} 
\ar[u] \ar@{<->}[d]^{\cong} \\
{\bigoplus_M W^1(E_M)} \ar[r] 
& {\bigoplus_M\prod_p W^1(\Q_p\otimes_\Q E_M)} \ar@{>->}[r] 
& {\prod_p\bigoplus_M W^1(\Q_p\otimes_\Q E_M)}
}
\end{equation*}
where $E_M=\End_{(R\dash \Q)}M$.
All the sums are indexed by the isomorphism classes of $\epsilon$-self-dual
simple representations~$M$ of~$R$ over $\Q$. 

When $E_M$ is not a quaternion algebra with non-standard involution
(i.e.~not in class 1c above), a local global principle applies.
In other words, the natural map
\begin{equation}
\label{local_global_map}
W^1(E_M)\rightarrow \prod_p W^1(\Q_p\otimes_\Q E_M)
\end{equation}
is injective. 
\subsection{Lewis $\theta$-invariant}\index{Lewis $\theta$-invariant|textbf}
Suppose $E=K(a,b)=K\langle i,j \bigm| i^2=a; j^2=b;
                       ij=-ji\rangle$ 
is a quaternion division algebra with non-standard involution $\widehat{i}=-i$,
$\widehat{j}=j$. 
The rank modulo $2$, the signatures and the discriminant are together
a complete set of local Witt invariants; that is to say, if 
two Witt classes cannot be distinguished by these three invariants then they
have the same image under the natural map~(\ref{local_global_map}).

One further relative invariant is needed to distinguish elements in the kernel
of~(\ref{local_global_map}). The first such invariant was constructed by
Bartels~\cite{Bar75,Bar76} using Galois cohomology. We shall
describe instead the more elementary invariant $\theta$ which is due to
Lewis~\cite{Lew82}. 

Let $L=K(\sqrt a)$ and let~$W^1(L)$ and~$W^1(L,\bar{\ })$ denote the
Witt groups of~$L$ with trivial involution and involution~$\sqrt
a\mapsto -\sqrt a$ fixing~$K$ respectively.
The definition of~$\theta$ involves the following commutative
diagram:
\begin{equation*}
\xymatrix{
0 \ar[r] & W^1(E,\bar{\ }) \ar[r] \ar[d] & W^1(L,\bar{\ }) \ar[r]^\iota \ar[d]
& W^1(E,\widehat{\ }) \ar[r] \ar[d]^\Lambda & W^1(L) \ar[d]^{\Lambda'}
\cdots \\
0 \ar[r] & {\prod W^1(E_{\mf{p}},\bar{\ })} \ar[r] & {\prod W^1(L_{\mf{p}},\bar{\ })}
\ar[r] & {\prod W^1(E_{\mf{p}},\widehat{\ })} \ar[r] & {\prod
W^1(L_{\mf{p}})} \cdots }
\end{equation*}
in which the horizontal sequences are exact
(Lewis~\cite{Lew79,Lew82,Lew82''}) and the products are indexed by the
primes $\mf{p}$ of $K$. The map $\iota$ is induced by the
inclusion of $L=K(\sqrt{a})\hookrightarrow E;\quad \sqrt{a}\mapsto
i$. Exactness of the upper row and the injectivity of~$\Lambda'$
is that $\Ker(\Lambda)\subset
\Im(\iota)$. Given any $\phi\in \Ker(\Lambda)$
one can choose $\psi\in W^1(L,\bar{\ })$ such that~$\iota(\psi)=\phi$
and let~$d$ be the discriminant~
\begin{equation*}
d=\Delta(\psi)\in
\frac{K^{\,\bullet}}{L^{\,\bullet}\overline{L^{\,\bullet}}}~.
\end{equation*} 
There is now a well-defined injection 
\begin{align*}
\theta : \Ker(\Lambda) &\to \{+1,-1\}^S \slash\sim \\
[\phi] &\mapsto \{(d,a)_{\mf{p}}\}_{\mf{p}\in S}
\end{align*}
where~$S$ is the set of primes at which~$E_{\mf{p}}$ is a division
algebra and where $\sim$ identifies each element $(\epsilon_1, \cdots,
\epsilon_{|S|})$ with its antipode $(-\epsilon_1, \cdots, -\epsilon_{|S|})$.

\section{Localization Exact Sequence}\index{Localization exact sequence}
\label{section:localization_exact_sequence}
Although we shall not attempt to extricate the subgroup
$W^\epsilon(P_\mu\dash \Z)$ from $W^\epsilon(P_\mu\dash \Q)$ in any detail,
a first step in that direction is the following localization exact
sequence which makes sense for any ring $R$ with involution:
\begin{equation}
\label{localization_sequence}
0\to W^\epsilon(R\dash \Z)\to W^\epsilon_\Z(R\dash \Q) \to
W^\epsilon(R\dash \Q/\Z)\to\cdots
\end{equation}
In this section we define and analyse the groups $W^\epsilon_\Z(R\dash \Q)$ and
$W^\epsilon(R\dash \Q/\Z)$. We omit to prove that~(\ref{localization_sequence})
is exact since the argument is quite standard - 
e.g.~Stoltzfus~\cite[pp16-19]{Sto77} or~Neumann~\cite[Theorem 6.5]{Neu77}). 

\begin{definition}
If $R$ is any ring with involution,
let $(R\dash \Q)_\Z\proj$\index{$(R\dash \Q)_\Z\proj$} denote the full
subcategory of $(R\dash \Q)\proj$
containing precisely those rational representations $(M,\rho)$ which
are induced up from integral representations,
i.e.~$(M,\rho)\cong\Q\otimes_\Z(M_0,\rho_0)$.
\end{definition}
\begin{definition}\index{$W^{\epsilon}_\Z(R\dash \Q)$}
Let $W^{\epsilon}_\Z(R\dash \Q)$ denote the Witt group of $(R\dash
\Q)_\Z\proj$.
\end{definition}
\begin{lemma}
\label{lemma:Zsummand}
Suppose we are given a short exact sequence in $(R\dash \Q)\proj$
\begin{equation*}
0\to M\to M'\to M''\to0.
\end{equation*}
If $M'\in (R\dash \Q)_\Z\proj$ then $M$ and $M''$ are also in $(R\dash \Q)_\Z\proj$.
It follows, in particular, that $(R\dash \Q)_\Z\proj$ is an abelian category.
\end{lemma}
\begin{proof}
If $M'\cong \Q\otimes_\Z M'_0$ then 
\begin{equation*}
M\cong\Q\otimes_\Z(M'_0\cap M)\quad\mbox{and}\quad
M''\cong \Q\otimes_\Z \frac{M'_0}{M\cap M'_0}. \qedhere
\end{equation*} 
\end{proof}

\begin{proposition}
There is an isomorphism
\begin{equation}
\label{eqn:algintdevissage}
W^\epsilon_\Z(R\dash \Q)\cong\bigoplus_M W^1(\End(M))
\end{equation}
with one direct summand for each simple $\epsilon$-self-dual
$M\in(R\dash \Q)_\Z\proj$. In particular $W^\epsilon_\Z(R\dash \Q)$ is a direct
summand of $W^\epsilon_\Z(R\dash \Q)$.
\end{proposition}
\begin{proof}
The first sentence is a consequence of 
theorems~\ref{hermitian_devissage}
and~\ref{hermitian_Morita_equivalence}. 
It follows from lemma~\ref{lemma:Zsummand} that simple objects in
$(R\dash \Q)_\Z\proj$ are simple in $(R\dash \Q)\proj$ which implies the second
sentence of the proposition.
\end{proof}

Proposition~\ref{character_in_Z} showed that a semisimple rational
representation $M\in(R\dash \Q)\proj$ lies in $(R\dash \Q)_\Z\proj$ if and only if the
character $\chi_M$ takes values in~$\Z$. Let us now give a criterion
which applies also to non-semisimple representations.
\begin{lemma}
\label{lemma:invariantlattice}
Suppose $(M,\rho)\in (R\dash \Q)\proj$. Then the following 
are equivalent: 
\begin{enumerate}
\item $M\in(R\dash \Q)_\Z\proj$.
\item $M$ is algebraically integral. In other words there exists
a rational representation $M'\in (R\dash \Q)\proj$ and an integral
representation $N_0\in(R\dash \Z)\proj$ such that 
$M\oplus M'\cong\Q\otimes_\Z M_0$.

\item The subring $\rho(R)\subset\End_\Q M$ is finitely generated as a
$\Z$-module.
\end{enumerate}
\end{lemma}
\begin{proof}
$1\Rightarrow 2$: Immediate. \\
$2\Rightarrow 1$: Suppose $M\oplus M'\simeq \Q\otimes N_0$
for some $M'\in (R\dash \Q)\proj$ and some $N_0\in (R\dash \Z)\proj$.
The intersection $M_0=N_0\cap M$ is invariant under the action of $R$
and satisfies $M\cong \Q\otimes_\Z M_0$. \\
$1\Rightarrow 3$: The representation $\rho$ factors through a
canonical inclusion
\begin{equation*}
\End_\Z M_0 \hookrightarrow \End_\Q M.
\end{equation*}
Now $\rho(R)$ is contained in the finitely generated free $\Z$-module
$\End_\Z M_0$. Since $\Z$ is Noetherian, $\rho(R)$ is
finitely generated. \\ 
$3\Rightarrow 1$: If $\rho(R)$ is finitely generated then $\rho(R)x$
is a finitely generated sub-$\Z$-module for every $x\in M$. If
$x_1,\cdots,x_m$ is a basis for $M$ over $\Q$ then we may set
$M_0=\sum_i \rho(R)x_i$.
\end{proof}

We turn next to the group $W^\epsilon(R\dash \Q/\Z)$.
\begin{definition}
We denote by $(R\dash \Q/\Z)\proj$\index{$(R\dash \Q/\Z)\proj$} the abelian category of representations
$(M,\rho)$ where $M$ is a finite abelian group and $\rho:R\to\End_\Z
M$ is a ring homomorphism. 
\end{definition}
The dual module $M^\wedge=\Hom_\Z(M,\Q/\Z)$
admits the usual $R$-action 
\begin{equation*}
\rho^\wedge(r)(\theta)= (m\mapsto \theta(\overline{r}m)) \quad \mbox{for
all $\theta\in M^\wedge$, $r\in R$ and $m\in M$}.
\end{equation*}
so we have defined a duality functor on~$(R\dash \Q/\Z)\proj$.
\begin{definition}
Let $W^\epsilon(R\dash \Q/\Z)$\index{$W^\epsilon(R\dash \Q/\Z)$}
denote the Witt group of the hermitian category $(R\dash
\Q/\Z)\proj$.
\end{definition}
The machinery described in chapters~\ref{chapter:morita_equivalence}
and~\ref{chapter:devissage} is general enough to compute this Witt group
since $(R\dash \Q/\Z)\proj$ is an abelian category. By
theorems~\ref{hermitian_devissage} and~\ref{hermitian_Morita_equivalence}
$W^\epsilon(R\dash \Q/\Z)$ is a direct sum of Witt groups of finite fields. These
are easy to compute since the group of units of a finite field is cyclic:
\begin{lemma}
If $K$ is a finite field with trivial involution then
\begin{align*}
W^1(K)&\cong\begin{cases}
\frac{\Z}{2\Z}\oplus\frac{\Z}{2\Z} &\mbox{if $|K| \equiv 1 \pmod{4}$} \\
\frac{\Z}{4\Z} &\mbox{if $|K| \equiv 3 \pmod{4}$} \\
\frac{\Z}{2\Z} &\mbox{if ${\rm char}(K) = 2$}. 
\end{cases}
\intertext{If the involution is non-trivial then}
W^1(K)&\cong\frac{\Z}{2\Z}.
\end{align*}
\end{lemma}
\begin{proof}
See for example Scharlau~\cite[p40]{Scha85} and Milnor and
Husemoller \cite[p117]{MilHus73}.
\end{proof}

Suppose $p$ is an odd prime. It is easy to construct examples of
self-dual objects in $(P_\mu\dash \Q/\Z)\proj$ whose endomorphism ring has $p$
elements. For example, let $s$ act on $\frac{\Z}{p\Z}$ as
multiplication by $\frac{p+1}{2}$. We obtain:
 \begin{proposition}
There is an isomorphism
\begin{equation*}
W^\epsilon(P_\mu\dash \Q/\Z)~\cong~\bigoplus_M W^1(\End(M))~\cong~\left(\frac{\Z}{4\Z}\right)^{\oplus\infty}\oplus\left(\frac{\Z}{2\Z}\right)^{\oplus\infty}.
\end{equation*}
There is one summand $W^1(\End(M))$ for each isomorphism class of
$\epsilon$-self-dual simple representations $M\in(P_\mu\dash \Q/\Z)\proj$.
\end{proposition}
\chapter{All Division Algebras Occur}
\label{chapter:endrings}
We prove here that every finite-dimensional division algebra with
involution is the endomorphism ring over $(P_\mu\dash \Q)$ of some simple
$\epsilon$-self-dual integral representation of $P_\mu$. Consequently,
all five of the classes of division algebras highlighted in
chapter~\ref{chapter:complete_invariants} arise in the computation of  
boundary link cobordism (whereas only the class 2a of number fields
with non-trivial involution occurs when $\mu=1$). We conclude the
chapter with a proof of theorem~\ref{the_answer_up_to_isomorphism}.

Recall that the quiver path ring $P_\mu$ has presentation
\begin{equation*}
\Z\left\langle s,\pi_1,\cdots,\pi_\mu \biggm| \sum_{i=1}^\mu \pi_i=1, \pi_i^2=\pi_i, \pi_i\pi_j=0 \text{\
for \mbox{$1\leq i,j\leq \mu$}} \right\rangle.
\end{equation*}
and involution $s\mapsto 1-s$, $\pi_i\mapsto\pi_i$ for $1\leq i\leq\mu$.
The endomorphism ring $\End_{(P_\mu\dash \Q)}M$ of a rational
representation $(M,\rho)$ is
\begin{equation}
\label{endrels}
\left\{(\alpha_1,\cdots,\alpha_\mu)\in
\bigoplus_{i=1}^\mu \End_\Q(\pi_iM) \biggm| \ \alpha_i \rho(s)_{ij} = \rho(s)_{ij}\alpha_j~\mbox{for $1\leq i,j\leq\mu$}\right\}
\end{equation}
where $\rho(s)_{ij}$ is the composite $\pi_jM\hookrightarrow M
 \xrightarrow{\rho(s)} M\twoheadrightarrow\pi_i M$.

\begin{proposition}
\label{OCCURRENCE_OF_ALGEBRA}
Suppose~$E$ is a finite-dimensional division algebra over~$\Q$
and let ~$\mu\geq 2$. There exists an integral representation $M_0$
of~$P_\mu$, finitely generated and free over $\Z$, such that $M=\Q\otimes_\Z
M_0$ is simple and $\End_{(P_\mu\dash \Q)}(M)$ is isomorphic to $E$.
\end{proposition}

Recall that the opposite ring $E^o$ is identical to $E$ as an additive
group but multiplication is reversed. 
\begin{proposition}
\label{OCCURRENCE_OF_ALGEBRA_WITH_INVOLUTION}
Suppose $E$ is a finite-dimensional division algebra over~$\Q$ with
involution $I:E\rightarrow E^o$. Let
$\epsilon=+1$ or $-1$ and let $\mu\geq 2$. There exists an $\epsilon$-self-dual
integral representation $M_0$ of $P_\mu$ and an $\epsilon$-hermitian form
$b_0:M_0\rightarrow M_0^*$ such that \\ \noindent
i) $M_0$ is finitely generated and free over $\Z$. \\ \noindent
ii) $M=\Q\otimes_\Z M_0$ is simple. \\ \noindent
iii) There is an isomorphism 
\begin{equation*}
(\End_{(P_\mu\dash \Q)}M\ ,\ \beta\mapsto b^{-1}\beta^*b)~\cong~(E,I)
\end{equation*}
of algebras with involution. Here, $b:M\to M^*$ is
the isomorphism induced by $b_0$.
\end{proposition}

The idea in the proof of proposition~\ref{OCCURRENCE_OF_ALGEBRA}
will be to express $E$ as an endomorphism ring via the identity
\begin{equation}
\label{every_ring_is_endo}
E\cong\End_{E^o}E^o
\end{equation}
and let $P_\mu$ mimic the action of $E^o$ on $E^o$.
Our proof of proposition~\ref{OCCURRENCE_OF_ALGEBRA_WITH_INVOLUTION}
uses, in addition, the observation that for every
involution
\begin{equation*}
I:E=\End_{E^o}E^o\to \End_{E^o}{E^o}^*= E^o
\end{equation*}
there exists $\delta:{E^o}^*\to {E^o}^*$ such that
$I(x)=\delta^{-1}x^*\delta$ for all $x\in \End_{E^o}E^o$. Such
$\delta$ exists because left vector spaces of equal
dimension over $E^o$ are isomorphic.
\section{Proof of Proposition~\ref{OCCURRENCE_OF_ALGEBRA}}
Let $l:E\to \End_\Q E$ and $r:E^o\to \End_\Q E$ denote the regular
representations `multiplication on the left' and `multiplication on
the right':
\begin{equation*}\index{$L(x)$@$l(x)$, $r(x)$}
l(x)(y)=xy\quad \mbox{and} \quad r(x)(y)=yx \quad \mbox{for all
$x,y\in E$}.
\end{equation*}
Note that $l(E)=Z_{\End_\Q E}(r(E^o))$ and $r(E^o)=Z_{\End_\Q E}(l(E))$ 
where 
\begin{equation*}
Z_{\End_\Q E}(S)=\{x\in {\End_\Q E}\ | xs= sx \ \text{for all}~
s\in S\}
\end{equation*}
denotes the commutator of a subset $S\subset\End_\Q E$.  

Let us consider first the case $\mu=2$. By theorem \ref{primthm} of
appendix~\ref{primappendix}, $E$ is generated as a $\Q$-algebra by two
elements $x_1$, $x_2$ say. We define a representation $M$ over $\Q$ as
follows: Let $\pi_1M\cong \pi_2M\cong E$ as $\Q$-vector spaces and let
$s\in P_2$ act via the matrix 
\begin{equation}
\label{saction}
\left(
 \begin{array}{c|c} 
 r(x_1) & 1 \\ \hline
 1 & r(x_2)
 \end{array}   
 \right).
\end{equation} 
By equation~(\ref{endrels}) we have
\begin{align*}
\End_{P_2}(M) &=
 \{
(\alpha,\alpha) \in {\End_\Q E}\times {\End_\Q E}\mid
\alpha r(x_i) = r(x_i)\alpha\ (i=1,2)\} \\
 &\cong \{(\alpha,\alpha)\in {\End_\Q E}\times {\End_\Q E} \mid
 \alpha\in l(E)\} \\
 &\cong E.
\end{align*}

Multiplying the generators $x_1$, $x_2$ by an integer if necessary we
can ensure that $M\cong\Q\otimes M_0$ where $M_0$ is an integral
representation, finitely generated over~$\Z$.

We must also check that $M$ is a simple representation. Indeed, if
$M'$ is a subrepresentation of $M$ then the action (\ref{saction}) of $s$
 implies that $\pi_1M'\cong \pi_2M'$ and moreover that $\pi_1M'$ is
 a sub-$E^o$-module of $E$. Thus $M'=0$ or $M$ as required.

\medskip
The construction extends easily to all $\mu\geq 2$.
Let $\pi_iM=E$ for $i=1,\cdots,\mu$ and define
$\rho(s)_{ij}:\pi_jM\to\pi_iM$ by
\begin{equation*}
\rho(s)_{ij}=\begin{cases}
r(x_i) & \mbox{if $i=j$} \\
1      & \mbox{if $i\neq j$}
\end{cases}
\end{equation*}
where $x_1$, $x_2$, $\cdots$, $x_\mu$ is any set of generators for
$E^o$ over $\Q$. This completes the proof of
proposition~\ref{OCCURRENCE_OF_ALGEBRA}.
\begin{remark}
\label{infinite_multiplicity}
This proof yields infinitely many representations
$M=\Q\otimes M_0$ with endomorphism ring $E$, no two of which are
isomorphic. Indeed, one may add any integer $a\in\Z$ to either of the
generators $x_1$,$x_2$ and continue to construct $M$ as
before. Distinct choices of $a$ yield representations with distinct
characters, so no two are isomorphic by
corollary~\ref{character_classification}. 
\end{remark}
\section{Proof of Proposition~\ref{OCCURRENCE_OF_ALGEBRA_WITH_INVOLUTION}}
\subsection{Construction of $M$ and $b$}
In outline, we begin with a representation $M$ constructed as in the previous
section which is {\it not} self-dual but which satisfies $\End_\Q
M\cong E$. The endomorphism ring of the direct sum 
$M\oplus M^*$ is $E\times E^o$ and the form 
$\left(\begin{matrix}
0 & 1 \\
\epsilon & 0
\end{matrix}\right):M\oplus M^* \rightarrow (M\oplus M^*)^*$
induces the transposition involution on $E\times E^o$:
\begin{equation}
\label{induce_transposition_involution}
\overline{\left(\begin{matrix}
\beta & 0 \\
0 & \beta'
\end{matrix}\right)}
=
\left(\begin{matrix}
0 & 1 \\
\epsilon & 0
\end{matrix}\right)^{-1}
\left(\begin{matrix}
\beta & 0 \\
0 & \beta'
\end{matrix}\right)^*
\left(\begin{matrix}
0 & 1 \\
\epsilon & 0
\end{matrix}\right)
=
\left(\begin{matrix}
\beta' & 0 \\
0 & \beta
\end{matrix}\right).
\end{equation}
We adjust the action of $P_\mu$ in such a way that the
endomorphism ring is confined precisely to the graph of the involution
\begin{equation*}
\{(\beta,I(\beta))\in E\times E^o\mid \beta\in E\}.
\end{equation*}

To explain the construction in more detail, let us
assume $\mu=2$ and let $x_1$, $x_2$ be generators for $E$ over $\Q$ as
before. We write
\begin{equation*}
\pi_1M\cong\pi_2M\cong E\oplus E^*
\end{equation*}
with $E^*=\Hom_\Q(E,\Q)$ and let $s$ act on $M$ via the matrix
\begin{equation}
\label{selfdual_saction}
\left(\begin{array}{cc|cc}
r(x_1)            & \phantom{-}0      & 1      & \phantom{-}\gamma \\
\phantom{-}0      & 1-r(x_1)^*        & \delta & -1 \\ \hline 
\phantom{-}1      & -\epsilon\gamma^* & r(x_2) & \phantom{-}0 \\
-\epsilon\delta^* & -1                & 0      & 1-r(x_2)^*
\end{array}\right).
\end{equation}
A very mild constraint on the choice of generators $x_1$ and $x_2$
will be imposed in the proof of lemma~\ref{decompose_endomorphism} and
the maps $\gamma$ and $\delta$ will be defined below. 
We define an $\epsilon$-hermitian form 
\begin{equation*}
b=\left(\begin{array}{cc|cc}
0    & 1 & 0    & 0 \\
\epsilon & 0 & 0    & 0 \\ \hline
0    & 0 & 0    & 1 \\
0    & 0 & \epsilon & 0 
\end{array}\right).
\end{equation*}
It is easy to check that $bs=(1-s)^*b$, so $b:M\to M^*$ is a homomorphism
of representations.
\subsection{Calculation of $\End_{(P_2\dash \Q)}M$}
\begin{lemma}
\label{decompose_endomorphism}
The generators $x_1$ and $x_2$ can be chosen such that 
every endomorphism $\alpha\in\End_{(P_2\dash \Q)}M$ can be written
\begin{equation*}
\alpha=
\left(
\begin{array}{cc|cc}
\beta_1 & 0        & 0       & 0 \\
0       & \beta'_1 & 0       & 0 \\ \hline
0       &   0      & \beta_2 & 0 \\
0       &   0      & 0       & \beta'_2
\end{array}
\right).
\end{equation*}
In other words
\begin{equation*}
\End_{P_2}(M)\subset \End_\Q(E)\times \End_\Q(E^*)\times\End_\Q(E)\times \End_\Q(E^*).
\end{equation*}
\end{lemma}
\begin{proof}
Certainly every endomorphism of $M$ is a direct sum of an endomorphism
of $\pi_1M$ and an endomorphism of $\pi_2M$.
We may think of the
component $\pi_1M$ as a representation of a polynomial
ring~$\Z[s_{11}]$ where $s_{11}$ acts via 
$\left(\begin{smallmatrix}
r(x_1) & 0 \\
0 & 1-r(x_1)^*
\end{smallmatrix}\right)\in\End_\Q(\pi_1M)$.

The subfield $L_1=\Q(x_1)\subset E$ generated by $x_1$ is a simple
representation of $\Z[s_{11}]$ in which the action of $s_{11}$ is
multiplication by $x_1$. Moreover, $E$ can be expressed as a direct
sum of copies of these
simple representations $E=L_1^{\oplus d}$. Dually, $E^*$ is a direct
sum of simple representations $L_1^*$ on which $s_{11}$ acts as $1-r(x_1)^*$. 

Multiplying $x_1$ by an integer if necessary, we can assume
that the endomorphism $r(x_1)\in\End_\Q E$ is {\it not} isomorphic to
$1-r(x_1)^*\in\End_\Q(E^*)$, so that $L_1$ and $L_1^*$ are not
isomorphic representations. It follows that
$\Hom_{\Z[s_{11}]}(L_1^{\oplus d},{L_1^*}^{\oplus d})=0$ whence
\begin{equation*}
\End_{\Z[s_{11}]}(\pi_1M)=\End_{\Z[s_{11}]}(E)\times\End_{\Z[s_{11}]}(E^*).
\end{equation*}
The same argument applies to $\pi_2(M)$, completing the proof of
lemma~\ref{decompose_endomorphism}. 
\end{proof}

Equation~(\ref{endrels}) imposes
relations on the components $\beta_1$,
$\beta'_1$, $\beta_2$ and $\beta'_2$ of an endomorphism
$\beta\in\End_{(P_2\dash\Q)}M$. In particular,
\begin{equation*}
\left(\begin{matrix}
\beta_1 & 0 \\
0       & \beta'_1
\end{matrix}
\right)
\left(
\begin{matrix}
1       & \phantom{-}\gamma \\
\delta  & -1
\end{matrix}
\right)
=
\left(\begin{matrix}
1      & \phantom{-}\gamma \\
\delta & -1
\end{matrix}\right)
\left(\begin{matrix}
\beta_2 & 0 \\
0 & \beta'_2
\end{matrix}\right)
\end{equation*}
so $\beta_1=\beta_2$ and $\beta'_1=\beta'_2$. If we write
$\beta=\beta_1=\beta_2$ and $\beta'=\beta'_1=\beta'_2$ then the 
relations imposed in~(\ref{endrels}) can be summarized as
follows:
\begin{align}
\label{betainE}
\beta r(x_i) &=r(x_i)\beta \ \ \mbox{for $i=1,2$}; \\
\label{beta'inE^o}
\beta' r(x_i)^* &=r(x_i)^*\beta' \ \ \mbox{for $i=1,2$};  
\end{align}
\begin{align}
\label{gammarels}
\beta\gamma &=\gamma \beta'; & \beta\gamma^* &=\gamma^*\beta'; \\
\label{deltarels}
\beta'\delta &=\delta \beta; & \beta'\delta^* &=\delta^*\beta.
\end{align}

Equation~(\ref{betainE}) is equivalent to the statement that
$\beta\in\End_\Q E$ is left multiplication by some 
element $x\in E$.  
Similarly, equation~(\ref{beta'inE^o}) is equivalent to the statement that
$\beta'\in\End_\Q(E^*)$ is dual to left multiplication by some
$x'\in E$:
\begin{equation*}
\beta=l(x); \quad\quad \beta'=l(x')^*.
\end{equation*}
The adjoint involution $\alpha\mapsto b^{-1}\alpha^*b$ is the
transposition $(l(x),l(x')^*)\mapsto (l(x'),l(x)^*)$ by
equation~(\ref{induce_transposition_involution}). We wish to choose
$\gamma$ and $\delta$ in such a way that equations~(\ref{gammarels})
and~(\ref{deltarels}) impose precisely the condition $x'=I(x)$.

If $\delta$ is invertible and satisfies $\delta l(x)=l(I(x))^*\delta$
for all $x\in E$ then the equations~(\ref{deltarels}) impose
the equivalent relations $x'=I(x)$ and $x=I(x')$ exactly as we require.
Suitable $\delta$ exists because the representations
\begin{align}
\label{isomorphic_Espaces}
E&\to \End_\Q E & E &\to \End_\Q(E^*) \\
x&\mapsto l(x)  & x &\mapsto l(I(x))^*
\end{align}
have the same dimension and are therefore isomorphic.
We can then set $\gamma=\delta^{-1}$ so that equations~(\ref{gammarels}) become
equivalent to equations~(\ref{deltarels}). 

To ensure that $(M,b)$ can be expressed as $\Q\otimes_\Z(M_0,b_0)$ we
can, if necessary, multiply $\gamma$ and $\delta$ by an integer.
\subsection{Proof that $M$ is Simple}
If $M'\subset M$ is a subrepresentation then 
$\left(\begin{matrix}
r(x_i) & 0 \\
0 & 1-r(x_i)^*
\end{matrix}\right)$
maps $\pi_iM'$ to itself for $i=1$ and for $i=2$ so in each case
$\pi_iM'=V_i\oplus W_i$ where
$V_i\subset E=L_i^{\oplus d}$ is invariant under right multiplication
by $x_i$ and
$W_i\subset E^*=(L_i^*)^{\oplus d}$ is invariant under the action of $1-r(x_i)^*$.  

The maps
\begin{align*}
\left(\begin{matrix}
1 & \phantom{-}\gamma \\
\delta & -1
\end{matrix}\right):V_2\oplus W_2 &\to V_1\oplus W_1 \\
\left(\begin{matrix}
1 & -\epsilon\gamma^* \\
-\epsilon\delta^* & -1
\end{matrix}\right):V_1\oplus W_1 &\to V_2\oplus W_2 \\
\end{align*}
imply that $V_1=V_2$ and $W_1=W_2$. Thus $V_1$ is invariant under
 multiplication by both $x_1$ and $x_2$, so $V_1=0$ or $E$. Since
 $\gamma$ and $\delta$ are isomorphisms, we have $V_1\cong W_2$ and so
 $M'=0$ or $M'=M$.
\subsection{The general case $\mu\geq2$}
If $\mu\geq 2$ one can perform a very similar construction. One
asserts that $\pi_iM=E\oplus E^*$ for $1\leq i\leq\mu$ and defines
$\rho(s):M\to M$ by 
\begin{equation*}
\rho(s)_{ij}=\begin{cases}
\left(\begin{matrix}
r(x_i) & 0 \\
0 & 1-r(x_i)^*
\end{matrix}\right) &\text{if $i=j$} \\
\left(\begin{matrix}
1 & \phantom{-}\gamma \\
\delta & -1
\end{matrix}\right) &\text{if $i\neq j$}.
\end{cases}
\end{equation*}
Here, $x_1$,$\cdots$, $x_\mu$ is a set of generators for $E$ over $\Q$
and $\gamma$ and $\delta$ are defined exactly as above.

\begin{remark}
\label{infinite_multiplicity_with_duality}
As in remark~\ref{infinite_multiplicity} above, one
may add integers to the generators $x_1$ and $x_2$ to obtain for any
$(E,I)$ an infinite family of integral representations $\{M_{i,E}\}_{i\in\N}$ 
and forms $b_{i,E}:M_{i,E}\to M_{i,E}^*$ 
such that, for each $i$, $(M_0,b_0)=(M_{i,E}, b_{i,E})$ satisfies
conditions i), ii) and iii) of
proposition~\ref{OCCURRENCE_OF_ALGEBRA_WITH_INVOLUTION}.
\end{remark}
\section{Computation of $C_{2q-1}(F_\mu)$ up to Isomorphism}
We are finally in a position to prove
theorem~\ref{the_answer_up_to_isomorphism} which states that if $\mu\geq2$ and $q>1$ then
$C_{2q-1}(F_\mu)$ is isomorphic to a countable direct sum
\begin{equation}
\label{the_isomorphism_answer_again}
C_{2q-1}(F_\mu)\cong \Z^{\oplus\infty}\oplus
\left(\frac{\Z}{2\Z}\right)^{\oplus\infty}\oplus\left(\frac{\Z}{4\Z}\right)^{\oplus\infty}
\oplus\left(\frac{\Z}{8\Z}\right)^{\oplus\infty}.
\end{equation}
We refer to Fuchs~\cite{Fuc70} for general results on infinitely
generated abelian groups.

We start, as usual, with the identification
$C_{2q-1}(F_\mu)\cong W^{(-1)^q}(P_\mu\dash \Z)$, assuming $q>2$. The case
$q=2$ of~(\ref{the_isomorphism_answer_again}) follows directly from
the other cases because $C_3(F_\mu)$ is an index $2^\mu$ subgroup of
$C_7(F_\mu)$.

Since $P_\mu$ is finitely generated, $W^\epsilon(P_\mu\dash \Z)$ is countable. 
Now the kernel of the natural map
\begin{equation}
\label{ZintoC_again}
W^\epsilon(P_\mu\dash \Z)\to W^\epsilon(P_\mu\dash \C^-)\cong\Z^{\oplus\infty}
\end{equation}
is $8$-torsion and is therefore a direct sum of cyclic
groups~\cite[Theorem 17.1]{Fuc70}. 
The image of~(\ref{ZintoC_again}) is a subgroup of a free abelian group
and hence is free abelian itself~\cite[Theorem 14.5]{Fuc70}.
Thus $C_{2q-1}(F_\mu)$ is a direct sum of cyclic groups of
orders $2$, $4$, $8$ and $\infty$. We must show that there are
infinitely many summands of each order.

By lemma~\ref{rationals_injection} and
corollary~\ref{general_devissage_and_ME} we have
\begin{equation}
\label{devissage_again}
W^\epsilon(P_\mu\dash \Z)\subset W^\epsilon(P_\mu\dash \Q) \cong \bigoplus_M
W^\epsilon_M(P_\mu\dash \Q).
\end{equation}
To each finite-dimensional division algebra $E$ with involution
let us associate a family $(M_{i,E},b_{i,E})_{i\in\N}$ of integral
representations as in remark~\ref{infinite_multiplicity_with_duality}.
Composing~(\ref{devissage_again}) with the Morita isomorphisms
\begin{equation*}
\Theta_{\Q\otimes (M_{i,E},b_{i,E})}: W^\epsilon_{\Q\otimes
M_{i,E}}(P_\mu\dash \Q)\to W^1(E)
\end{equation*}
the element $[M_{i,E}]\in W^\epsilon(P_\mu\dash \Z)$ is mapped to
$\langle1\rangle\in W^1(E)$. 

It will suffice for our purposes to consider commutative fields $E=K$
with {\it trivial} involution. By a theorem of A.Pfister, the order of
$\langle1\rangle$ in $W^1(K)$ is twice the level of $K$ (e.g.~Milnor
and~Husemoller~\cite[p75]{MilHus73} or Scharlau~\cite[pp71-73]{Scha85}):
\begin{definition}\index{Level of a field}
Let $K$ be any (commutative) field. If $-1$ is a sum of squares in $K$
then the {\it level} of $K$ is the 
smallest integer $s$ such that $-1$ is a sum of $s$ squares. If $-1$
is not a sum of squares, so that $K$ is formally real, let $s=\infty$.
\end{definition}
Pfister showed further that the level of any field is either infinity or a
power of~$2$ (cf~theorem~\ref{pfister_generalized} above). 
In the case of an algebraic number field, caveat~\ref{addition_in_witt_ring_of_number_field} confirms that
$\text{order}(\langle1\rangle)=2s$ and implies moreover that $s$ is
$1$, $2$ or $4$ or $\infty$ (compare Lam~\cite[p299]{Lam73}).

Let $T_K\subset W^1(K)$ denote the additive subgroup of $W^1(K)$
generated by $\langle1\rangle$. We have $|T_K|=2s=2$, $4$, $8$ or $\infty$.
If $|T_K|=2$, $4$ or $8$ then 
it is easy to check that $T_K$ is a pure subgroup of $W^1(K)$:
\begin{equation*}
T_K\cap mW^1(K)=mT_K ~\mbox{for all $m\in\Z$}.
\end{equation*}
It follows that $T_K$ is a direct summand of $W^1(K)$~\cite[Proposition
27.1]{Fuc70}. On the other hand, if $T_K$ is infinite then any signature
is a split surjection
$\sigma:W^1(K)\to\Z;\langle1\rangle\mapsto1$ so, again, $T_K$ is a
direct summand. The preimage 
\begin{equation*}
\bigoplus_{i,K}\Theta_{\Q\otimes (M_{i,K},b_{i,K})}^{-1}(T_K)\subset
W^\epsilon(P_\mu\dash \Z).
\end{equation*}
generated by the $[M_{i,K},b_{i,K}]$ must therefore be a direct summand.

It remains to exhibit number fields $K$ of level $1$, $2$, $4$ and $\infty$. 
\begin{itemize}
\item[$1$)] If $K=\Q(\sqrt{-1})$ then $-1$
is a square so $|T_K|=2$ and
$\left(\frac{\Z}{2\Z}\right)^{\oplus\infty}$ is a summand of
$C_{2q-1}(F_\mu)$. 
\item[$2$)] If $K=\Q(\sqrt{-3})$ then it is easy to
check that $-1$ is not square, but is a sum of two squares: 
\begin{equation*}
-1=\left(\frac{1}{2}(1+\sqrt{-3})\right)^2 +
\left(\frac{1}{2}(1-\sqrt{-3})\right)^2.
\end{equation*}
Thus $|T_K|=4$ and $\left(\frac{\Z}{4\Z}\right)^{\oplus\infty}$ is a
summand of of $C_{2q-1}(F_\mu)$.
\item[$4$)] In $K=\Q(\sqrt{-7})$, we have $-1=(\sqrt{-7})^2+2^2+1^2+1^2$.
Let us show that $-1$ is not a sum of
two squares. Let $\omega=\frac{1}{2}(1+\sqrt{-7})$ and let 
$\Q(\sqrt{-7})_\omega$ denote the completion of $\Q(\sqrt{-7})$ at the
dyadic prime ideal $(w)\vartriangleleft\Z[\omega]$. 
If $-1$ is a sum of two squares in $\Q(\sqrt{-7})_\omega$
then there exists a non-trivial equation $a^2+b^2+c^2=0$ where 
$a,b,c$ are elements of the valuation ring
$\Z[\omega]_\omega$ not all of which lie in the maximal ideal
$(\omega)\vartriangleleft\Z[\omega]_\omega$. Working modulo $\omega^2$
we obtain a non-trivial equation 
\begin{equation*}
a^2+b^2+c^2=0\quad\mbox{with}~a,b,c\in \frac{\Z[\omega]_\omega}{(\omega^2)}\cong
\frac{\Z}{4\Z}.
\end{equation*}
Only $0$ and $1$ are squares in $\frac{\Z}{4\Z}$ so we have reached a 
contradiction. Thus $|T_K|=8$ and $\left(\frac{Z}{8\Z}\right)^{\oplus\infty}$
is a summand of $C_{2q-1}(F_\mu)$.
\item[$\infty$)] The field $\Q$ is formally real; $-1$ is not
a sum of squares. Thus $T_\Q$ is infinite and $\Z^{\oplus\infty}$ is a
summand of $C_{2q-1}(F_\mu)$.
\end{itemize}
\begin{example}
Setting $E=\Q(\sqrt{-7})$ with trivial involution and
$\epsilon=(-1)^q$,
the proof of proposition~\ref{OCCURRENCE_OF_ALGEBRA_WITH_INVOLUTION}
implies that the following integral representation $(\Z^8,\rho,\phi)$
has order $8$ in $W^\epsilon(P_2\dash \Z)$: 
\begin{align*}
\rho(s) &=\left(
\begin{array}{cccc|cccc}
\phantom{-}1 & \phantom{-}0 & \phantom{-}0 & \phantom{-}0  & 1  & \phantom{-}0 & \phantom{-}0  & \phantom{-}1  \\
\phantom{-}0 & \phantom{-}1 & \phantom{-}0 & \phantom{-}0  & 0  & \phantom{-}1 & \phantom{-}1  & \phantom{-}0  \\
\phantom{-}0 & \phantom{-}0 & \phantom{-}0 & \phantom{-}0  & 0  & \phantom{-}1 & -1 & \phantom{-}0  \\
\phantom{-}0 & \phantom{-}0 & \phantom{-}0 & \phantom{-}0 & 1  & \phantom{-}0 & \phantom{-}0  & -1 \\ \hline 
\phantom{-}1 & \phantom{-}0 & \phantom{-}0 & -\epsilon     & 0 & -7 & \phantom{-}0 & \phantom{-}0 \\ 
\phantom{-}0 & \phantom{-}1 & -\epsilon    & \phantom{-}0  & 1 & \phantom{-}0 & \phantom{-}0 & \phantom{-}0 \\
\phantom{-}0 & -\epsilon    & -1           & \phantom{-}0  & 0 & \phantom{-}0 & \phantom{-}1 & -1 \\
-\epsilon    & \phantom{-}0 & \phantom{-}0 & -1            & 0 & \phantom{-}0 & \phantom{-}7 & \phantom{-}1
\end{array}
\right), \\[5mm]
\phi &=\left(
\begin{array}{cccc|cccc}
0 & 0 & 1 & 0 & 0  & 0 & 0  & 0  \\
0 & 0 & 0 & 1 & 0  & 0 & 0  & 0  \\
\epsilon & 0 & 0 & 0 & 0  & 0 & 0  & 0  \\
0 & \epsilon & 0 & 0 & 0  & 0 & 0  & 0 \\ \hline 
0 & 0 & 0 & 0 & 0 & 0 & 1 & 0 \\ 
0 & 0 & 0 & 0 & 0 & 0 & 0 & 1 \\
0 & 0 & 0 & 0 & \epsilon & 0 & 0 & 0 \\
0 & 0 & 0 & 0 & 0 & \epsilon & 0 & 0
\end{array}
\right). \\
\intertext{Any $F_\mu$-link (in dimension $2q-1>1$) which has the
  corresponding Seifert matrix} 
\lambda=\phi\rho(s) &=\left(\begin{array}{cccc|cccc}
\phantom{-}0 & \phantom{-}0 & \phantom{-}0 & \phantom{-}0 & 0 & \phantom{-}1 & -1 & \phantom{-}0 \\
\phantom{-}0 & \phantom{-}0 & \phantom{-}0 & \phantom{-}0 & 1 & \phantom{-}0 & \phantom{-}0  & -1 \\
\phantom{-}\epsilon & \phantom{-}0 & \phantom{-}0 & \phantom{-}0 & \epsilon & \phantom{-}0 & \phantom{-}0 & \phantom{-}\epsilon \\
\phantom{-}0 & \phantom{-}\epsilon & \phantom{-}0 & \phantom{-}0 & 0 & \phantom{-}\epsilon & \phantom{-}\epsilon & \phantom{-}0 \\ \hline
\phantom{-}0 & -\epsilon & -1 & \phantom{-}0 & 0 & \phantom{-}0 & \phantom{-}1 & -1 \\
-\epsilon & \phantom{-}0 & \phantom{-}0 & -1 & 0 & \phantom{-}0 & \phantom{-}7 & \phantom{-}1 \\
\phantom{-}\epsilon & \phantom{-}0 & \phantom{-}0 & -1 & 0 & -7\epsilon & \phantom{-}0 & \phantom{-}0 \\
\phantom{-}0 & \phantom{-}\epsilon & -1 & \phantom{-}0 & \epsilon & \phantom{-}0 & \phantom{-}0 & \phantom{-}0
\end{array}\right)
\end{align*}
is therefore of order $8$ in $C_{2q-1}(F_\mu)$.
\end{example}
\begin{appendix}
\def\thechapter{\Roman{chapter}}
\chapter{Primitive Element Theorems}
\label{primappendix} 
The aim of this appendix is to prove the following:
\begin{theorem}
Suppose $k$ is a commutative field of characteristic zero and $E$ is a
division $k$-algebra which is finite-dimensional over $k$.
Then there exist elements $\alpha,\beta\in E$ such that $E$ is
generated as a $k$-algebra by $\alpha$ and $\beta$.
\label{primthm}
\end{theorem}

Suppose $E\subset F$ are division rings.
If $F$ is finitely generated as a left $E$-module and $S\subset F$ is
any subset then the ring $E\langle S\rangle$ generated in $F$ by $S$ over $E$ (i.e.
the intersection of subrings of $F$ containing $E$ and $S$) is a
division ring. Indeed, any non-zero $\gamma\in E\langle S\rangle$ must
satisfy some equation 
\begin{equation*}
a_0+a_1\gamma + \cdots + a_m\gamma^m = 0
\end{equation*}
with each $a_i\in E$ and $a_0\neq0$ so we have
\begin{equation*}
\gamma^{-1}=-a_0^{-1}(a_1+a_2\gamma+\cdots+a_m\gamma^{m-1})\in
E\langle S\rangle.
\end{equation*}
To prove theorem \ref{primthm} we therefore need only show that
$\alpha$ and $\beta$ generate $E$ as a division ring over $k$.
\begin{definition}\index{$I(F/E)$@$\I(F/E)$}
 We shall use the notation
$\I(F/E)$ to denote the set of intermediate division rings
\begin{equation*}
\I(F/E)=\{G\ |\ E\subset G \subset F\}.
\end{equation*}
\end{definition}
We shall assume the following theorems:
\begin{theorem} {\rm (e.g. \cite[p243]{Lan93})}.
\label{bog_standard_primitive_element_theorem}
Suppose $L$ is a commutative field of finite-dimension over a subfield
$K$. There exists an element $\alpha\in L$ such that $K(\alpha)=L$ if
and only if $\I(L/K)$ is finite. If $L/K$ is separable then $\I(L/K)$
is finite.
\label{commprimthm}
\end{theorem}
\begin{theorem} {\rm (\cite[p110]{Coh95})}.
Suppose $E\subset F$ are division rings and $E$ is infinite. If
$\I(F/E)$ is finite then there exists $\alpha\in F$ such that
$F=E(\alpha)$.
\label{ncommprimthm}
\end{theorem}
\begin{theorem}[Centralizer Theorem] {\rm (e.g. \cite[p42]{Scha85} or
\cite[p42]{Dra83})}
Let $A$ be a simple algebra with center $K$ such that $\dim_KA<\infty$
and let $B\subset A$ be a simple sub-$K$-algebra. Then
$$\dim_KA=\dim_KB\dim_K(Z_A(B))$$
where
$$Z_A(B)=\{a\in A\ | ab=ba \ \forall b\in B\}$$
is the centralizer of $B$ in $A$.
\label{centrthm}
\end{theorem}
\begin{proof}[Proof of theorem \ref{primthm}]
Let $K$ be the center of $E$ and let $L$ be a maximal commutative
subfield of $E$, so that $Z_E(L)=L$. Now $L$ is separable over $k$ so by
theorem~\ref{commprimthm}, $L=k(\alpha)$ for some $\alpha\in L$ and the
set $\I(L/k)$ of fields intermediate between $L$ and $k$ is finite. It
follows from the centralizer theorem \ref{centrthm} that
\begin{equation*}
Z_E(Z_E(C))=C
\end{equation*}
for all $C\in\I(E/K)$ and hence that
there is a bijection
\begin{equation*}
Z_E(\functor):\I(L/K)\rightarrow\I(E/L).
\end{equation*}
Therefore $\I(E/L)$ is finite and, by theorem
\ref{ncommprimthm}, there exists $\beta\in E$ such that $E=L(\beta)$ whence
$E=k(\alpha, \beta)$ as required.
\end{proof}
\chapter{Hermitian Categories}
\label{chapter:hermitian_cat}
In this appendix we prove that a duality preserving
functor\index{Duality preserving functor}
$(F,\Phi,\eta)$  between hermitian categories is an equivalence if and
only if the underlying functor $F$ is an equivalence.  

\section{Functors}
We first recall some definitions of category theory.
Let $\cy{C}$ and $\cy{D}$ be categories and let $F:\cy{C}\to\cy{D}$ be
a functor. $F$ is said to be {\it faithful} if for all
objects $M$, $M'$ in $\cy{C}$, the induced map
\begin{equation}
\label{functor_on_morphisms}
F:\Hom_\cy{C}(M,M')\to \Hom_\cy{D}(F(M),F(M'))
\end{equation}
is injective. If this induced map~(\ref{functor_on_morphisms}) is
surjective for all $M$ and $M'$ one says that $F$ is full.
$F$ is an equivalence of categories if there exists a functor
$F':\cy{D}\to \cy{C}$ and natural isomorphisms $\alpha:F'F\to
\id_\cy{C}$ and $\beta:FF'\to \id_\cy{D}$. 
\begin{lemma}
\label{eq_is_fully_faithful}
An equivalence of categories is faithful and full.
\end{lemma}
In fact, $F$ is an equivalence of categories if and only if $F$ is
faithful and full and every object in $\cy{D}$ is isomorphic to $F(M)$
for some object $M\in \cy{C}$ (see~Bass~\cite[p4]{Bas68}). 
\begin{lemma}
\label{relate_nat_isos}
If $F:\cy{C}\to\cy{D}$ is an equivalence of categories with inverse
$F'$ then one can choose natural isomorphisms $\alpha:F'F\to
\id_\cy{C}$ and $\beta:FF'\to \id_\cy{D}$ such that
$F'(\beta_N)=\alpha_{F(N)}$ for all objects $N\in\cy{D}$ and
$F(\alpha_M)=\beta_{F(M)}$ for all objects $M\in\cy{C}$.
\end{lemma}
\begin{proof}
Fix a natural isomorphism $\alpha:F'F\to \id_\cy{C}$. For each
$N\in\cy{D}$ define $\beta_N:FF'(N)\to N$ by
$F'(\beta_N)=\alpha_{F'(N)}$. We must check that $\beta$ is a natural
isomorphism and that $\beta_{F(M)}=F(\alpha_M)$.

Suppose $N$ and $N'$ are objects of $\cy{D}$ and 
$f\in\Hom_\cy{D}(N,N')$.  Using the naturality of $\alpha$ we have
$F'(f)F'(\beta_N)=F'(f)\alpha_{F'(N)}=\alpha_{F'(N')}F'FF'(f):F'FF'(N)\to
F'(N')$. Since $F'$ is faithful by
lemma~\ref{eq_is_fully_faithful}, $f\beta_N=\beta_{N'}FF'(f)$ so
$\beta$ is a natural transformation. 

Since $F'$ is full and $F(\beta_N)=\alpha_{F'(N)}$ has an inverse for
each $N$, $\beta_N$ is also invertible. Thus $\beta$ is a natural isomorphism.

Finally, since $\alpha$ is natural there is a commutative square
\begin{equation*}
\xymatrix{
F'FF'F(M) \ar[r]^{\alpha_{F'F(M)}} \ar[d]_{F'F(\alpha_M)} & F'F(M) \ar[d]^{\alpha_M} \\
F'F(M) \ar[r]^{\alpha_M} & M
}
\end{equation*}
so $F'F(\alpha_M)=\alpha_{F'F(M)}$ and therefore
$F'(\beta_{F(M)})=\alpha_{F'F(M)}= F'F(\alpha_M)$.
Since $F'$ is faithful, $\beta_{F(M)}=F(\alpha_M)$ as required.
\end{proof}
\section{Duality Preserving Functors}
In definition~\ref{duality_preserving_functor} the notion was
introduced of a duality preserving functor
$(F,\Phi,\eta):\cy{C}\to\cy{D}$ between two hermitian categories.
\begin{definition}\index{Duality preserving functor!composition of}
Composition of duality preserving functors is given by 
\begin{equation*}
(F',\Phi',\eta')\circ(F,\Phi,\eta)= (F'\circ F, \Phi'F'(\Phi),\eta'\eta)
\end{equation*}
with $(\Phi'F'(\Phi))_M=\Phi'_{F(M)}F'(\Phi_M)$.
\end{definition}
The identity morphism on a hermitian category $\cy{C}$ is
$(\id_\cy{C}, \{\id_M^*\}_{M\in\cy{C}},1)$, which will be abbreviated 
$\id_{\cy{C}}$. 
\begin{definition}\index{Duality preserving functor!natural
transformation between}
Suppose $(F,\Phi,\eta)$, $(F',\Phi',\eta):\cy{C}\to\cy{D}$ are duality
preserving functors. A natural transformation $\alpha:(F,\Phi,\eta)\to
(F',\Phi',\eta)$ is a natural transformation $\alpha:F\to F'$ such that 
\begin{equation}
\label{naturaltransformation_DPF}
\alpha_M^*\Phi'_M\alpha_{M^*}=\Phi_M : F(M^*)\to F(M)^*
\end{equation}
for all $M\in\cy{C}$. 
\end{definition}
\begin{lemma}
\label{pullback_DPF}
If $(F',\Phi',\eta):\cy{C}\to\cy{D}$ is a duality preserving functor and
$\alpha:F\to F'$ is a natural transformation then
equation~(\ref{naturaltransformation_DPF}) defines a pull-back $\Phi$
such that $(F,\Phi,\eta)$ is a duality preserving functor and
$\alpha:(F,\Phi,\eta)\to (F',\Phi',\eta)$ is a natural transformation.
\end{lemma}
\begin{proof}
It suffices to check that
$\Phi_M^*i_{F(M)}=\eta\Phi_{M^*}F(i_M)$. Indeed,
\begin{align*}
\Phi_M^*i_{F(M)} &= (\alpha_M^*\Phi'_M\alpha_{M^*})^*i_{F(M)} \\
		 &= \alpha^*_{M^*}(\Phi'_M)^*\alpha_M^{**}i_{F(M)} \\
		 &= \alpha^*_{M^*}(\Phi'_M)^*i_{F'(M)}\alpha_M
		 \hspace{3mm} \text{by naturality of $i$} \\
		 &= \eta \alpha^*_{M^*} \Phi'_{M^*}F'(i_M)\alpha_M \\
		 &= \eta
		 \alpha^*_{M^*}\Phi'_{M^*}\alpha_{M^{**}}F(i_M) \\
		 &= \eta \Phi_{M^*}F(i_M).\qedhere
\end{align*}
\end{proof}
\begin{definition}
\label{equivalence_of_hermitian_categories}%
\index{Hermitian!category!equivalence of}%
A duality preserving functor $(F,\Phi,\eta):\cy{C}\to \cy{D}$ is said to be an
equivalence of hermitian categories if there exists
$(F',\Phi',\eta):\cy{D}\to \cy{C}$ such that
$(F',\Phi',\eta)\circ(F,\Phi,\eta)$ is naturally isomorphic to
$\id_\cy{C}=(\id_\cy{C}, \{\id_M^*\}_{M\in\cy{C}},1)$ and
$(F,\Phi,\eta)\circ(F',\Phi',\eta)$ is naturally isomorphic to $\id_\cy{D}$.
Hermitian categories $\cy{C}$ and $\cy{D}$ are said to be
$\eta$-equivalent if there is an equivalence $(F,\Phi,\eta):\cy{C}\to\cy{D}$.
\end{definition}
\begin{proposition}
\label{cat_eq_is_hermitian_eq}
A duality preserving functor $(F,\Phi,\eta):\cy{C}\to \cy{D}$ is an
equivalence of hermitian categories if and only if $F:\cy{C}\to
\cy{D}$ is an equivalence of categories. 
\end{proposition}
\begin{proof}
The `only if' part follows directly from the definitions.
Conversely, suppose $(F,\Phi,\eta)$ is a duality preserving functor,
and $F$ is an equivalence. As in lemma~\ref{relate_nat_isos},
there exists $F':\cy{D}\to \cy{C}$ and natural
isomorphisms $\alpha:F'F\to \id_\cy{C}$ and $\beta:FF'\to \id_\cy{D}$ 
such that $F'(\beta_N)=\alpha_{F(N)}$ for all objects $N\in\cy{D}$ and
$F(\alpha_M)=\beta_{F(M)}$ for all objects $M\in\cy{C}$.

We aim to define $\Phi':F'(\functor^*)\to F'(\functor)^*$ such that
$(F',\Phi',\eta)$ is a duality preserving functor and
$\alpha:(F',\Phi',\eta)\circ(F,\Phi,\eta)\to \id_\cy{C}$ and
$\beta:(F,\Phi,\eta)\circ(F',\Phi',\eta)\to \id_\cy{D}$ are natural
isomorphisms.  $\Phi'$ must satisfy the following equations for all
$M\in\cy{C}$ and $N\in\cy{D}$:
\begin{align}
(\Phi'_N)^*i_{F'(N)} &= \eta \Phi'_{N^*} F'(i_N) : F'(N) \to
F'(N^*)^*;
\label{Phi'_is_DPF}
\\
\Phi'_{F(M)}F'(\Phi_M) &= \alpha_M^*\alpha_{M^*} : F'F(M^*)\to
F'F(M)^*;
\label{Phi'_is_inverse}
\\
\Phi_{F'(N)}F(\Phi'_N) &= \beta_N^*\beta_{N^*}:FF'(N^*)\to
FF'(N)^*.
\label{define_Phi'}
\end{align}
Since $F$ is faithful, equation~(\ref{define_Phi'}) serves as a definition of
$\Phi'_N$. \smallskip \\ \noindent
{\it Proof of equation}~(\ref{Phi'_is_DPF}): Since $F$ is faithful, it suffices to check that
\begin{equation*}
F((\Phi'_N)^*)F(i_{F'(N)})=\eta F(\Phi'_{N^*})FF'(i_N).
\end{equation*}
Composing on the left by $\Phi_{F'(N^*)}$ and applying
equation~(\ref{define_Phi'})
to the right hand side we aim to show 
\begin{equation*}
\Phi_{F'(N^*)}F((\Phi'_N)^*)F(i_{F'(N)}) =
\eta\beta_{N^*}^*\beta_{N^{**}}FF'(i_N).
\end{equation*}
Indeed,
\begin{align*}
\Phi_{F'(N^*)}F((\Phi'_N)^*)F(i_{F'(N)}) &=
F(\Phi'_N)^*\Phi_{F'(N)^*}F(i_{F'(N)}) \hspace{3mm} \text{by
naturality of $\Phi$} \\
&= \eta F(\Phi'_N)^*\Phi_{F'(N)}^*i_{FF'(N)} \hspace{3mm} \text{by
equation~(\ref{duality_functor})} \\
&= \eta (\beta_N^*\beta_{N^*})^*i_{FF'(N)} \hspace{3mm} \text{by
equation~(\ref{define_Phi'})} \\
&= \eta \beta_{N^*}^*\beta_N^{**}i_{FF'(N)} \\
&= \eta \beta_{N^*}^*i_N\beta_N \hspace{3mm} \text{by naturality of
$i$} \\
&= \eta \beta_{N^*}^*\beta_{N^{**}}FF'(i_N) \hspace{3mm} \text{by naturality of $\beta$.}
\end{align*}
\smallskip \\ \noindent
{\it Proof of equation}~(\ref{Phi'_is_inverse}):
Applying $F$ and composing on the left by $\Phi_{F'F(M)}$ it suffices
to check that
\begin{equation*}
\Phi_{F'F(M)}F(\Phi'_{F(M)})FF'(\Phi_M)=\Phi_{F'F(M)}F(\alpha_M^*)F(\alpha_{M^*}).
\end{equation*}
By equation~(\ref{define_Phi'}) the left hand side is
$\beta_{F(M)}^*\beta_{F(M)^*}FF'(\Phi_M)$ while 
\begin{align*}
\Phi_{F'F(M)}F(\alpha_M^*)F(\alpha_{M^*}) &= F(\alpha_M)^*\Phi_M
F(\alpha_{M^*}) \hspace{3mm} \text{by naturality of $\Phi$} \\
&= \beta_{F(M)}^*\Phi_M\beta_{F(M^*)}  \\
&= \beta_{F(M)}^*\beta_{F(M)^*}F'F(\Phi_M) \hspace{3mm} \text{by
naturality of $\beta$}. \qedhere
\end{align*}
\end{proof} 
\end{appendix}
\providecommand{\bysame}{\leavevmode\hbox to3em{\hrulefill}\thinspace}

\address{Dept of Mathematics, UC Riverside, \\
California 92521, USA.\\}
\email{des@sheiham.com}
\begin{theindex}

  \item $(F,\Phi,\eta)$, 41
  \item $(M,\phi)$, 38
  \item $(M,\rho)$, 35
  \item $(R\dash A)\proj$, 35
  \item $(R\dash \Q)_\Z\proj$, 98
  \item $(R\dash \Q/\Z)\proj$, 99
  \item $A\proj$, 37
  \item $A^o$, 38
  \item $A_\mu$, 12
  \item $B(n,\mu)$, 11
  \item $C(n,\mu)$, 9
  \item $\C^-$, $\C^+$, 24, 39
  \item $C_n(F_\mu)$, 12
  \item $\End(M)$, 28
  \item $\GL(\alpha)$, 61
  \item $G^{\epsilon,\mu}(A)$, 15, 20, \textbf{44}
  \item $\Gamma_n(\Z[\pi]\to\Z)$, 18, 19
  \item $H^\epsilon(A)$, 38
  \item $H^\epsilon(\cy{C})$, 40
  \item $\I(F/E)$, 109
  \item $K_0(A)$, 37
  \item $K_0(R\dash A)$, 38
  \item $K_0(\cy{C})$, 37
  \item $l(x)$, $r(x)$, 102
  \item $L_n(A)$, 17, 39
  \item $\ssdM(\cy{C},\epsilon)$, 54
  \item $\cy{M}(R)=\cy{M}(R\dash \C)$, 59
  \item $\sdM(R)$, $\cy{M}^s(R)$, $\ssdM(R)$, 59
  \item $\cy{M}(R,\alpha)$, $\sdM(R,\alpha)$, 26, 61, 62
  \item $\cy{O}$, $\cy{O}_K$, 83
  \item $P_\mu$, 25, \textbf{37}, 44
  \item $\sigma_{M,b}(L,\theta)$, 32
  \item $W^\epsilon(A)$, 39
  \item $W^\epsilon(R\dash \Q/\Z)$, 99
  \item $W^\epsilon(\cy{C})$, $W^\epsilon(R\dash A)$, 42
  \item $W^{\epsilon}_\Z(R\dash \Q)$, 98

  \indexspace

  \item Addition of knots and $F_\mu$-links, 12
  \item Admissible subobject, 41
  \item Alexander polynomial, 16, 24
  \item Algebra
    \subitem Artin, 73
    \subitem cyclic, 92
    \subitem quaternion, 28, \textbf{92}

  \indexspace

  \item Blanchfield form, 16, 21, 25
  \item Blanchfield-Duval form, 21, 23
  \item Boundary link, 11

  \indexspace

  \item Character, 26, \textbf{73--75}
    \subitem independence of, 74
  \item Cobordism
    \subitem of boundary links, 11
    \subitem of $F_\mu$-links, 12
    \subitem of links, 9
    \subitem of Seifert surfaces, 11, 12

  \indexspace

  \item Devissage, 30, \textbf{53--56}
    \subitem hermitian, 54
  \item Dimension vector, 26, \textbf{36}
  \item Discriminant, 28, \textbf{95}
  \item Duality functor, 26, \textbf{40}
  \item Duality preserving functor, \textbf{41}, 111
    \subitem composition of, 112
    \subitem natural transformation between, 112

  \indexspace

  \item $\epsilon$-hermitian form, 38, 40
  \item $\epsilon$-self-dual, 40

  \indexspace

  \item $F_\mu$-link, 11
  \item Frobenius reciprocity, 69, 71

  \indexspace

  \item Grothendieck group, 37

  \indexspace

  \item Hasse-Minskowski map, 92
  \item Hasse-Witt invariant, 28, \textbf{96}
  \item Hermitian
    \subitem category, 40
      \subsubitem equivalence of, 112
    \subitem form, 38

  \indexspace

  \item Involution, 26, 38
    \subitem non-standard, 28, \textbf{93}
    \subitem of the first kind, 93
    \subitem of the second kind, \textbf{93}, 94
    \subitem standard, 93
  \item Isotopy, 8

  \indexspace

  \item Jordan-H\"older theorem, 30, \textbf{53}

  \indexspace

  \item Knot, 8

  \indexspace

  \item Level of a field, 106
  \item Lewis $\theta$-invariant, 28, \textbf{97}
  \item Link, 8
    \subitem boundary, 11
    \subitem cobordism, 9
    \subitem $F_\mu$-, 11
    \subitem split, 13, 26
  \item Localization exact sequence, 18, 21, 33, 98
  \item Luna stratum, 26, 62

  \indexspace

  \item Metabolic, 15, 20, \textbf{42}, 43
  \item Metabolizer (=Lagrangian), 39, \textbf{42}
  \item Morita equivalence, 31, \textbf{47--51}

  \indexspace

  \item Null-cobordant, 9

  \indexspace

  \item Ordered field, 69

  \indexspace

  \item Path ring, 25, \textbf{36}
  \item Pfister's theorem, 69

  \indexspace

  \item Quiver, 25, \textbf{35}
    \subitem complete, 26

  \indexspace

  \item Real
    \subitem algebraic set, 62
    \subitem closed field, 70
    \subitem closure of a field, 70
    \subitem nullstellensatz, 63
    \subitem radical, 63
    \subitem ring (=Formally real ring), 63
    \subitem variety, 63
  \item Representation, 26, \textbf{35}
    \subitem algebraic, 82
    \subitem algebraically integral, 27, \textbf{83}
    \subitem conjugate, 82
    \subitem induced, 77
    \subitem integral, 33
    \subitem of a quiver, 36
    \subitem restriction of, 78
    \subitem self-dual, 26
    \subitem simple (=irreducible), 26, \textbf{74}
    \subitem type, 62

  \indexspace

  \item Seifert
    \subitem form, 14, 15, 20--21, \textbf{43}
    \subitem surface, 9, 11
  \item Self-dual, 40
  \item Signature, 24, 29, 32
  \item Slice, 9, 15
    \subitem boundary-, 11, 14
  \item Sublagrangian, 42
  \item Surgery, \textbf{16}, 13--21
    \subitem homology, 16, 19
    \subitem on a Seifert surface, 14--16
  \item Sylvester's theorem, 29

  \indexspace

  \item Variety of representations, 26, \textbf{59}

  \indexspace

  \item Witt group
    \subitem of a   hermitian category, 41--43
    \subitem of a ring with involution, 38
    \subitem of Blanchfield forms, 16
    \subitem of Seifert forms, 15

\end{theindex}

\begin{thebibliography}{100}

\bibitem{Alb39}
A.~A. Albert, \emph{{Structure of Algebras}}, Amer. Math. Soc. Colloq. Publ.,
  vol.~24, American Mathematical Society, New York, 1939.

\bibitem{ARS95}
M.~Auslander, I.~Reiten, and S.~O. Smal\o, \emph{{Representation Theory of
  {Artin} Algebras}}, Cambridge Stud. Adv. Math., 36, Cambridge University
  Press, 1995.

\bibitem{Bar75}
H.-J. Bartels, \emph{Invarianten hermitescher {Formen} \"uber
  {Schiefk\"orpern}}, Math. Ann. \textbf{215} (1975), 269--288.

\bibitem{Bar76}
\bysame, \emph{Zur {Klassifikation} {Schiefhermitescher} {Formen} \"uber
  {Zahlk\"orpern}}, Math. Ann. \textbf{219} (1976), no.~1, 13--19.

\bibitem{Bas68}
H.~Bass, \emph{{Algebraic {$K$-theory}}}, W. A. Benjamin Inc, New
  York-Amsterdam, 1968.

\bibitem{Ben95}
D.~J. Benson, \emph{{Representations and cohomology.{I}. {Basic} representation
  theory of finite groups and associative algebras}}, Cambridge Stud. Adv.
  Math., 30, Cambridge University Press, 1995.

\bibitem{Bla57}
R.~C. Blanchfield, \emph{Intersection theory of manifolds with operators with
  applications to knot theory}, Ann. of Math. (2) \textbf{65} (1957), 340--356.

\bibitem{BCR98}
J.~Bocknak, M.~Coste, and M.-F. Roy, \emph{{Real Algebraic Geometry}},
  Springer, Berlin, 1998.

\bibitem{Bou58}
N.~Bourbaki, \emph{{{\'El\'ements} de Math\'ematique, Book 2 Ch.8}}, Hermann,
  Paris, 1958.

\bibitem{CapSha73}
S.~E. Cappell and J.~L. Shaneson, \emph{Topological knots and knot cobordism},
  Topology \textbf{12} (1973), 33--40.

\bibitem{CapSha74}
\bysame, \emph{The codimension two placement problem, and homology equivalent
  manifolds}, Ann. of Math. (2) \textbf{99} (1974), 277--348.

\bibitem{CapSha80}
\bysame, \emph{Link cobordism}, Comment. Math. Helv. \textbf{55} (1980),
  20--49.

\bibitem{CasGor78}
A.~J. Casson and C.~McA. Gordon, \emph{On slice knots in dimension three},
  Algebraic and geometric topology, Proc. Sympos. Pure Math., no. XXXII, Part
  2, American Mathematical Society, Providence RI, 1978, pp.~39--53.

\bibitem{CasGor86}
\bysame, \emph{Cobordism of classical knots}, \`A la recherche de la topologie
  perdue, Progr. Math., 62, Birkh\"auser Boston, Boston, MA, 1986,
  pp.~181--199.

\bibitem{CocOrr90}
T.~D. Cochran and K.~E. Orr, \emph{Not all links are concordant to boundary
  links}, Bull. Amer. Math. Soc. (N.S.) \textbf{23} (1990), no.~1, 99--106.

\bibitem{CocOrr93}
\bysame, \emph{Not all links are concordant to boundary links}, Ann. of Math.
  (2) \textbf{138} (1993), 519--554.

\bibitem{COT99}
T.~D. Cochran, K.~E. Orr, and P.~Teichner, \emph{Knot concordance, {Whitney}
  towers and {$L^2$}-signatures}, e-print math.GT/9908117, 1999.

\bibitem{Coh95}
P.~M. Cohn, \emph{{Skew Fields, Theory of General Division Rings}},
  Encyclopedia Math. Appl., 57, Cambridge University Press, 1995.

\bibitem{CurRei81}
C.~W. Curtis and I.~Reiner, \emph{{Methods of Representation Theory. {Vol I}.
  {With} Applications to Finite Groups and Orders}}, John Wiley {\&} Sons, New
  York, 1981.

\bibitem{Dra83}
P.~K. Draxl, \emph{{Skew Fields}}, London Math. Soc. Lecture Note Ser., 81,
  Cambridge University Press, 1983.

\bibitem{Dub74}
D.~Dubois and G.~Efroymson, \emph{A dimension theorem for real primes}, Canad.
  J. Math. \textbf{26} (1974), no.~1, 108--114.

\bibitem{Duv86}
J.~Duval, \emph{Forme de {Blanchfield} et cobordisme d'entrelacs bords},
  Comment. Math. Helv. \textbf{61} (1986), no.~4, 617--635.

\bibitem{Far83}
M.~Farber, \emph{The classification of simple knots}, Uspekhi Mat. Nauk.
  \textbf{38} (1983), no.~5, 59--106, Russian Math. Surveys 38:5 (1983) 63-117.

\bibitem{Far91}
\bysame, \emph{Hermitian forms on link modules}, Comment. Math. Helv.
  \textbf{66} (1991), no.~2, 189--236.

\bibitem{Far92}
\bysame, \emph{Noncommutative rational functions and boundary links}, Math.Ann.
  \textbf{293} (1992), no.~3, 543--568.

\bibitem{Far92B}
\bysame, \emph{Stable-homotopy and homology invariants of boundary links},
  Trans. Amer. Math. Soc. \textbf{334} (1992), no.~1, 455--477.

\bibitem{FarVog92}
M.~Farber and P.~Vogel, \emph{The {Cohn} localization of the free group ring},
  Math. Proc. Cambridge Philos. Soc. \textbf{111} (1992), no.~3, 433--443.

\bibitem{For86}
E.~Formanek, \emph{Generating the ring of matrix invariants}, Ring theory
  (Antwerp 1985), Lecture Notes in Math., 1197, Springer, Berlin, 1986,
  pp.~73--82.

\bibitem{FoxMil66}
R.~H. Fox and J.~Milnor, \emph{Singularities of $2$-spheres in $4$-space and
  cobordism of knots}, Osaka J. Math. \textbf{3} (1966), 257--267.

\bibitem{Fuc70}
L.~Fuchs, \emph{{Infinite Abelian Groups. {Volume I}}}, Pure and Applied
  Mathematics, 36, Academic Press, New York and London, 1970.

\bibitem{Gil83}
P.~Gilmer, \emph{Slice knots in {$S^3$}}, Quart. J. Math. Oxford Ser. (2)
  \textbf{34} (1983), no.~135, 305--322.

\bibitem{Gil93}
\bysame, \emph{Classical knot and link concordance}, Comment. Math. Helv.
  \textbf{68} (1993), no.~1, 1--19.

\bibitem{Gut72}
M.~A. Guti\'errez, \emph{Boundary links and an unlinking theorem}, Trans. Amer.
  Math. Soc. \textbf{171} (1972), 491--499.

\bibitem{HamMad93}
I.~Hambleton and I.~Madsen, \emph{On the computation of the projective surgery
  obstruction groups}, $K$-theory \textbf{7} (1993), no.~6, 537--574.

\bibitem{Hil81}
J.~A. Hillman, \emph{{Alexander ideals of links}}, Lecture Notes in Math., 895,
  Springer, Berlin-New York, 1981.

\bibitem{Hil02}
J.~A. Hillman, \emph{{Algebraic Invariants of Links}}, World
Scientific Publishing, to appear.

\bibitem{Hir76}
M.~W. Hirsch, \emph{{Differential Topology}}, Grad. Texts in Math., No. 33,
  Springer, 1976.

\bibitem{Kea75'}
C.~Kearton, \emph{Blanchfield duality and simple knots}, Trans. Amer. Math.
  Soc. \textbf{202} (1975), 141--160.

\bibitem{Kea75}
\bysame, \emph{Cobordism of knots and {Blanchfield} duality}, J. London Math.
  Soc. (2) \textbf{10} (1975), no.~4, 406--408.

\bibitem{KerWeb78}
M.~Kervaire and C.~Weber, \emph{A survey of multidimensional knots}, Proc. 1977
  Plans Conf. Knot Theory, Lecture Notes in Math., 685, Springer, 1978,
  pp.~61--134.

\bibitem{Ker65}
M.~A. Kervaire, \emph{Les noeuds de dimensions sup\'erieures}, Bull. Soc. Math.
  France \textbf{93} (1965), 225--271.

\bibitem{Ker71}
\bysame, \emph{Knot cobordism in codimension two}, Manifolds--{Amsterdam} 1970,
  Lecture Notes in Math., 197, Springer, Berlin, 1971, pp.~83--105.

\bibitem{KirLiv99}
Paul Kirk and Charles Livingston, \emph{Twisted {Alexander} invariants,
  {Reidemeister} torsion and {Casson-Gordon} invariants}, Topology \textbf{38}
  (1999), no.~3, 635--661.

\bibitem{Knu91}
M.-A. Knus, \emph{{Quadratic and Hermitian Forms over Rings}}, Grundlehren der
  Mathematischen Wissenschaften, 294, Springer, Berlin, 1991.

\bibitem{Ko87}
K.~H. Ko, \emph{Seifert matrices and boundary link cobordisms}, Trans. Amer.
  Math. Soc. \textbf{299} (1987), no.~2, 657--681.

\bibitem{Ko89}
\bysame, \emph{A {Seifert}-matrix interpretation of {Cappell} and {Shaneson's}
  approach to link cobordisms}, Math. Proc. Cambridge Philos. Soc. \textbf{106}
  (1989), 531--545.

\bibitem{Lam73}
T.~Y. Lam, \emph{{The Algebraic Theory of Quadratic Forms}}, W.A.Benjamin,
  Reading, Massachusetts, 1973.

\bibitem{Lam84}
\bysame, \emph{An introduction to real algebra}, Rocky Mountain J. Math.
  \textbf{14} (1984), no.~4, 767--814.

\bibitem{Lan36}
W.~Landherr, \emph{{\"Aquivalenz} {Hermitescher} {Formen} \"uber einem
  beliebigen algebraischen {Zahlk\"orper}}, Abh. Math. Sem. Univ. Hamburg
  \textbf{11} (1936), 245--248.

\bibitem{Lan93}
S.~Lang, \emph{Algebra}, 3rd ed., Addison-Wesley, 1993.

\bibitem{LeBPro90}
L.~Le~Bruyn and C.~Procesi, \emph{Semisimple representations of quivers},
  Trans. Amer. Math. Soc. \textbf{317} (1990), no.~2, 585--598.

\bibitem{LeD88}
J.-Y. Le~Dimet, \emph{Cobordisme d'enlacements de disques}, Bull. Soc. Math.
  France \textbf{116} (1988), ii+92, M\'emoire no. 32.

\bibitem{Let00}
C.~F. Letsche, \emph{An obstruction to slicing knots using the eta invariant},
  Math. Proc. Cambridge Philos. Soc. \textbf{128} (2000), no.~2, 301--319.

\bibitem{Lev65}
J.~Levine, \emph{Unknotting spheres in codimension two}, Topology \textbf{4}
  (1965), 9--16.

\bibitem{Lev69B}
\bysame, \emph{Invariants of knot cobordism}, Invent. Math. \textbf{8} (1969),
  98--110, Addendum, 8:355.

\bibitem{Lev69}
\bysame, \emph{Knot cobordism groups in codimension two}, Comment. Math. Helv.
  \textbf{44} (1969), 229--244.

\bibitem{Lev77}
\bysame, \emph{Knot modules {I}}, Trans. Amer. Math. Soc. \textbf{229} (1977),
  1--50.

\bibitem{Lev92}
\bysame, \emph{Signature invariants of homology bordism with applications
  to links}, Knots 90 (Osaka 1990), de Gruyter, Berlin, 1992, pp.~395--406.

\bibitem{Lev94}
\bysame, \emph{Link invariants via the eta invariant}, Comment. Math. Helv.
  \textbf{69} (1994), no.~1, 82--119.

\bibitem{LevOrr00}
J.~Levine and K.~E. Orr, \emph{A survey of applications of surgery to knot and
  link theory}, Surveys on surgery theory, Vol 1, Ann. of Math. Stud., 145,
  Princeton University Press, Princeton, NJ, 2000, pp.~345--364.

\bibitem{Lew79}
D.~W. Lewis, \emph{A note on {Hermitian} and quadratic forms}, Bull. London
  Math. Soc. \textbf{11} (1979), no.~3, 265--267.

\bibitem{Lew82'}
\bysame, \emph{The isometry classification of {Hermitian} forms over division
  algebras}, Linear Algebra Appl. \textbf{43} (1982), 245--272.

\bibitem{Lew82''}
\bysame, \emph{New improved exact sequences of {Witt} groups}, J. Algebra
  \textbf{74} (1982), no.~1, 206--210.

\bibitem{Lew82}
\bysame, \emph{Quaternionic skew-{Hermitian} forms over a number field}, J.
  Algebra \textbf{74} (1982), no.~1, 232--240.

\bibitem{Lic97}
W.~B.~R. Lickorish, \emph{{An introduction to knot theory}}, Grad. Texts in
  Math., 175, Springer, New York, 1997.

\bibitem{Lun73}
D.~Luna, \emph{Slices \'etales}, Bull. Soc. Math. France M\'emoirs \textbf{33}
  (1973), 81--105.

\bibitem{Lun75}
\bysame, \emph{Adh\'erances d'orbite et invariants}, Invent. Math. \textbf{29}
  (1975), no.~3, 231--238.

\bibitem{Mil68}
J.~W. Milnor, \emph{Infinite cyclic coverings}, Conference on the Topology of
  Manifolds, Prindle, Weber {\&} Schmidt, Boston, Mass, 1968, pp.~115--133.

\bibitem{Mil69}
\bysame, \emph{On isometries of inner product spaces}, Invent. Math. \textbf{8}
  (1969), 83--97.

\bibitem{MilHus73}
J.~W. Milnor and D.~Husemoller, \emph{{Symmetric Bilinear Forms}}, Springer,
  1973.

\bibitem{Mio87}
W.~Mio, \emph{{On boundary-link cobordism}}, Math. Proc. Cambridge Philos. Soc.
  \textbf{101} (1987), 259--266.

\bibitem{Mum94}
D.~Mumford, J.~Fogarty, and F.~Kirwan, \emph{{Geometric Invariant Theory}}, 3rd
  ed., Springer, Berlin, 1994.

\bibitem{Neu77}
W.~D. Neumann, \emph{{Equivariant {Witt} Rings}}, Bonner Math. Schriften, 100,
  Universit\"at Bonn, Mathematisches Institut, 1977.

\bibitem{Nov68}
S.~P. Novikov, \emph{On manifolds with free {Abelian} fundamental groups and
  their applications}, Izvestiya Akademii Nauk SSSR. Seriya Matematicheskaya
  (1966), 207--246, Amer. Math. Soc. Transl. Ser. 2, 71 (1968), 1-42.

\bibitem{O'M63}
O.~T. O'Meara, \emph{{Introduction to quadratic forms}}, Die Grundlehren der
  mathematischen Wissenschaften, 117, Academic Press, New York; Springer,
  Berlin, 1963.

\bibitem{Par76}
W.~Pardon, \emph{Local surgery and applications to the theory of quadratic
  forms}, Bull. Amer. Math. Soc. \textbf{82} (1976), no.~1, 131--133.

\bibitem{Par77}
\bysame, \emph{Local surgery and the exact sequence of a localization for
  {Wall} groups}, Mem. Amer. Math. Soc. \textbf{12} (1977), no.~196, iv+171.

\bibitem{Pro76}
C.~Procesi, \emph{The invariant theory of $n\times n$ matrices}, Adv. Math.
  \textbf{19} (1976), 306--381.

\bibitem{Pro87}
\bysame, \emph{A formal inverse to the {Cayley}-{Hamilton} theorem}, J. Algebra
  \textbf{107} (1987), 63--74.

\bibitem{QSS79}
H.-G. Quebbemann, W.~Scharlau, and M.~Schulte, \emph{Quadratic and {Hermitian}
  forms in additive and abelian categories}, J. Algebra \textbf{59} (1979),
  no.~2, 264--289.

\bibitem{Ran80}
A.~A. Ranicki, \emph{The algebraic theory of surgery {I}. {Foundations}}, Proc.
  London Math. Soc. (3) \textbf{40} (1980), no.~1, 87--192.

\bibitem{Ran80'}
\bysame, \emph{The algebraic theory of surgery {II}. {Applications} to
  topology}, Proc. London Math. Soc. (3) \textbf{40} (1980), no.~2, 193--283.

\bibitem{Ran81}
\bysame, \emph{{Exact Sequences in the Algebraic Theory of Surgery}}, Princeton
  University Press, New Jersey; University of Tokyo Press, Tokyo, 1981.

\bibitem{Ran98}
\bysame, \emph{{High-dimensional Knot Theory}}, Springer, Berlin, 1998.

\bibitem{RRV99}
V.~Retakh, C.~Reutenauer, and A.~Vaintrob, \emph{Noncommutative rational
  functions and {Farber's} invariants of boundary links}, Differential
  topology, infinite-dimensional Lie algebras, and applications, Amer. Math.
  Soc. Transl., Ser. 2, 194, American Mathematical Society, Providence, RI,
  1999, pp.~237--246.

\bibitem{RibZal00}
L.~Ribes and P.~Zalesskii, \emph{{Profinite Groups}}, Ergeb. Math. Grenzgeb.
  (3), 40, Springer, Berlin, 2000.

\bibitem{Rol76}
D.~Rolfsen, \emph{{Knots and Links}}, Mathematics Lecture Series, 7, Publish or
  Perish, Inc, Berkeley, California, 1990, Corrected reprint of the 1976
  original.

\bibitem{Ros94}
J.~Rosenberg, \emph{{Algebraic {$K$}-theory and its applications}}, Grad. Texts
  in Math., 147, Springer, New York, 1994.

\bibitem{Sat81'}
N.~Sato, \emph{Algebraic invariants of boundary links}, Trans. Amer. Math. Soc.
  \textbf{265} (1981), no.~2, 359--374.

\bibitem{Sat81}
\bysame, \emph{Free coverings and modules of boundary links}, Trans. Amer.
  Math. Soc. \textbf{264} (1981), no.~2, 499--505.

\bibitem{Sat84}
\bysame, \emph{Alexander modules}, Proc. Amer. Math. Soc. \textbf{91} (1984),
  no.~1, 159--162.

\bibitem{Scha70}
W.~Scharlau, \emph{Induction theorems and the structure of the {Witt} group},
  Invent. Math. \textbf{11} (1970), 37--44.

\bibitem{Scha85}
\bysame, \emph{{Quadratic and {Hermitian} forms}}, Grundlehren Math. Wiss.,
  270, Springer, Berlin, 1985.

\bibitem{Ser77}
J.-P. Serre, \emph{{Linear Representations of Finite Groups}}, Grad. Texts in
  Math., 42, Springer, New York-Heidelberg, 1977.

\bibitem{Smi81}
J.~R. Smith, \emph{Complements of codimension-two submanifolds - {III} -
  cobordism theory}, Pacific J. Math. \textbf{94} (1981), no.~2, 423--484.

\bibitem{Smy66}
N.~Smythe, \emph{Boundary links}, Topology Seminar, Wisconsin, 1965 (R.H. Bing
  and R.J. Bean, eds.), Ann. of Math. Stud., 60, Princeton University Press,
  Princeton, N.J., 1966, pp.~69--72.

\bibitem{Sto77}
N.~W. Stoltzfus, \emph{{Unraveling the integral knot concordance group}}, Mem.
  Amer. Math. Soc., 192, vol. 12, Issue 1, American Mathematical Society,
  Providence, Rhode Island, 1977.

\bibitem{Tri69}
A.~G. Tristram, \emph{Some cobordism invariants for links}, Proceedings of the
  Cambridge Philosophical Society \textbf{66} (1969), 251--264.

\bibitem{Vog80}
P.~Vogel, \emph{Localisation in algebraic {$L$}-theory}, Proc. 1979 Siegen
  Topology Conf., Lecture Notes in Math., 788, Springer, 1980, pp.~482--495.

\bibitem{Vog82}
\bysame, \emph{On the obstruction group in homology surgery}, Publ. Math.
  I.H.E.S. \textbf{55} (1982), 165--206.

\bibitem{Wall70}
C.~T.~C. Wall, \emph{{Surgery on Compact Manifolds}}, 2nd ed., Math. Surveys
  Monogr., vol.~69, American Mathematical Society, 1999, (1st edition published
  1970).

\end{thebibliography}
\end{document}